\documentclass{amsart}

\usepackage{amssymb}
\usepackage{eepic}
\usepackage{color}
\usepackage{enumitem}

\allowdisplaybreaks

\numberwithin{equation}{section}

\definecolor{lgray}{gray}{0.7}

\newtheorem{theorem}{Theorem}[section]
\newtheorem{lemma}[theorem]{Lemma}
\newtheorem{proposition}[theorem]{Proposition}
\newtheorem{corollary}[theorem]{Corollary}
\newtheorem{remark}[theorem]{Remark}

\newtheorem{TheoA}{Theorem A1}
\newtheorem{TheoC}{Theorem A2}
\newtheorem{TheoB}{Theorem B1}
\newtheorem{TheoD}{Theorem B2}

\newcommand{\N}{\mathbb{N}}
\newcommand{\Z}{\mathbb{Z}}
\newcommand{\R}{\mathbb{R}}
\newcommand{\C}{\mathbb{C}}

\newcommand{\summ}{\sum\nolimits}

\newenvironment{rmk}{\begin{remark}\rm}{\end{remark}}

\def\ten{\otimes}
\def\G{\mathrm{G}}
\def\1{\mathbf{1}}
\def\H{\mathcal{H}}
\def\Q{\mathcal{Q}}
\def\NN{\mathcal{N}}
\def\M{\mathcal{M}}
\def\E{\mathsf{E}}
\def\A{\mathcal{A}}
\def\RR{\mathcal{R}}
\def\P{\mathcal{P}}
\def\T{\mathcal{S}}
\def\S{\mathcal{S}}
\def\Mn{\mathcal{N}}
\def\V{\mathrm{\mathcal{L}(G)}}

\newcommand{\dem}{\noindent {\bf Proof. }}

\newcommand{\fin}{\hspace*{\fill} $\square$ \vskip0.2cm}

\def\esssup{\mathop{\mathrm{ess \, sup}}}

\def\mean{- \hskip-10.6pt \int}

\begin{document}

\title[Algebraic Calder{\'o}n-Zygmund theory] {Algebraic Calder{\'o}n-Zygmund theory}

\author[Junge, Mei, Parcet and Xia]
{Marius Junge, Tao Mei, \\ Javier Parcet and Runlian Xia}

\maketitle

\begin{abstract}
Calder\'on-Zygmund theory has been traditionally developed on metric measure spaces satisfying additional regularity properties. In the lack of good metrics, we introduce a new approach for general measure spaces which admit a Markov semigroup satisfying purely algebraic assumptions. We shall construct an abstract form of \lq Markov metric\rq${}$ governing the Markov process and the naturally associated BMO spaces, which interpolate with the $L_p$-scale and admit endpoint inequalities for Calder\'on-Zygmund operators. Motivated by noncommutative harmonic analysis, this approach gives the first form of Calder\'on-Zygmund theory for arbitrary von Neumann algebras, but is also valid in classical settings like Riemannian manifolds with nonnegative Ricci curvature or doubling/nondoubling spaces. Other less standard commutative scenarios like fractals or abstract probability spaces are also included. Among our applications in the noncommutative setting, we improve recent results for quantum Euclidean spaces and group von Neumann algebras, respectively linked to noncommutative geometry and geometric group theory.
\end{abstract}

%\tableofcontents

\addtolength{\parskip}{+1ex}

\section*{\bf Introduction}

The analysis of linear operators associated to singular kernels is a central topic in harmonic analysis and partial differential equations. A large subfamily of these maps is under the scope of Calder\'on-Zygmund theory, which exploits the relation between metric and measure in the underlying space to provide sufficient conditions for $L_p$ boundedness. The Calder\'on-Zygmund decomposition \cite{CZ} or the H\"ormander smoothness condition for the kernel \cite{Ho} combine the notions of proximity in terms of the metric with that of smallness in terms of the measure. The doubling and/or polynomial growth conditions between metric and measure yield more general forms of the theory \cite{CW3,NTV1,NTV2,To,T2}. To the best of our knowledge, the existence of a metric in the underlying space is always assumed in the literature.  

In this paper, we introduce a new form of Calder\'on-Zygmund theory for general measure spaces admitting a Markov semigroup which only satisfies purely algebraic assumptions. This is especially interesting for measure spaces where the geometric information is poor. It includes abstract probability spaces or fractals like the Sierpinski gasket, where a Dirichlet form is defined. It is also worth mentioning that our approach recovers Calder\'on-Zygmund theory for classical spaces and provides alternative forms over them. In spite of these promising directions |very little explored here| our main motivation has been to develop a noncommutative form of Calder\'on-Zygmund theory for noncommutative measure spaces (von Neumann algebras) where the notions of quantum metric \cite{KW,Ri1,Ri2} seem inefficient.

A great effort has been done over the last years to produce partial results towards a noncommutative Calder\'on-Zygmund theory \cite{GPJP, HLMP, JMP4, MP, Pa1}. The model cases considered so far are all limited to (different) noncommutative forms of Euclidean spaces, described as follows: 

\noindent \textbf{A) Tensor products.} Let $f = (f_{\hskip-1pt jk}): \R^n \to \mathcal{B}(\ell_2)$ be a matrix-valued function and consider the tensor product extension of a standard Calder\'on-Zygmund operator acting on $f$, formally given by $$T \hskip-1pt f(x) = \int_{\R^n} k(x,y) f(y) \, dy = \Big( T \hskip-1pt f_{\hskip-1pt jk} (x) \Big)_{jk} \qquad \mbox{for} \qquad x \notin \mathrm{supp} f.$$ The $L_p$-boundedness of this map in the associated (tensor product) von Neumann algebra $\M = L_\infty(\R^n) \bar\otimes \mathcal{B}(\ell_2)$ trivially follows for $p>1$ from the vector-valued theory, due to the UMD nature of Schatten $p$-classes. On the contrary, no endpoint estimate for $p=1$ is possible using vector-valued methods. The original argument in \cite{Pa1} |also in a recent simpler form \cite{Ca}| combines noncommutative martingales with a pseudolocalization principle for classical Calder\'on-Zygmund operators. More precisely, a quantification of how much $L_2$-mass of a singular integral is concentrated around the support of the function on which it acts. This inequality has been the key tool in the recent solution of the Nazarov-Peller conjecture \cite{CPSZ}, a strengthening of the celebrated Krein conjecture \cite{PS} on operator Lipschitz functions.

\vskip3pt

\noindent \textbf{B) Crossed products.} New $L_p$ estimates for Fourier multipliers in group von Neumann algebras have recently gained considerable momentum for its connections to geometric group theory. The first H\"ormander-Mikhlin type theorem in this direction \cite{JMP4} exploited finite-dimensional cocycles of the given group $\G$ to transfer the problem to the cocycle Hilbert space $\H = \R^n$. To find sufficient regularity for $L_p$-boundedness  amounts to study Calder\'on-Zygmund operators in the crossed products $L_\infty(\R^n) \rtimes \G$ induced by the cocycle action. Nonequivariant extensions of CZOs on these von Neumann algebras were investigated in \cite{JMP4}, after identifying the right BMO space for the length function determined by the cocycle. These operators have the form $$\int_\G f_g \rtimes \lambda(g) \, d\mu(g) \, \mapsto \, \int_\G T_g(f_g) \rtimes \lambda(g) \, d\mu(g).$$ Here $\mu, \lambda$ respectively denote the Haar measure and left regular representation on the (unimodular) group $\G$, whereas $T_g = \alpha_g T \alpha_{g^{-1}}$ is a twisted form of a classical CZO $T$ on $\R^n$ by the cocycle action $\alpha$. We refer to \cite{JMP5,PRdlS} for further results. 

\vskip3pt

\noindent \textbf{C) Quantum deformations.} PDEs in matrix algebras and \lq noncommutative manifolds\rq${}$ appear naturally in theoretical physics. Pseudodifferential operators were introduced by Connes in 1980 to study a quantum form of the Atiyah-Singer index theorem over these algebras. These techniques have been underexploited over the last 30 years, due to fundamental obstructions to understand singular integral theory in this context. The core of singular integrals and pseudodifferential operator $L_p$-theory was developed in \cite{GPJP} over the archetypal algebras of noncommutative geometry. It includes quantum tori, Heisenberg-Weyl algebras and other quantum deformations of $\R^n$ of great interest in quantum field theory, string theory and quantum probability. This was the first approach to a \lq fully noncommutative\rq${}$ Calder\'on-Zygmund theory for CZOs not acting on copies of $\R^n$ as tensor or crossed product factors, but still related to Euclidean methods.

We introduce in this paper the first form of Calder\'on-Zygmund theory valid for general (semifinite) von Neumann algebras. As we explained above, the main difficulty arises from the lack of very standard geometric tools, like the existence of a nice underlying metric or the construction of suitable covering lemmas. We shall circumvent it using a very different approach based on algebraic properties of a given Markov process. Our applications cover a wide variety of scenarios which will be discussed, giving especial emphasis to noncommutative forms of Euclidean spaces and locally compact abelian groups, which are our main classical models. In the first case, we shall weaken/optimize the CZ conditions on quantum Euclidean spaces \cite{GPJP}. In the second case, LCA groups correspond to quantum groups which are both commutative and cocommutative \cite{Ta0}. We shall give CZ conditions for convolution maps over quantum groups. In the cocommutative (non necessarily commutative) context, this includes group von Neumann algebras.

\vskip-5pt

\null

\begin{center}
{\emph{Calder\'on-Zygmund extrapolation}}
\end{center}

Based on the behavior of the Hilbert transform in the real line, the main goal of Calder\'on-Zygmund theory is to establish regularity properties on the kernel of a singular integral operator, so that $L_2$-boundedness automatically extrapolates to $L_p$ boundedness for $1 < p < \infty$. A singular integral operator in a Riemannian manifold $(\mathrm{X},\mathrm{d}, \mu)$ admits the kernel representation $$T_kf(x) = \int_\mathrm{X} k(x,y) f(y) \, d\mu(y) \quad \mbox{for} \quad x \notin \mbox{supp} f.$$ Namely, $T_k$ is only assumed a priori to send test functions into distributions, so that it admits a distributional kernel in $\mathrm{X} \times \mathrm{X}$ which coincides in turn with a locally integrable function $k$ away from the diagonal $x=y$, where the kernel presents certain singularity. This already justifies the assumption $x \notin \mbox{supp} f$ in the kernel representation. The work in \cite{CZ,Ho} culminated in the following sufficient conditions on a singular integral operator in $\R^n$ for its $L_p$-boundedness:
\begin{itemize}
\item[i)] $L_2$-boundedness $$\big\| T_k: L_2(\R^n) \to L_2(\R^n) \big\| < \infty.$$

\item[ii)] H\"ormander kernel smoothness $$\esssup_{x,y \in \R^n} \int_{|x-z| >2 |x-y|} \big| k(x,z) - k(y,z) \big| + \big| k(z,x) - k(z,y) \big| \, dz <\infty.$$
\end{itemize}
Historically, this was used to prove a weak endpoint inequality in $L_1$. The same holds for Riemannian manifolds with nonnegative Ricci curvature \cite{B}. Alternatively it is simpler to use $L_2$-boundedness and the kernel smoothness condition to prove $L_\infty\to \rm BMO$ boundedness. The result then follows by well-known duality and interpolation arguments. Our strategy resembles this approach:      
\begin{itemize}
\item[P1.] Identify the appropriate BMO spaces.

\item[P2.] Prove the expected interpolation results with $L_p$ spaces.

\item[P3.] Provide conditions on CZO's which yield $L_\infty \to \mathrm{BMO}$ boundedness.
\end{itemize}
In the classical setting, we typically find $\mathrm{H}_1$/BMO spaces associated to a metric or a martingale filtration. Duong and Yan \cite{DY1,DY2} extended this theory replacing some averages over balls in the metric space by semigroups of positive operators, although the existence of a metric was still assumed. This assumption was removed in \cite{JM2,Me3} providing a theory of semigroup type BMO spaces with no further assumptions on the given space. In particular, we could say that Problems 1 and 2 were solved in \cite{JM2}, but it has been unclear since then how to provide natural CZ conditions which imply $L_\infty \to \mathrm{BMO}$ estimates. In this paper we solve P3 by splitting it into:
\begin{itemize}
\item[P3\hskip0.5pta.] Construct a \lq metric\rq${}$ governing the Markov process.

\item[P3b.] Define \lq metric BMO\rq${}$ spaces which still interpolate with the $L_p$ scale.

\item[P3\hskip1ptc.] Provide CZ conditions giving $L_\infty \to \mathrm{BMO}$ boundedness for metric BMO's.
\end{itemize}

\noindent \textbf{P3a. Markov metrics.} Given a Markov semigroup $\S=(S_t)_{t\geq 0}$ on the semifinite von Neumann algebra $(\M,\tau)$ |in other words, formed by normal self-adjoint cpu maps $S_t$| we introduce a \emph{Markov metric} for it as any family
$$
\Q = \Big\{ \big( R_{j,t}, \sigma_{j,t}, \gamma_{j,t} \big) : \ (j,t) \in \Z_+ \times \R_+ \Big\}
$$ 
composed of completely positive unital (cpu) maps $R_{j,t}: \M \to \M$ and elements $\sigma_{j,t}, \gamma_{j,t}$ of $\M$ with $\gamma_{j,t} \ge \1_\M$, such that the following estimates (which show how the Markov metric governs the Markov semigroup in a controlled way) hold:
\begin{itemize}
\item[i)] \emph{Hilbert module majorization}: $\ \displaystyle \big\langle \xi, \xi \big\rangle_{S_t} \, \le \, \sum_{j \ge 1} \sigma_{j,t}^* \big\langle \xi, \xi \big\rangle_{R_{j,t}} \sigma_{j,t}$.

\item[ii)] \emph{Metric integrability condition}: $\ \displaystyle \mathrm{k_\Q} \, = \, \sup_{t > 0} \Big\| \sum_{j \ge 1}  \sigma_{j,t}^* \gamma_{j,t}^2 \sigma_{j,t} \Big\|_\M^{\frac12} \, < \, \infty$.
\end{itemize}
Here $\langle\, ,\,\rangle _{\Phi}$ is the $\M$-valued inner product on $\M \bar\ten \M$ for any cpu map $\Phi$, given by $\langle a \otimes b, a' \otimes b' \rangle_{\Phi} =  b^* \Phi(a^*a') b'$. Markov metrics are a priori unrelated to Rieffel's quantum metric spaces \cite{Ri1,Ri2}. They present on the contrary vague similarities with abstract formulations of classical CZ theory in the absence of CZ kernels and/or doubling measures \cite{BK,To}. We shall explain what motivates our definition below and we shall also illustrate how Euclidean and other classical metrics fit in. 

\noindent \textbf{P3b. Metric type BMO spaces.} Let $$\|f\|_{\mathrm{BMO}_\S^c} = \sup_{t \ge 0} \Big\| \Big( S_t (f^*f) - (S_t f)^*(S_tf) \Big)^\frac12 \Big\|_\M$$ and $\|f\|_{\mathrm{BMO}_\S} = \max \{ \|f\|_{\mathrm{BMO}_\S^c}, \|f^*\|_{\mathrm{BMO}_\S^c} \}$. We shall define the semigroup type BMO space $\mathrm{BMO}_\S(\M)$ as the weak-$*$ closure of $\M$ in certain direct sum of Hilbert modules determined by $\S = (S_t)_{t \ge 0}$. These spaces interpolate with the $L_p$ scale \cite{JM2}.
Given a Markov metric $\Q$ associated to this semigroup, let us define in addition
$$\|f\|_{\mathrm{BMO}_\Q} = \max \Big\{ \|f\|_{\mathrm{BMO}_\Q^c}, \|f^*\|_{\mathrm{BMO}_\Q^c} \Big\},$$
$$\|f\|_{\mathrm{BMO}_\Q^c} = \sup_{t > 0} \inf_{\begin{subarray}{c} \stackrel{\null}{M_t} \, \mathrm{cpu} \\ M_t: \M \to \M \end{subarray}} \sup_{j \ge 1}  \Big\| \Big( \gamma_{j,t}^{-1} \big[ R_{j,t} |f|^2 - |R_{j,t} f|^2 + |R_{j,t} f - M_t f|^2 \big] \gamma_{j,t}^{-1} \Big)^\frac12 \Big\|_\M.$$ 

\begin{TheoA} 
Let $(\M,\tau)$ be a semifinite von Neumann algebra equipped with a Markov semigroup $\T = (S_t)_{t \ge 0}$. Let us consider a Markov metric $\Q$ associated to $\T = (S_t)_{t \ge 0}$. Then, we find $$\|f\|_{\mathrm{BMO}_\T} \, \lesssim \, \mathrm{k}_\Q \, \|f\|_{\mathrm{BMO}_\Q}.$$ Thus, defining $\mathrm{BMO}_\Q(\M)$ as a subspace of $\mathrm{BMO}_\S(\M)$, it interpolates with $L_p(\M)$. 
\end{TheoA} 

Theorem A1 solves P3b. Its proof is not hard after having defined the right notion of Markov metric and the right BMO norm. Let us note in passing that the term $R_{j,t}f - M_t f$ is there to accommodate nondoubling spaces to our definition in the spirit of Tolsa's RBMO space \cite{To}. As a consequence of Theorem A1, proving $L_\infty \to \mathrm{BMO}$ boundedness for metric BMO's (Problem 3c) implies the same result for semigroup BMO spaces (Problem 3). Of course, one could try to prove such a statement directly, but it seems that the metric/measure relation found with these new notions is crucial for a noncommutative CZ theory. 

\noindent \textbf{P3c. Calder\'on-Zygmund operators.} The commutative idea behind the notion of Markov metric (explained in more detail in the body of the paper) is to find pointwise majorants of the integral kernels of our semigroup $\S = (S_t)_{t \ge 0}$, so that we can dominate $S_t$ by certain sum of averaging operators over a distinguished family of measurable sets $\Sigma_{j,t}(x)$. These sets may be considered as the \lq balls\rq${}$ in the Markov metric. In the noncommutative setting, this pointwise estimates must be written in terms of the given Hilbert module majorization and the cpu maps $R_{j,t}$ must be averages over certain projections $q_{j,t}$. Making this precise in full generality is one of the challenges of our algebraic approach and too technical to be explained at this point of the paper. A simple model case is given by 
\begin{equation} \tag{Avg} \label{Eq-Avg}
R_{j,t}f \, = \, (id \otimes \tau)(q_{j,t})^{-\frac12} (id \otimes \tau) \big( q_{j,t} (\1 \otimes f) q_{j,t} \big) (id \otimes \tau)(q_{j,t})^{- \frac12}
\end{equation}
for certain family of projections $q_{j,t} \in \M \bar\otimes \M$. The linear map $\widehat{R}_{j,t}(\1 \otimes f) = R_{j,t} f$ trivially amplifies to $\M \bar\otimes \M$. We may also consider similar formulas for the cpu maps $M_t$ in the metric BMO norm. \eqref{Eq-Avg} allows to identify the Markov metric in terms of the \lq balls\rq${}$ $q_{j,t}$ instead of the corresponding averaging maps $R_{j,t}$.  

As it happens in classical Calder\'on-Zygmund theory, we need to impose some additional properties in the Markov metric to establish a good relation with the underlying (noncommutative) measure. We have split these into \emph{algebraic} and \emph{analytic} conditions, further details will be given in the text. Let us just mention that the algebraic ones are inherent to noncommutativity and hold trivially in commutative cases. The analytic ones provide forms of Jensen's inequality and a measure/metric growth condition. Once we know the Markov metric satisfies these conditions, we may introduce Calder\'on-Zygmund operators. Assume that $T(\A_\M) \subset \M$ for a map $T$ acting on a weak-$*$ dense subalgebra $\A_\M \subset \M$. The goal is to establish sufficient \emph{Calder\'on-Zygmund conditions} on $T$ for $L_\infty \to \mathrm{BMO}_c$ boundedness. These are noncommutative forms of standard properties. Again, it is unnecessary to introduce them here in full generality, we do it in Section \ref{Sect-CZ}. In the model case above, our CZ conditions are:
\begin{itemize}
\item[i)] \emph{$L_\infty(L_2^c)$-boundedness} 
$$\Big\| (id \otimes \tau) \Big( (id \otimes T)(x)^*(id \otimes T)(x) \Big)^\frac12 \Big\|_\M \, \lesssim \, \big\| (id \otimes \tau) (x^*x)^\frac12 \big\|_\M.$$

\vskip5pt

\item[ii)] \emph{Size \lq kernel' conditions} 

\vskip3pt

\begin{itemize}
\item[$\bullet$] $ \displaystyle \widehat{M}_{t} \Big( \big| (id \otimes T) \big( (\1 \otimes f) (A_{j,t} - a_{t}) \big) \big|^2 \Big) \, \lesssim \, \gamma_{j,t}^2 \|f\|_\M^2$, 

\vskip3pt

\item[$\bullet$] $\displaystyle \widehat{R}_{j,t} \Big( \big| (id \otimes T) \big( (\1 \otimes f) (A_{j,t} - a_{j,t}) \big) \big|^2 \Big) \, \lesssim \,  \gamma_{j,t}^2 \|f\|_\M^2$,
\end{itemize}

\vskip3pt

\noindent for certain family of operators $A_{j,t}, a_{j,t} \in \M \bar\otimes \M$ with $A_{j,t} \ge a_{j,t}$.

\vskip5pt

\item[iii)] \emph{H\"ormander \lq kernel' conditions}

\vskip3pt

\begin{itemize}
\item[$\bullet$] $\displaystyle \Phi_{j,t} \Big( \big| \delta \big( (id \otimes T) \big( (\1 \otimes f)(\mathbf{1} \hskip1.5pt - \hskip1.5pt a_{j,t}) \big) \big) \big|^2 \Big) \, \lesssim \, \gamma_{j,t}^2 \|f\|_\M^2$,

\vskip3pt

\item[$\bullet$] $\displaystyle \Psi_{j,t} \Big( \big| \delta \big( (id \otimes T) \big( (\1 \otimes f)(\mathbf{1} - A_{j,t}) \big) \big) \big|^2 \Big) \, \lesssim \, \gamma_{j,t}^2 \|f\|_\M^2,$
\end{itemize}

\vskip3pt

\noindent for certain family of cpu linear maps $\Phi_{j,t}, \Psi_{j,t}: \M \bar\otimes \M \to \M$. 
\end{itemize}

In condition ii), $A_{j,t}$ and $a_{j,t}$ play the role of \lq dilated balls\rq${}$ from $q_{j,t}$. In the last condition, $\delta$ is the derivation $x \mapsto x \otimes \1 - \1 \otimes x$ acting on the second leg of the tensor product. In the Euclidean case, these conditions reduce to $L_2$-boundedness and the classical size/smoothness conditions for the kernel. Our general conditions include many more amplification algebras and derivations, other than $\M \bar\otimes \M$ and $\delta$. Any map $T: \A_\M \to \M$ satisfying the above CZ-conditions will be called a \emph{column} CZ-\emph{operator}. 

\begin{TheoC} \label{ThmA2}
Let $(\M,\tau)$ be a semifinite von Neumann algebra equipped with a Markov semigroup $\S = (S_t)_{t \ge 0}$ with associated Markov metric $\Q $ fulfilling our algebraic and analytic assumptions. Then, any column \emph{CZ}-operator $T$ defines a bounded operator $$T \hskip-1pt : \A_{\M} \to \mathrm{BMO}_\Q^c(\M).$$ Interpolation and duality give similar $($symmetrized$)$ conditions for $L_p$-boundedness.  
\end{TheoC}

A generalized form of Theorem A2 is the main result of this paper. It is easy to recover Euclidean CZ-extrapolation from it. In the Euclidean and many other doubling scenarios, the size kernel condition ii) does not play any role. Our next goal is to explore how the general form of Theorem A2 applies in concrete von Neumann algebras with specific Markov metrics. 

\vskip5pt

\null

\begin{center}
\emph{Applications}
\end{center}

Algebraic Calder\'on-Zygmund theory applies in classical and noncommutative measure spaces. In the commutative context, we shall limit ourselves to prove that algebraic and classical theories match in three important cases: Euclidean spaces with both Lebesgue or Gaussian measures and Riemannian manifolds with non-negative Ricci curvature. We shall not explore further implications in new commutative scenarios, like abstract probability spaces or fractals equipped with specific Dirichlet forms. In the noncommutative context, we start by analyzing the model case of matrix-valued functions from a very general viewpoint. We also consider Calder\'on-Zygmund operators over matrix algebras, generalizing triangular truncations as the archetype of singular integral operator. Most importantly, our abstract theory applies to quantum Euclidean spaces and quantum groups, which constitute our main motivations in this paper. 

It will be useful to specify the form that our Calder\'on-Zygmund operators take when come associated to a concrete kernel. Our applications below include CZ conditions on the kernel. In the basic model case above, we set  
\begin{equation} \tag{Ker 1} \label{Eq-Ker1}
T_k f \, = \, (id \otimes \tau) \big( k (\1 \otimes f) \big)
\end{equation}
for some kernel $k$ affiliated to $\M \bar\otimes \M_{\mathrm{op}}$. Recall that the opposite structure ($\M_{\mathrm{op}}$ is the same algebra $\M$ endowed with the reversed product $a \cdot b = ba$) in the second tensor leg of the kernel for this (standard) model was already justified in \cite{GPJP}. It is a feature of CZ theory which can only be witnessed in noncommutative algebras. It will also be useful to generalize a bit our model case before analyzing any concrete application. Consider an auxiliary von Neumann algebra $\A$ equipped with a n.s.f. trace $\varphi$, a $*$-homomorphism $\sigma: \M \to \A \bar\otimes \M$ and the representation 
\begin{equation} \tag{Ker 2} \label{Eq-Ker2}
S_{\tilde{k}} f \, = \, (id \otimes \varphi) \big( \tilde{k} \ \mathrm{flip} \circ \sigma(f) \big)
\end{equation}
for some kernel $\tilde{k}$ affiliated to $\M \bar\otimes \A_{\mathrm{op}}$. Of course, when $\A = \M$ and $\sigma(f) = \1 \otimes f$ we recover our model case above, with kernel representation \eqref{Eq-Ker1}. This more general framework requires to redefine $R_{j,t}$ in \eqref{Eq-Avg} and the CZ conditions, as we shall do in the body of the paper. The advantage is to take $\A$ as an elementary (commutative) algebra, from which we can transfer metric information. One may think of $\sigma$ as a corepresentation in the terminology of quantum groups. Theorem A2 still holds in this case. We shall refer to \emph{intrinsic} or \emph{transferred} theories when using the model case $\A = \M$ or its generalization respectively.    

\noindent \textbf{Quantum Euclidean spaces.} As geometrical spaces with noncommuting spatial coordinates, quantum Euclidean spaces have appeared frequently in the literature of mathematical physics, in the contexts of string theory  and noncommutative field theory. These algebras  play the role of a central and testing example in noncommutative geometry as well. The singular integral operators on quantum Euclidean spaces naturally appear in the recent study of  Connes'  quantized calculus  \cite{LSZ17, MSX19, SZ18} and noncommutative harmonic analysis \cite{CXY, GPJP, HLP18, XXY18}. Let $$\Theta \in M_n(\R)$$ be anti-symmetric. Briefly, the quantum Euclidean space $\RR_\Theta$ is the von Neumann algebra generated by certain family of unitaries $\{ u_j(s) : 1\leq j\leq n, s\in \R \}$ satisfying 
$$\begin{array}{c} u_j(s)u_j(t)=u_j(s+t), \\ [3pt] u_j(s)u_k(t)=e^{2\pi i \Theta_{jk}st}u_k(t)u_j(s). \end{array}$$
Define  $\lambda_\Theta (\xi)=u_1(\xi_1)u_2(\xi_2)\cdots u_n(\xi_n)$ and set   
$$
f = \int_{\R^n} \check{f}_\Theta(\xi) \lambda_\Theta (\xi)\,d\xi = \lambda_\Theta (\check{f}_\Theta).
$$
for $\check{f}_\Theta \in \mathcal{C}_c(\R^n)$. The trace on $\RR_\Theta$ is determined by
$$
\tau_\Theta(f) = \tau_\Theta \left( \int_{\R^n} \check{f}_\Theta(\xi) \lambda_\Theta (\xi)\, d\xi \right) = \check{f}_\Theta (0).
$$
When $\Theta = 0$, $L_p(\RR_\Theta,\tau_\Theta)$ reduces to $L_p(\R^n)$ with the Lebesgue measure. Precise definitions and a theory of singular integrals for $\RR_\Theta$ appears in \cite{GPJP}. The main result relies on gradient kernel conditions for the intrinsic model \eqref{Eq-Ker1}. \hskip-1pt Remarkably, we show in this paper that the transference model \eqref{Eq-Ker2} $$\sigma_\Theta: \RR_\Theta \ni \lambda_\Theta(\xi) \mapsto \exp_\xi \otimes \lambda_\Theta(\xi) \in L_\infty(\R^n) \bar\otimes \RR_\Theta$$ goes further, since it just requires H\"ormander type smoothness for the kernel. Here $\exp_\xi$ stands for the $\xi$-th character $\exp(2 \pi i \langle \xi, \cdot \rangle)$ in $\R^n$. There is a close relation between both models in this case $$T_k(f) = S_{\tilde{k}}(f) \quad \mbox{for} \quad k = \tilde{\pi}_\Theta(\tilde{k}) \quad \mbox{and} \quad \tilde{\pi}_{\Theta}(m \otimes \exp_\xi) = m \lambda_\Theta(\xi)^* \otimes \lambda_\Theta(\xi).$$ Another crucial map is the $*$-homomorphism
$$\pi_\Theta: L_\infty(\R^n) \ni \exp_\xi \mapsto \lambda_\Theta(\xi)\otimes \lambda_\Theta (\xi)^* \in \mathcal{R}_\Theta \bar{\otimes} \mathcal{R}_\Theta^{\rm{op}}.$$ If $\mathrm{B}_R$ denotes the Euclidean $R$-ball centered at the origin, define the projections $a_R = \pi_\Theta(1_{5\mathrm{B}_R})$ and $a_R^\perp = \1 - a_R$. Set $k_\sigma = (\sigma_\Theta \otimes id_{\RR_\Theta^{\mathrm{op}}})(k) \in L_\infty(\R^n) \bar\otimes \RR_\Theta \bar\otimes \RR_\Theta^{\mathrm{op}}$ and define the derivation $\delta \varphi (x,y) = \varphi(x) - \varphi(y)$ to set the kernel condition in $L_\infty(\R^n) \bar\otimes \RR_\Theta \bar\otimes \RR_\Theta^{\mathrm{op}}$
\begin{equation} \tag{H\"or} \label{Eq-Hor}
\sup_{|x| \le R, |y| \le R} \Big| \delta \Big((id \ten id \ten \tau_\Theta) \big[ k_\sigma(\1 \otimes \1 \otimes f) (\1 \otimes a_R^\perp) \big] \Big) (x,y) \Big|  
\lesssim \|f\|_{\RR_\Theta}.
\end{equation}
As we shall justify in the paper, \eqref{Eq-Hor} is the right form of H\"ormander kernel condition in this framework. The column BMO-norm admits in $\RR_\Theta$ an equivalent form
$$\|f\|_{\mathrm{BMO}_c(\RR_\Theta)} \approx \|\sigma _\Theta (f)\|_{{\rm BMO}_c(\R^n; \RR_\Theta)}$$
for the operator-valued BMO space ${\rm BMO}_c(\R^n; \RR_\Theta)$ from \cite{Me2}. These are all the ingredients to obtain Calder\'on-Zygmund extrapolation over quantum Euclidean spaces. Namely, the general form of Theorem A2 then yields the following theorem. 

\begin{TheoB}
$T_k$ is bounded from $\RR_\Theta$ to $\mathrm{BMO}_c(\RR_\Theta)$ provided$\hskip1pt :$
\begin{itemize}
\item[\emph{i)}]$T_k$ is bounded on $L_2(\RR_\Theta)$.

\item[\emph{ii)}] The kernel condition \eqref{Eq-Hor} holds.
\end{itemize}
Interpolation and duality give similar $($symmetrized$)$ conditions for $L_p$-boundedness.
\end{TheoB}

Theorem B1 improves the main CZ extrapolation theorem in \cite{GPJP} by reducing the gradient kernel condition there to the (more flexible) H\"ormander integral condition above, as we shall prove along the paper. In fact, the result which we shall finally prove is slightly more general than the statement above. 

\vskip3pt

\noindent \textbf{Quantum groups.} Let $\G$ be a locally compact group with a left invariant Haar measure $\mu$. When $\G$ is abelian, the Fourier transform carries the convolution algebra $L_1(\G,\mu)$ into the multiplication algebra $L_\infty (\widehat{\G},\widehat{\mu})$ associated to the dual group with its (normalized) Haar measure. However, when $\G$ is not abelian, we can not construct the dual group and the multiplication algebra above becomes the group von Neumann algebra which is generated by the left regular representation of $\G$. These algebras are basic models of (noncommutative, but still cocommutative) quantum groups, over which we shall study singular integrals. 

Let $\mathbb{G}$ be a locally compact quantum group  |precise definitions in the body of the paper| with comultiplication $\Delta$ and left-invariant and right-invariant Haar weights $\psi$, $\varphi$. Given a weak-$*$ dense subspace $\A$ of $L_\infty(\mathbb{G})$ and a linear map $T$ satisfying $T(\A) \subset  L_\infty(\mathbb{G})$, it is is a convolution map when
$$(T \otimes id_{\mathbb{G}}) \circ \Delta = \Delta \circ T = (id_{\mathbb{G}} \otimes T) \circ \Delta.$$
To simplify the problem, we shall consider the case where $\mathbb{G}$ admits an $\alpha$-doubling intrinsic Markov metric. That is, the projections which generate the cpu maps $R_{j,t}$'s satisfy $$\frac{\psi (q_{\alpha(j),t})}{\psi (q_{j,t})} \leq c_\alpha$$ for a strictly increasing function $\alpha: \N \rightarrow \N$ with $\alpha(j) > j$ and a constant $ c_\alpha$.  
 
%define the auxiliary map $\Phi_{\mathrm{flip}} = \mathrm{flip} \circ (id_\Mn \otimes \Phi) \circ \mathrm{flip}$ on $\A_\Mn \otimes \Mn$. Then  

\begin{TheoD} 
Let $\mathbb{G}$ be a locally compact quantum group and assume it comes equipped with a convolution semigroup $\T = (S_t)_{t \ge 0}$ which admits an
$\alpha$-doubling intrinsic Markov metric. Let $T: \A \to L_\infty(\mathbb{G})$ be a convolution map defined on a weakly dense $*$-subalgebra $\A$ of $L_\infty(\mathbb{G})$ such that
\begin{itemize}
\item[\emph{i)}] $T$ is bounded on $L_2(\mathbb{G})$.

\vskip3pt

\item[\emph{ ii)}]  $\displaystyle \frac{1}{|\psi(q_{j,t})|^2} (\psi\ten \psi) \Big( (q_{j,t} \otimes q_{j,t}) \big| \delta \big(
T( f q_{\alpha(j),t}^{\perp}\big)\big|^2 \Big) 
\lesssim  \|f\|_{L_\infty(\mathbb{G})}^2$.
\end{itemize}
Then, the linear map $T$ extends to a  bounded map $T: L_\infty(\mathbb{G}) \to \mathrm{BMO}_{\T}^c(L_\infty(\mathbb{G})).$
\end{TheoD}

As usual, $L_p$ estimates follow from symmetrized conditions by interpolation and duality. In fact, we shall prove a more general statement which incorporates tensor products with an additional algebra $(\M,\tau)$. Theorem B2 is proved one more time from Theorem A2. In fact, it is conceivable to remove the $\alpha$-doubling restriction and still make the convolution map bounded under an additional size kernel condition as Theorem A2 indicates. 

\vskip3pt

\noindent \textbf{Noncommutative transference.} In a different direction, we shall finish this paper with a section devoted to noncommutative forms of Calder\'on-Cotlar method of transference \cite{C, CW, Co}. The basic idea is to transfer $L_p$ estimates of convolution maps on quantum groups to a much wider class of maps which arise by transference. We refer to \cite{CPPR,CS,CXY,NR,PRdlS,R2} for other forms of transference in the context of group von Neumann algebras and quantum tori. 

\section{\bf Markov metrics}
\label{Sect1}

An abstract form of Calder\'on-Zygmund theory incorporating noncommutative algebras lacks standard geometrical tools. Given a Markov semigroup on a von Neumann algebra |a semigroup of normal cpu self-adjoint maps on the given algebra| we shall construct a \lq metric\rq${}$ governing the Markov process. Our model case in a commutative measure space $(\Omega, \mu)$ is a Markov semigroup of linear maps of the form $$S_t f(x) = \int_{\Omega} s_t(x,y) f(y) \, d\mu(y).$$ The idea is to find pointwise majorants of the form 
\begin{equation} \label{Eq-UpperHM}
s_t(x,y) \, \le \, \sum_{j=1}^\infty \frac{|\sigma_{j,t}(x)|^2}{\mu(\Sigma_{j,t}(x))} \chi_{\Sigma_{j,t}(x)}(y),
\end{equation}
so that $S_t f(x)$ is dominated by a given combination of averaging operators over certain measurable sets $\Sigma_{j,t}(x)$. These sets will determine some sort of metric on $(\Omega,\mu)$ under additional integrability properties. Naively, we may think of them as balls or coronas around $x$ in the hidden metric with radii depending on $(j,t)$. In this section we  formalize this idea and construct BMO spaces with respect to the associated \lq Markov metric\rq${}$ which satisfy the expected interpolation results. 

\subsection{Hilbert modules}

A noncommutative measure space is a pair $(\M,\tau)$ formed by a semifinite von Neumann algebra $\M$ and a n.s.f. trace $\tau$. We assume in what follows that the reader is familiar with basic terminology from noncommutative integration theory \cite{KR,Ta}. Nonexpert readers may proceed by fixing a measure space $(\Omega,\mu)$ with $\M = L_\infty(\Omega)$ and $\tau$ the integral operator associated to $\mu$. Given a cpu map $\Phi: \M \to \M$ we may construct the Hilbert module $\M \bar\otimes_\Phi \M$. Namely consider the seminorm on $\M \otimes \M$ $$\|\xi\|_{\M \bar\otimes_\Phi \M} = \big\| \sqrt{\langle \xi, \xi \rangle_\Phi} \big\|_\M$$ determined by the $\M$-valued inner product $$\Big\langle \summ_j a_j \otimes b_j, \summ_k a'_k \otimes b'_k \Big\rangle_\Phi = \summ_{j,k} b_j^* \Phi(a_j^*a'_k) b'_k.$$ 
Then $\M \bar\otimes_\Phi \M$ will stand for the completion in the topology determined by $\xi_\alpha \to \xi$ when $\tau ( \langle \xi - \xi_\alpha, \xi - \xi_\alpha \rangle_\Phi \, g ) \to 0$ for all $g \in L_1(\M)$. When $\Phi$ is normal, the abstract characterization of Hilbert modules \cite{Ps} yields a weak-$*$ continuous right $\M$-module map $\rho: \M \bar\otimes_\Phi \M \to \H^c\bar\otimes \M$ satisfying $\langle \xi, \eta \rangle_\Phi = \rho(\xi)^*\rho(\eta)$. Let us collect a few properties which will be instrumental along this paper. 

\begin{lemma} \label{ModuleLemma}
Given a cpu map $\Phi: \M \to \M \hskip-1pt:$
\begin{itemize}
\item[\emph{i)}]$\displaystyle \big\langle \xi_1 + \xi_2, \xi_1 + \xi_2 \big\rangle_\Phi \le 2 \big\langle \xi_1, \xi_1 \big\rangle_\Phi + 2 \big\langle \xi_2, \xi_2 \big\rangle_\Phi$,

\vskip3pt

\item[\emph{ii)}] $\big\| f \otimes \1_\M - \1_\M \otimes \Phi f \big\|_{\M \bar\otimes_\Phi \M} = \big\| \Phi |f|^2 - |\Phi f|^2 \big\|_\M^{\frac12}$,

\vskip3pt

\item[\emph{iii)}] $\big| \Phi f - g \big|^2 \le \big\langle f \otimes \1_\M - \1_\M \otimes g, f \otimes \1_\M - \1_\M \otimes g \big\rangle_\Phi$,

\vskip6pt

\item[\emph{iv)}] $\displaystyle \big\| f \otimes \1_\M - \1_\M \otimes \Phi f \big\|_{\M
\bar\otimes_\Phi \M} \sim \inf_{g \in \M} \big\| f \otimes \1_\M
- \1_\M \otimes g \big\|_{\M \bar\otimes_\Phi \M}$,

\vskip3pt

\item[\emph{v)}] If $\Phi \le _{cp} \sum_k \beta_k \Psi_k$, then $$\|\xi\|_{\M \bar\otimes_\Phi \M} \le \Big( \summ_k \beta_k\|\xi\|_{\M \bar\otimes_{\Psi_k} \M}^2 \Big)^{\frac12}.$$
\end{itemize}
\end{lemma}

\dem The first inequality follows from hermitianity of the inner product and the identity $\langle \xi, \eta \rangle_\Phi = \rho(\xi)^* \rho(\eta)$ explained above. The second one is straightforward from the definition of $\M \bar\otimes_\Phi \M$. The third inequality follows from Kadison-Schwarz inequality after expanding both sides. The lower estimate in iv) holds trivially with constant 1, while the upper estimate holds with constant 2 since $$f \otimes \1_\M - \1_\M \otimes \Phi f = (f \otimes \1_\M - \1_\M \otimes g) - (\1_\M \otimes (\Phi f-g))$$ and the second term on the right hand side is estimated using iii). Finally, for the last inequality let $\xi = \sum_k a_k \otimes B_j$ and define the column matrices $A^* = \sum_k a_k^* \otimes e_{k1}$ and $B= \sum_k B_j \otimes e_{k1}$. Then we find \\ \vskip-8pt $\hskip-10pt \begin{array}{rcl} \displaystyle \big\langle \xi, \xi \big\rangle_\Phi \, = \, \summ_{j,k} b_j^* \Phi(a_j^* a_k) B_j & = & \displaystyle  B^* \Phi(A^*A) B \\ & \le & \displaystyle  \summ_k \beta_k B^* \Psi_k(A^*A) B \, = \, \summ_k \beta_k \big\langle \xi, \xi \big\rangle_{\Psi_k}. \hskip10pt \square \end{array}$ 

\vskip3pt

Let $(\M,\tau)$ denote a noncommutative measure space equipped with a Markov semigroup $\T = (S_t)_{t \ge 0}$ acting on it. A \emph{Markov metric} associated to $(\M,\tau)$ and $\T$ is determined by a family 
$$
\Q = \Big\{ \big( R_{j,t}, \sigma_{j,t}, \gamma_{j,t} \big) :  \ (j,t) \in \Z_+ \times \R_+ \Big\}
$$ 
where $R_{j,t}: \M \to \M$ are completely positive unital maps and $\sigma_{j,t}, \gamma_{j,t}$ are elements of the von Neumann algebra $\M$ with $\gamma_{j,t} \ge \1_\M$, so that the estimates below hold: 
\begin{itemize}
\item[i)] \emph{Hilbert module majorization}: $\displaystyle \big\langle \xi, 
\xi \big\rangle_{S_t} \, \le \, \sum_{j \ge 1} \sigma_{j,t}^* \big\langle \xi, \xi \big\rangle_{R_{j,t}} \sigma_{j,t}$,

\item[ii)] \emph{Metric integrability condition}: \ $\displaystyle \mathrm{k_\Q} \, = \, \sup_{t > 0} \Big\| \sum_{j \ge 1}  \sigma_{j,t}^* \gamma_{j,t}^2 \sigma_{j,t} \Big\|_\M^{\frac12} \, < \, \infty$.
\end{itemize}

Our notion of Markov metric is easily understood for our (commutative) model case above. Let $\T = (S_t)_{t \ge 0}$ be a Markov semigroup on $(\Omega,\mu)$ with associated kernels $s_t(x,y)$ satisfying the pointwise estimate \eqref{Eq-UpperHM}. Given $\xi: \Omega \times \Omega \to \C$ essentially bounded, we have 
\begin{equation} \label{Eq-CommKernel}
\big\langle \xi, \xi \big\rangle_{S_t} \hskip-2pt = \hskip-2pt \int_\Omega s_t(x,y) |\xi(x,y)|^2 \, d\mu(y) \le \sum_{j=1}^\infty \frac{|\sigma_{j,t}(x)|^2}{\mu(\Sigma_{j,t}(x))} \int_{\Sigma_{j,t}(x)} |\xi(x,y)|^2 \, d\mu(y).
\end{equation}
This means that $R_{j,t}f(x)$ is the average of $f$ over the set $\Sigma_{j,t}(x)$. Reciprocally, if we take $\xi_k(x,y) = \phi_k(y_0-y)$ to be an approximation of identity around $y_0$, we recover the pointwise estimates for the kernel $s_t(x,y_0)$. In other more general contexts, the upper bounds for the kernel or even the kernel description of the semigroup might not have the same form. As we shall see, many of these cases can still be handled via Hilbert module majorization. We shall provide along the paper a wide variety of examples which fall into these possible classes. 

\subsection{Semigroup $\mathrm{BMOs}$}

Given a noncommutative measure space $(\M,\tau)$ and a Markov semigroup $\T = (S_t)_{t \ge 0}$ acting on $(\M,\tau)$, we may define the \emph{semigroup} $\mathrm{BMO}_\T$-\emph{norm} as $$\|f\|_{\mathrm{BMO}_\T} = \max \Big\{ \|f\|_{\mathrm{BMO}_\T^r}, \|f\|_{\mathrm{BMO}_\T^c} \Big\},$$ where the row and column BMO norms are given by
\begin{eqnarray*}
\|f\|_{\mathrm{BMO}_\T^r} & = & \sup_{t \ge 0} \Big\| \Big( S_t (ff^*) - (S_t f)(S_tf)^*
\Big)^\frac12 \Big\|_\M, \\ \|f\|_{\mathrm{BMO}_\T^c} & = & \sup_{t \ge 0} \Big\| \Big( S_t (f^*f) - (S_t f)^*(S_tf)
\Big)^\frac12 \Big\|_\M.
\end{eqnarray*}
This definition makes sense since we know from the Kadison-Schwarz inequality that $|S_tf|^2 \le S_t|f|^2$. The null space of this seminorm is $\mathrm{ker} A_\infty$, the fixed-point subspace of our semigroup. Indeed, if $\|f\|_{\mathrm{BMO}_\T} = 0$ we know from \cite{Ch} that $f$ belongs to the multiplicative domain of $S_t$, so that $$\tau(gf) = \tau(S_{t/2}(gf)) = \tau(S_{t/2}(g) S_{t/2}(f)) = \tau(g S_t(f)).$$ This proves that $f$ is fixed by the semigroup. Reciprocally, $\mathrm{ker} A_\infty$ is a $*$-subalgebra of $\M$ by \cite{JX2}. Thus, the seminorm vanishes on $\mathrm{ker} A_\infty$. In particular, we obtain a norm after quotienting out $\mathrm{ker} A_\infty$. Letting $w_t(f) = f \otimes \1- \1 \otimes S_tf$, this provides us with a map 
$$f \in \M \stackrel{w}{\longmapsto} \big( w_t(f) \big)_{t \ge 0} \in \bigoplus_{t \ge 0} \M \bar\otimes_{S_t} \M$$ which becomes isometric when we equip $\M$ with the norm in $\mathrm{BMO}_\T^c$. Define $\mathrm{BMO}_\T^c$ as the weak-$*$ closure of $w(\M)$ in the latter space. Similarly, we may define $\mathrm{BMO}_\T$ as the intersection $\mathrm{BMO}_\T^r \cap \mathrm{BMO}_\T^c$, where the row BMO follows by taking adjoints above. The natural operator space structure is given by $$M_m(\mathrm{BMO}_\T(\M)) = \mathrm{BMO}_{\widehat{\T}}(M_m(\M)) \quad \mbox{with} \quad \widehat{S}_t = id_{M_m} \otimes S_t.$$ 

\begin{remark}
\emph{Incidentally, we note that $\mathrm{BMO}_\T$ is written as $bmo(\T)$ in \cite{JM2}.}
\end{remark}

It will be essential for us to provide interpolation results between semigroup type $\mathrm{BMO}$ spaces and the corresponding noncommutative $L_p$ spaces. It is a hard problem to identify the minimal regularity on the semigroup $\T = (S_t)_{t \ge 0}$ which suffices for this purpose. The first substantial progress was announced in a preliminary version of \cite{JM}, where the gradient form $2 \Gamma(f_1,f_2) = A(f_1^*)f_2 + f_1^*A(f_2) - A(f_1^*f_2)$ with $A$ the infinitesimal generator of the semigroup, was a key tool in finding sufficient regularity conditions in terms of \emph{nice enough} Markov dilations. However we know after \cite{JM2} an even sharper condition. Consider the sets
\begin{eqnarray*}
\A_\T f & = & \Big\{ B_tf = \frac{1}{t} \big( S_t(f^2) + f^2 - S_t(f)f - fS_t(f) \big) : \ t > 0 \Big\}, \\ \Gamma_1 \M & = & \Big\{ f \in \M_{\mathrm{s.a.}} : \ \A_\T f \ \mbox{is relatively compact in} \ L_1(\M) \Big\},
\end{eqnarray*}
where $\M_{\mathrm{s.a.}}$ denotes the self-adjoint part of $\M$. The family $\A_\T f $ is called uniformly integrable in $L_1(\M)$ if for all $\varepsilon > 0$ there exists $\delta > 0$ such that $\|(B_tf)q\|_1 < \varepsilon$ for every projection $q$ satisfying $\tau(q) < \delta$. It is well-known that $\A_\T f$ is relatively compact in $L_1(\M)$ if and only if it is bounded and uniformly integrable. Let us also recall that $B_tf \to 2 \Gamma(f,f)$ as $t \to 0$. Define $$L_p^\circ(\M) = \Big\{ f \in L_p(\M) \, : \ \lim_{t \to \infty} S_t f = 0 \Big\}.$$ 
As it was explained in \cite{JM2}, the space $L_p^\circ(\M)$ is complemented in $L_p(\M)$ and $[ L_1^\circ(\M), \mathrm{BMO}_\T]$ form an interpolation couple. A Markov semigroup $\T = (S_t)_{t \ge 0}$ satisfying that $\Gamma_1 \M$ is weak-$*$ dense in $\M_{\mathrm{s.a.}}$ is called \emph{regular}. All the semigroups that we handle in this paper are regular. The following result will be crucial in what follows, we refer the reader to \cite{JM2} for a detailed proof. 

\begin{theorem} \label{Interpolation}
If $\T = (S_t)_{t \ge 0}$ is regular on $(\M,\tau)$ \vskip-7pt $$\big[ \mathrm{BMO}_\T, L_p^\circ(\M) \big]_{p/q} \, \simeq_{cb} \,
L_q^\circ(\M) \quad \mbox{for all} \quad 1 \le p < q < \infty.$$ 
\end{theorem}

Note that interpolation against the full space $L_p(\M)$ is meaningless since $\mathrm{BMO}_\T$ does not distinguish the fixed-point space of the semigroup. Very roughly, we shall typically apply the above result to a CZO which is bounded on $L_2(\M)$ and sends a weak-$*$ dense subalgebra $\A$ of $\M$ to $\mathrm{BMO}_\T$. Recalling the projection map $J_p: L_p(\M) \to L_p^\circ(\M)$ and letting $T$ denote the CZO, we find by interpolation and the weak-$*$ density of $\A$ that $$J_p T: L_p(\M) = \big[ \A, L_2(\M) \big]_{2/p} \to \big[ \mathrm{BMO}_\T, L_2^\circ(\M) \big]_{2/p} = L_p^\circ(\M) \subset L_p(\M).$$ To obtain $L_p$ boundedness of $T$, it suffices to assume that $T$ leaves the fixed-point space invariant and is bounded on it. It should be noticed though, that in many cases the $L_p$ boundedness of the CZO follows automatically. For instance, in $\R^n$ with the Lebesgue measure and the heat semigroup, it turns out that $L_p = L_p^\circ$. On the other hand, the fixed-point space for the Poisson semigroup on the $n$-torus is just composed of constant functions and the corresponding projection can be estimated apart regarded as a conditional expectation. Moreover, the same applies for Fourier multipliers on arbitrary discrete groups. The $L_p$ boundedness for $1 < p < 2$ will follow by taking adjoints under certain symmetry on the hypotheses. 

\subsection{Markov metric $\mathrm{BMOs}$} \label{PapagraphMetricBMO}

Let us now introduce a Markov metric type BMO space for von Neumann algebras and relate it with the semigroup type BMO spaces defined above. Given a Markov semigroup $\T = (S_t)_{t \ge 0}$ acting on $(\M,\tau)$, consider a Markov metric $\Q = \{ ( R_{j,t}, \sigma_{j,t}, \gamma_{j,t} ) : (j,t) \in \Z_+ \times \R_+ \}$ as defined above and define $\|f\|_{\mathrm{BMO}_\Q} = \max \{ \|f\|_{\mathrm{BMO}_\Q^c}, \|f^*\|_{\mathrm{BMO}_\Q^c} \}$, where the column BMO-norm is given by $$\sup_{t > 0} \inf_{\stackrel{\null}{M_t} \, \mathrm{cpu}} \sup_{j \ge 1}  \Big\| \Big( \gamma_{j,t}^{-1} \big[ R_{j,t} |f|^2 - |R_{j,t} f|^2 + |R_{j,t} f - M_t f|^2 \big] \gamma_{j,t}^{-1} \Big)^\frac12 \Big\|_\M$$ 
and the infimum runs over cpu maps $M_t: \M \to \M$. Since $\gamma_{j,t} \ge \mathbf{1}_\M$, the inverses exist and $L_\infty(\M)$ embeds in $\mathrm{BMO}_\Q$. Indeed, using that $R_{j,t}$ and $M_t$ are cpu, the square braket above is bounded by $5 \|f\|_\infty^2 \mathbf{1}_\M$ and $\gamma_{j,t}^{-2} \le \mathbf{1}_\M$. The row norm is estimated in the same way. Now, recalling the value of the constant $\mathrm{k}_\Q$ in our definition of Markov metric, we prove that $\mathrm{BMO}_\Q$ embeds in $\mathrm{BMO}_\T$. 

\begin{theorem} \label{MetricBMOInterp}
Let $(\M,\tau)$ be a noncommutative measure space equipped with a Markov semigroup $\T = (S_t)_{t \ge 0}$. Let us consider a Markov metric $\Q$ associated to $\T = (S_t)_{t \ge 0}$. Then, we find $$\|f\|_{\mathrm{BMO}_\T} \, \lesssim \, \mathrm{k}_\Q \, \|f\|_{\mathrm{BMO}_\Q}.$$ In particular, we see that $L_\infty(\M) \subset \mathrm{BMO}_\Q \subset \mathrm{BMO}_\S$ and $$\big[ \mathrm{BMO}_\Q, L_p^\circ(\M) \big]_{p/q} \, \simeq \, L_q^\circ(\M) \quad \mbox{for all} \quad 1 \le p < q < \infty$$ for any Markov metric $\Q$ associated to a regular semigroup $\T = (S_t)_{t \ge 0}$ on $(\M,\tau)$.
\end{theorem}

\dem Let us set
\begin{eqnarray*}
\xi_t & = & f \otimes \mathbf{1}_\M - \mathbf{1}_\M \otimes M_tf \\ & = & (f \otimes \mathbf{1}_\M - \mathbf{1}_\M \otimes R_{j,t}f) + (\mathbf{1}_\M \otimes (R_{j,t}f - M_tf) ) \ = \ \xi_{j,t}^1 + \xi_{j,t}^2.
\end{eqnarray*}
The assertion follows from Lemma \ref{ModuleLemma} and our definition of Markov metric
\begin{eqnarray*}
\|f\|_{\mathrm{BMO}_\T^c} & \!\! = \!\! & \sup_{t > 0} \Big\| \Big( S_t|f|^2 - |S_tf|^2 \Big)^\frac12 \Big\|_\M \\ [5pt] & \!\! = \!\! & \sup_{t > 0} \big\| f \otimes \mathbf{1}_\M - \mathbf{1}_\M \otimes \hskip1.5pt S_t \hskip0.5pt f \hskip1.5pt \big\|_{\M \bar\otimes_{S_t} \M} \\ [3pt] & \!\! \lesssim \!\! & \sup_{t > 0} \big\| f \otimes \mathbf{1}_\M - \mathbf{1}_\M \otimes M_tf \big\|_{\M \bar\otimes_{S_t} \M} \, = \, \sup_{t > 0} \| \langle \xi_t, \xi_t \rangle_{S_t} \|_\M^{\frac12} \\ 
& \!\! \lesssim \!\! & \sup_{t > 0} \Big\| \sum_{j \ge 1} \sigma_{j,t}^* \big[ \langle \xi_{j,t}^1, \xi_{j,t}^1 \rangle_{R_{j,t}} + \langle \xi_{j,t}^2, \xi_{j,t}^2 \rangle_{R_{j,t}} \big] \sigma_{j,t} \Big\|_\M^{\frac12} \\
& \, \le \,& \mathrm{k}_\Q \|f\|_{\mathrm{BMO}_\Q^c}.
\end{eqnarray*}
The identities are clear. The first inequality follow from Lemma \ref{ModuleLemma} iv), the second one from the Hilbert module majorization associated to the Markov metric and Lemma \ref{ModuleLemma} i). To justify the last inequality, note that the square bracket inside the term on the left equals $R_{j,t}|f|^2 - |R_{j,t}f|^2 + |R_{j,t}f - M_tf|^2$. Hence, left multiplication by $\gamma_{j,t}^{\null}\gamma_{j,t}^{-1}$ and right multiplication by $\gamma_{j,t}^{-1} \gamma_{j,t}^{\null}$ yields the given inequality with $\mathrm{k}_\mathcal{Q}$ the metric integrability constant. The interpolation result follows from Theorem \ref{Interpolation} and the embeddings $L_\infty(\M) \subset \mathrm{BMO}_\Q \subset \mathrm{BMO}_\T$. The proof is complete. \fin

\begin{remark} \label{ossMarkovBMO}
\emph{Let $\xi = \sum_j A_j \otimes B_j$, with $$A_j = \Big( a^j_{\alpha \beta} \Big) \quad \mbox{and} \quad B_j = \Big( b^j_{\alpha \beta} \Big)$$ elements of $M_m(\M)$. If $\widehat{S}_t = id_{M_m} \otimes S_t$, it turns out that $$\big\langle \xi, \xi \big\rangle_{\widehat{S}_t} \, = \, \Bigg( \sum_{\alpha=1}^m \Big\langle \underbrace{\sum_{j,\beta} a_{\alpha \beta}^j \otimes b_{\beta \gamma_1}^j}_{\eta_{\alpha, \gamma_1}}, \underbrace{\sum_{j,\beta} a_{\alpha \beta}^j \otimes b_{\beta\gamma_2}^j}_{\eta_{\alpha, \gamma_2}} \Big\rangle_{S_t} \Bigg)_{\gamma_1, \gamma_2} \in M_m(\M).$$ This can be used to provide an operator space structure on $\mathrm{BMO}_\Q$. Namely, the canonical choice for the matrix norms is $M_m(\mathrm{BMO}_\Q(\M)) = \mathrm{BMO}_{\widehat{\Q}}(M_m(\M))$ where the Markov metric on $M_m(\M)$ $$\widehat{\Q} = \Big\{ \big( id_{M_m} \otimes R_{j,t}, \mathbf{1}_{M_m} \otimes \sigma_{j,t}, \mathbf{1}_{M_m} \otimes \gamma_{j,t} \big) \Big\}$$ is associated to the extended semigroup $(\widehat{S}_t)_{t \ge 0}$. Then, we trivially obtain that $\mathrm{k}_{\widehat{\Q}} = \mathrm{k}_\Q < \infty$. However, according to the identity above for $\langle \xi, \xi \rangle_{\widehat{S}_t}$, the Hilbert module majorization takes the form 
$$\Bigg( \sum_{\alpha=1}^m \big\langle \eta_{\alpha,\gamma_1}, \eta_{\alpha,\gamma_2} \big\rangle_{S_t} \Bigg)_{\gamma_1, \gamma_2} \le \hskip8pt \sum_{j \ge 1} \Bigg( \sigma_{j,t}^* \sum_{\alpha=1}^m \big\langle \eta_{\alpha,\gamma_1}, \eta_{\alpha, \gamma_2} \big\rangle_{R_{j,t}} \sigma_{j,t} \Bigg)_{\gamma_1, \gamma_2}.$$ This gives a matrix-valued generalization of our Hilbert module majorization for $\T=(S_t)_{t \ge 0}$ on $\M$, to be checked when we use this o.s.s. Theorem \ref{MetricBMOInterp} yields a cb-embedding of $\mathrm{BMO}_\Q$ into $\mathrm{BMO}_\T$ under this assumption. According to the characterization \eqref{Eq-CommKernel}, it holds for Markov metrics on commutative spaces $(\Omega,\mu)$.}
\end{remark}

\subsection{The Euclidean metric} \label{SectEuclideanMetric}

Before using Markov metrics in our approach to Calder\'on-Zygmund theory, it is illustrative to recover the Euclidean metric from a suitable Markov semigroup. Let $\T = (H_t)_{t \ge 0}$ denote the classical heat semigroup on $\R^n$, with kernels $$h_t(x,y) = \frac{1}{(4\pi t)^{\frac{n}{2}}} \exp \Big(\frac{-|x-y|^2}{4t} \Big).$$ Take $\Q = \big\{ ( R_{j,t}, \sigma_{j,t}, \gamma_{j,t} ) : (j,t) \in \Z_+ \times \R_+ \big\}$ determined by

\vskip8pt

\begin{itemize}
\item $\sigma_{j,t}^2 \equiv \frac{2e}{\sqrt{\pi}} j^{\frac{n}{2}} e^{-j}$ and $\gamma_{j,t}^2 \equiv j^{\frac{n}{2}} \ge 1$,

\vskip5pt

\item $R_{j,t}f(x) = \displaystyle \frac{1}{|\mathrm{B}_{\sqrt{4jt}}(x)|} \int_{\mathrm{B}_{\sqrt{4jt}}(x)} f(y) \, dy$.
\end{itemize}
Note that $\sigma_{j,t}$ and $\gamma_{j,t}$ are allowed to be essentially bounded functions in $\R^n$, but in this case it suffices to take constant functions. In the definition of $R_{j,t}$, we write $\mathrm{B}_r(x)$ to denote the Euclidean ball in $\R^n$ centered at $x$ with radius $r$. It is clear that $R_{j,t}$ defines a cpu map on $L_\infty(\R^n)$. To show that $\Q$ defines a Markov metric, we need to check that it provides a Hilbert module majorization of the heat semigroup and the metric integrability condition holds. The latter is straightforward, while the Hilbert module majorization reduces to check that 
$$h_t(x,y) \, \le \, \frac{2e}{\sqrt{\pi}} \sum_{j \ge 1} \frac{j^{\frac{n}{2}} e^{-j}}{|\mathrm{B}_{\sqrt{4jt}}(x)|} \chi_{\mathrm{B}_{\sqrt{4jt}}(x)}(y).$$ 
This can be justified by determining the unique corona centered at $x$ with radii $\sqrt{4(j-1)t}$ and $\sqrt{4jt}$ where $y$ lives, details are left to the reader. Note that we could have taken $\gamma_{j,t} \equiv 1$ and still obtain a Markov metric. Our choice will be justified below and also in the next section, where we shall need $\gamma_{j,t} \equiv j^{\frac{n}{2}}$ to compare $\mathrm{BMO}_\Q$ with other $\mathrm{BMO}$ spaces which interpolate. Before that, our only evidences that this is the right Markov metric in the Euclidean case are the fact that the $R_{j,t}$'s are averages over Euclidean balls and the isomorphism $$\mathrm{BMO}_\Q \, = \, \mathrm{BMO}_{\R^n},$$ where the latter space is the usual BMO space in $\R^n$ 
$$\|f\|_{\mathrm{BMO}_{\R^n}} \, = \, \sup_{\mathrm{B} \subset \R^n} \Big( \frac{1}{|\mathrm{B}|} \int_\mathrm{B} \big| f(x) - f_\mathrm{B} \big|^2 \, dx \Big)^\frac12.
$$ 
Here, the supremum is taken over all Euclidean balls $\mathrm{B}$ in $\R^n$ and $f_\mathrm{B}$ stands for the average of $f$ over $\mathrm{B}$. Let us justify this isomorphism. If we pick $M_tf(x) = R_{1,t}f(x)$ it follows from a standard calculation that 
\begin{eqnarray} \label{Eq-BMODoubling}
\lefteqn{\hskip20pt \big| R_{j,t}f (x) - M_tf(x) \big|^2} \\ \nonumber & = & \Big| \frac{1}{|\mathrm{B}_{\sqrt{4t}}(x)|} \int_{\mathrm{B}_{\sqrt{4t}}(x)} \big( f(y) - f_{\mathrm{B}_{\sqrt{4jt}}(x)} \big) \, dy \Big|^2 \\ \nonumber & \le & \frac{j^{\frac{n}{2}}}{|\mathrm{B}_{\sqrt{4jt}}(x)|} \int_{\mathrm{B}_{\sqrt{4jt}}(x)} \big| f(y) - f_{\mathrm{B}_{\sqrt{4jt}}(x)} \big|^2 \, dy \ = \ j^{\frac{n}{2}} \big( R_{j,t}|f|^2 - |R_{j,t}f|^2 \big) (x).
\end{eqnarray}
This automatically yields the following inequality $$\|f\|_{\mathrm{BMO}_\Q}^2 \, \lesssim \, \sup_{j,t} \esssup_{x \in \R^n} \frac{1}{|\mathrm{B}_{\sqrt{4jt}}(x)|} \int_{\mathrm{B}_{\sqrt{4jt}}(x)} \big| f(y) - f_{\mathrm{B}_{\sqrt{4jt}}(x)} \big|^2 \, dy \, \le \, \|f\|_{\mathrm{BMO}_{\R^n}}^2.$$ The converse is even simpler, since taking $j=1$ we obtain 
\begin{eqnarray*}
\|f\|_{\mathrm{BMO}_{\R^n}}^2 & = & \sup_{t>0} \esssup_{x \in \R^n} \frac{1}{|\mathrm{B}_{\sqrt{4t}}(x)|} \int_{\mathrm{B}_{\sqrt{4t}}(x)} \big| f(y) - f_{\mathrm{B}_{\sqrt{4t}}(x)} \big|^2 \, dy \\ & = & \sup_{t > 0} \Big\| \gamma_{1,t}^{-1} \big[ R_{1,t}|f|^2 - |R_{1,t}f|^2 \big] \gamma_{1,t}^{-1} \Big\|_{L_\infty(\R^n)} \ \le \ \|f\|_{\mathrm{BMO}_\Q}^2.
\end{eqnarray*}

\begin{remark} \label{Rem-NDTerm}
\emph{The term $|R_{j,t}f - M_tf|$ did not play a significant role at this point. More generally, the above argument also works for any doubling metric space $\Omega$ equipped with a Borel measure $\mu$: $\mu(B(x,2r)) \leq C \mu(B(x,r))$ for every $x \in \Omega$ and $r>0$, with $B(x,r)=\{y\in \Omega : {\rm{dist}}(x,y)\leq r \}$. As we shall see later, the additional term $|R_{j,t}f - M_tf|$ in the $\mathrm{BMO}_\Q$-norm appears to include Tolsa's RBMO spaces \cite{To} in those measure spaces $(\Omega,\mu)$ for which we can find an appropriate Dirichlet form which provides us with a Markov semigroup acting on $(\Omega,\mu)$.}
\end{remark}

\begin{remark} \label{Remark-EquivBMO} 
\emph{A related semigroup BMO norm is $$\|f\|_{\mathbf{BMO}_\T^c} = \sup_{t \ge 0} \Big\| \Big( H_t \big[ |f \hskip1pt - \hskip0.5pt H_t f|^2 \big] \Big)^\frac12 \Big\|_\infty.$$ All the norms consider so far are equivalent for the heat semigroup $\T = (H_t)_{t \ge 0}$ on $\R^n$, generated by the Laplacian $\Delta = \sum_{j=1}^n \partial_{x_j}^2$. In fact, we may also consider by subordination the Poisson semigroup $\mathcal{P} = (P_t)_{t \ge 0}$ on $\R^n$ generated by the square root $\sqrt{- \Delta}$, or even other subordinations \cite{MeietalAIM} . Then, elementary calculations give the following norm equivalences up to dimensional constants}
$$\|f\|_{\mathrm{BMO}_{\R^n}} \sim \|f\|_{\mathrm{BMO}_\P} \sim \|f\|_{\mathbf{BMO}_\P}
\sim \|f\|_{\mathrm{BMO}_\T} \sim \|f\|_{\mathbf{BMO}_\T}
\sim \|f\|_{\mathrm{BMO}_\Q}.$$
\emph{Moreover, let $\RR = L_\infty(\R^n) \bar\otimes \M$ denote the von Neumann algebra tensor product of $L_\infty(\R^n)$ with a noncommutative measure space $(\M,\tau)$. Define the norm in $\mathrm{BMO}_\RR$ as $\|f\|_{\mathrm{BMO}_\RR} = \max \{ \|f\|_{\mathrm{BMO}_\RR^c}, \|f^*\|_{\mathrm{BMO}_\RR^c} \}$, where $$\|f\|_{\mathrm{BMO}_\RR^c} \ = \ \sup_{\mathrm{B} \ \mathrm{balls}} \Big\| \Big( \frac{1}{|\mathrm{B}|} \int_\mathrm{B} \big| f(x) - f_\mathrm{B} \big|^2 dx \Big)^\frac12 \Big\|_\M.$$ Then, the same norm equivalences hold in the semicommutative case $$\|f\|_{\mathrm{BMO}_\RR} \sim \|f\|_{\mathrm{BMO}_{\P_\otimes}} \sim \|f\|_{\mathbf{BMO}_{\P_{\otimes}}}
\sim \|f\|_{\mathrm{BMO}_{\T_\otimes}} \sim \|f\|_{\mathbf{BMO}_{\T_\otimes}},$$ where $S_{\otimes,t} = S_t \otimes id_\M$ and $P_{\otimes,t} = P_t \otimes id_\M$. Moreover, by Remark \ref{ossMarkovBMO}, all these norms are in turn equivalent to the norm in $\mathrm{BMO}_{\Q_\RR}$, with the Markov metric which arises tensorizing the canonical one with the identity/unit of $\M$.}
\end{remark}

\section{\bf Algebraic CZ theory} \label{Sect-CZ}

In classical Calder\'on-Zygmund theory, $L_p$ boundedness of CZOs follows from $L_2$ boundedness under a smoothness condition on the kernel. Our next goal is to identify which are the analogues of these conditions for semifinite von Neumann algebras equipped with a Markov metric, and to show $L_p$ boundedness of CZOs fulfilling them. Our new conditions are certainly surprising. The boundedness for $p=2$ must be replaced by a certain mixed-norm estimate (which reduces in the classical theory to $L_2$ boundedness), while H\"ormander kernel smoothness will be formulated intrinsically without any reference to the kernel. These abstract assumptions will adopt a more familiar form in the specific cases that we shall consider in the forthcoming sections.  

In order to give a Calder\'on-Zygmund framework for von Neumann algebras we start with some initial considerations, which determine the general form of Markov metrics that we shall work with. Consider a Markov metric $\Q$ associated to a Markov semigroup $\T = (S_t)_{t \ge 0}$ acting on $(\M,\tau)$. Then, we shall assume that the cpu maps $R_{j,t}$ from $\Q$ are of the following form 
\begin{equation} \label{Initial}
\begin{array}{c} \M \stackrel{\rho_j}{\longrightarrow} \NN_\rho \stackrel{\E_\rho}{\longrightarrow} \rho_1(\M) \simeq \M, \\ [7pt] R_{j,t}f \, = \, \E_\rho(q_{j,t})^{-\frac12} \E_\rho \big( q_{j,t} \rho_2(f) q_{j,t} \big) \E_\rho(q_{j,t})^{- \frac12}, %\\ [5pt] \hskip32pt \, = \, \E_\rho \Big( \E_\rho(q_{j,t})^{-\frac12} q_{j,t} \rho_2(f) q_{j,t} \E_\rho(q_{j,t})^{-\frac12} \Big), 
\end{array}
\end{equation}
where $\rho_1, \rho_2: \M \to \NN_\rho$ are $*$-homomorphisms into certain von Neumann algebra $\NN_\rho$, the map $\E_\rho: \NN_\rho \to \rho_1(\M)$ is an operator-valued weight and the $q_{j,t}$'s are projections in $\NN_\rho$. In particular, we shall assume that our formula for $R_{j,t}f$ makes sense so that $q_{j,t}$ and $q_{j,t} \rho_2(f) q_{j,t}$ belong to the domain of $\mathsf{E}_\rho$, see Section \ref{SectO-VWeights} for further details. Our model provides a quite general form of Markov metric which includes the Markov metric for the heat semigroup considered before. Indeed, take $\NN_\rho = L_\infty(\R^n \times \R^n)$ with $\rho_1f(x,y) = f(x)$ and $\rho_2f(x,y) = f(y)$. Let $\E_\rho$ be the integral in $\R^n$ with respect to the variable $y$ and set $$q_{j,t}(x,y) \, = \, \chi_{\mathrm{B}_{\sqrt{4jt}}(x)}(y) \, = \, \chi_{\mathrm{B}_{\sqrt{4jt}}(y)}(x) \, = \, \chi_{|x-y| < \sqrt{4jt}}.$$ Then, it is straightforward to check that we recover from \eqref{Initial} the $R_{j,t}$'s for the heat semigroup. Note that the $q_{j,t}(x,\cdot)$'s reproduce in this case all the Euclidean balls in $\R^n$. Morally, this is why we call $\Q$ a Markov metric, since it codifies some sort of underlying metric in $(\M,\tau)$. According to our definition of $\mathrm{BMO}_\Q$, we shall also consider projections $q_t$ in $\NN_\rho$ and cpu maps 
\begin{equation} \label{Mt}
M_tf \, = \, \E_\rho(q_t)^{-\frac12} \E_\rho \big( q_t \rho_2(f) q_t \big) \E_\rho(q_t)^{- \frac12}.
\end{equation}

\subsection{Operator-valued weights}
\label{SectO-VWeights}

In this subsection we briefly review the definition and basic properties of operator-valued weights from \cite{H1,H2}. A unital, weakly closed $*$-subalgebra is called a von Neumann subalgebra. A conditional expectation $\mathcal{E}_\M: \NN \to \M$  onto a von Neumann subalgebra $\M$ is a positive unital projection satisfying the bimodular property $\mathcal{E}_\M(a_1 f \hskip1pt a_2) = a_1 \mathcal{E}_\M(f) \hskip1pt a_2$ for all $a_1, a_2 \in \M$. It is called normal if $\sup_\alpha \mathcal{E}_\M(f_\alpha) = \mathcal{E}_\M(\sup_\alpha f_\alpha)$ for bounded increasing nets $(f_\alpha)$ in $\NN_+$. A normal conditional expectation exists if and only if the restriction of $\tau$ to the von Neumann subalgebra $\M$ remains semifinite \cite{Ta}. Any such conditional expectation is trace preserving $\tau \circ \mathcal{E}_\M = \tau$.

The \emph{extended positive part} $\widehat{\M}_+$ of the von Neumann algebra $\M$ is the set of lower semicontinuous maps $m: \M_{*,+} \to [0,\infty]$ which are linear on the positive cone, $m(\lambda_1 \phi_1 + \lambda_2 \phi_2) = \lambda_1 m(\phi_1) + \lambda_2 m(\phi_2)$ for $\lambda_j \ge 0$ and $\phi_j \in \M_{*,+}$. The extended positive part is closed under addition, increasing limits and is fixed by the map $x \mapsto a^* x a$ for any $a \in \M$. It is clear that $\M_+$ sits in the extended positive part. When $\M$ is abelian, we find $\M \simeq L_\infty(\Omega,\mu)$ for some measure space $(\Omega,\mu)$ and the extended positive part corresponds in this case to the set of $\mu$-measurable functions on $\Omega$ (module sets of zero measure) with values in $[0,\infty]$. A harder characterization of the extended positive part for arbitrary von Neumann algebras was found by Haagerup in \cite{H1}. Assume that $\M$ acts on $\H$ and consider a positive operator $A$ affiliated with $\M$ with spectral resolution $A = \int_{\R_+} \lambda d e_\lambda$. Then, we may construct an associated element in $\widehat{\M}_+$ $$m_A(\phi) \, = \, \int_{\R_+} \lambda d(\phi(e_\lambda)).$$ In general, any $m \in \widehat{\M}_+$ has a unique spectral resolution $$m(\phi) \, = \, \int_{\R_+} \lambda d (\phi(e_\lambda)) + \infty \phi(p)$$ where the $e_\lambda$'s form an increasing family of projections in $\M$ and $p$ is the projection $\mathbf{1}_\M - \lim_{\lambda} e_\lambda$. Moreover, the map $\lambda \mapsto e_\lambda$ is strongly continuous from the right and we find that $e_0 = 0$ iff $m$ does not vanish on $\M_*^+ \setminus \{0\}$, while $p=0$ iff the family of $\phi \in \M_*^+$ with $m(\phi) < \infty$ is dense in $\M_*^+$. 

Operator-valued weights appear as \lq\lq unbounded conditional expectations" and the simplest nontrivial model is perhaps a partial trace $\E_\M = \mathrm{tr}_\A \otimes id_\M$ with $\NN = \A \bar\otimes \M$ and $\A$ a semifinite non-finite von Neumann algebra. In general, an \emph{operator-valued weight} from $\NN$ to $\M$ is just a linear map $$\E_\M: \NN_+ \to \widehat{\M}_+ \quad \mbox{satisfying} \quad \E_\M(a^*fa) \, = \, a^* \E_\M(f) a$$ for all $a \in \M$. As usual, $\E_\M$ is called normal when $\sup_\alpha \E_\M(f_\alpha) = \E_\M(\sup_\alpha f_\alpha)$ for bounded increasing nets $(f_\alpha)$ in $\NN_+$. Since $a^*fb = \frac14 \sum_{k=0}^3 i^{-k} (a+i^kb)^*f (a+i^kb)$ by polarization, we see that bimodularity of conditional expectations is equivalent to $\mathcal{E}_\M(a^*fa) = a^*\mathcal{E}_\M(f)a$ for $a \in \M$. In particular, the fundamental properties which operator-valued weights loose with respect to conditional expectations are unitality and the fact that unboundedness is allowed for the image. Additionally, when $\M = \C$ the map $\E_\M$ becomes an ordinary weight on $\NN$. In analogy with ordinary weights, we take 
 $$L_\infty^c(\NN; \E_\M) \, = \, \Big\{ f \in \NN : \, \big\| \E_\M(f^*f) \big\|_\M < \infty \Big\}.$$
Note that when $\E_\M = \mathrm{tr}_\A \otimes id_\M$ with $\NN = \A \bar\otimes \M$, $L_\infty^c(\NN; \E_\M)$ are the Hilbert space valued noncommutative $L_\infty$ spaces defined in \cite{JLX}, which we denote by $L_2^c(\A)\bar\ten \M$.
Let $\NN_{\E_\M}$  be the linear span of $f_1^*f_2$ with $f_1, f_2 \in L_\infty^c(\NN; \E_\M)$. Then we find
\begin{itemize}
\item[i)] $\NN_{\E_\M} = \mathrm{span} \{f \in \NN_+ : \, \|\E_\M f\| < \infty\}$,

\vskip3pt

\item[ii)] $L_\infty^c(\NN; \E_\M)$ and $\NN_{\E_\M}$ are two-sided modules over $\M$,

\vskip3pt

\item[iii)] $\E_\M$ has a unique linear extension $\E_\M: \NN_{\E_\M} \to \M$ satisfying $$\E_\M(a_1 f a_2) \, = \, a_1 \E_\M(f) a_2 \quad \mbox{with} \quad f \in \NN_{\E_\M} \ \mbox{and} \ a_1, a_2 \in \M.$$
\end{itemize}
In particular, if $\E_\M(\mathbf{1}) = \mathbf{1}$ we recover a conditional expectation onto $\M$. An operator-valued weight $\E_\M$ is called faithful if $\E_\M(f^*f) = 0$ implies $f=0$ and semifinite when $L_\infty^c(\NN; \E_\M)$ is $\sigma$-weakly dense in $\NN$. It is of interest to determine for which pairs $(\NN,\M)$ we may construct n.s.f. operator-valued weights. Among other results, Haagerup proved in \cite{H2} that this is the case when both von Neumann algebras are semifinite and there exists a unique trace preserving one. Note that conditional expectations do not always exist in this case. He also proved that given $\E_{\M_j}$ n.s.f. operator-valued weights in $(\NN_j, \M_j)$ for $j=1,2$, there exists a unique n.s.f. operator-valued weight $\E_{\M_1 \otimes \M_2}$ associated to $(\NN_1 \bar\otimes \NN_2, \M_1 \bar\otimes \M_2)$ such that $(\phi_1 \otimes \phi_2) \circ \E_{\M_1 \otimes \M_2} \, = \, (\phi_1 \circ \E_{\M_1}) \otimes (\phi_2 \circ \E_{\M_2})$ for any pair $(\phi_1, \phi_2)$ of normal semifinite faithful weights on $(\M_1, \M_2)$.

\subsection{Algebraic/analytic conditions} \label{Sect-AAConditions}

The identity $$Tf(x) \, = \, \int_\Omega k(x,y) f(y) \, d\mu(y)$$ is just a vague expression to consider classical Calder\'on-Zygmund operators. It is well-known that particular realizations as above are only meaningful outside the support of $f$ and understanding $k$ as a distribution which coincides with a locally integrable function on $\R^n \times \R^n \setminus \Delta$. Instead of that, we shall not specify any kernel representation of our CZOs since our conditions below will be formulated in a more intrinsic way. These kernel representations will appear later on in this paper with the concrete examples that we shall consider. 

Let $T$ be a densely defined operator on $\M$, which means that $Tf \in \M$ for all $f$ in a weak-$*$ dense subalgebra $\A_\M$ of $\M$. Our assumption does not necessarily hold for classical Calder\'on-Zygmund operators defined in abelian von Neumann algebras $(\M,\tau) = L_\infty (\Omega,\mu)$, but it is true for the truncated singular integral operators satisfying the standard size condition for the kernel, take for instance $\A_\M = \M \cap L_1(\M)$. In particular, this is not a crucial restriction since we shall be able to take $L_p$-limits as far as our estimates below are independent of $T$. Our aim is to settle conditions on $T$ of CZ type assuring that $T: L_\infty(\M) \to \mathrm{BMO}_\Q^c$, provided $(\M,\tau)$ comes equipped with a Markov metric $\Q$. In this paragraph, we establish some preliminary algebraic and analytic conditions on the Markov metric and the CZO. Consider $*$-homomorphisms $\pi_1, \pi_2: \M \to \NN_\pi$ and an operator-valued weight $\E_\pi: \NN_\pi \to \pi_1(\M)$ which may or may not coincide with $\rho_1$, $\rho_2$ and $\E_\rho$ from \eqref{Initial}. Assume there exists a (densely defined) map
\begin{equation} \label{Amplification}
\begin{array}{c} \widehat{T}: \A_{\NN_\pi} \subset \NN_\pi \to \NN_\rho \\ \mbox{satisfying} \quad \widehat{T} \circ \pi_2 = \rho_2 \circ T \quad \mbox{on} \quad \A_\M. \end{array}
\end{equation}

\vskip3pt

\noindent \textbf{Algebraic conditions:}

\begin{itemize}
\item[i)] \emph{$\Q$-monotonicity of $\E_\rho$} $$\E_\rho(q_{j,t} |\xi|^2 q_{j,t}) \le \E_\rho(|\xi|^2)$$ for all $\xi \in \NN_\rho$ and every projection $q_{j,t}$ determined by $\Q$ via  the identity in \eqref{Initial}. Similarly, we assume the same inequality holds when we replace the $q_{j,t}$'s by the $q_t$'s appearing in \eqref{Mt}.

\vskip5pt

\item[ii)] \emph{Right $\mathcal{B}$-modularity of $\widehat{T}$} $$\widehat{T} \big( \eta \, \pi_1 \rho_1^{-1}(b) \big) \, = \, \widehat{T}(\eta) b$$ for all $\eta \in \A_{\NN_\pi}$ and all $b$ lying in some von Neumann subalgebra $\mathcal{B}$ of $\rho_1(\M)$ which includes $\E_\rho(q_t)$, $\E_\rho(q_{j,t})$ and $\rho_1(\gamma_{j,t})$ for every projection $q_t$ and $q_{j,t}$ determined by $\Q$ via the identities in \eqref{Initial} and \eqref{Mt}. 
\end{itemize}

\vskip3pt

As we shall see both conditions trivially hold in the classical theory, where the first condition essentially says that integrating a positive function over a \lq\lq Markov metric ball" is always smaller than integrating it over the whole space, while the second condition allows to take out $x$-dependent functions from the $y$-dependent integral defining $T$. Our conditions remain true in many other situations, which will be explored below in this paper. Nevertheless, condition i) suggests that certain amount of commutativity might be necessary to work with Markov metrics.

To state our analytic conditions we introduce an additional von Neumann algebra $\NN_\sigma$ containing $\NN_\rho$ as a von Neumann subalgebra. Then, we consider derivations $\delta: \NN_\rho \to \NN_\sigma$ given by the difference $\delta = \sigma_1 - \sigma_2$ of two $*$-homomorphisms, so that $\delta(ab) = \sigma_1(a) \sigma_1(b) - \sigma_2(a) \sigma_2(b) = \delta(a) \sigma_1(b) + \sigma_2(a) \delta(b)$ as expected. We also consider the natural amplification maps  $$\widehat{R}_{j,t}:  \NN_\rho \ni \xi \mapsto \E_\rho(q_{j,t})^{-\frac12} \E_\rho(q_{j,t} \xi q_{j,t}) \E_\rho(q_{j,t})^{- \frac12} \in \rho_1(\M),$$$$\widehat{M}_{t}:  \NN_\rho \ni \xi \mapsto \E_\rho(q_{t})^{-\frac12} \E_\rho(q_{t} \xi q_{t}) \E_\rho(q_{t})^{- \frac12} \in \rho_1(\M).$$

\vskip5pt

\noindent \textbf{Analytic conditions:}

\begin{itemize}
\item[i)] \emph{Mean differences conditions} 

\vskip3pt

\begin{itemize}
\item[$\bullet$] $\displaystyle \widehat{R}_{j,t}(\xi^*\xi) - \widehat{R}_{j,t}(\xi)^*\widehat{R}_{j,t}(\xi) \, \le \, \Phi_{j,t} \big( \delta(\xi)^*\delta(\xi) \big)$,

\vskip5pt

\item[$\bullet$] $\displaystyle \big[\widehat{R}_{j,t}(\xi) - \widehat{M}_t(\xi) \big]^* \big[\widehat{R}_{j,t}(\xi) - \widehat{M}_t(\xi) \big] \, \le \, \Psi_{j,t} \big( \delta(\xi)^*\delta(\xi) \big)$,
\end{itemize}

\vskip5pt

\noindent for some derivation $\delta: \NN_\rho \to \NN_\sigma$ and cpu maps $\Phi_{j,t}, \Psi_{j,t}: \NN_\sigma \to \rho_1(\M)$. 

\vskip8pt

\item[ii)] \emph{Metric/measure growth conditions} 

\vskip3pt

\begin{itemize}
\item[$\bullet$] $\displaystyle \mathbf{1} \, \le \, \pi_1 \rho_1^{-1} \E_\rho(q_{t})^{-\frac12} \E_\pi (a_{t}^* a_{t}) \pi_1 \rho_1^{-1} \E_\rho(q_{t})^{-\frac12}  \, \lesssim \, \pi_1 \rho_1^{-1}(\gamma_{j,t}^2)$,

\vskip5pt

\item[$\bullet$] $\displaystyle \mathbf{1} \, \le \, \pi_1 \rho_1^{-1} \E_\rho(q_{j,t})^{-\frac12} \E_\pi (a_{j,t}^* a_{j,t}) \pi_1 \rho_1^{-1} \E_\rho(q_{j,t})^{-\frac12} \, \lesssim \, \pi_1 \rho_1^{-1}(\gamma_{j,t}^2)$,
\end{itemize}

\vskip5pt

\noindent for some family of operators $a_t, a_{j,t} \in \NN_\pi$ to be determined later on. 
\end{itemize}

\vskip3pt
 
\noindent A complete determination of the operators $a_t$ and $a_{j,t}$ is only possible after imposing additional size and smoothness conditions in our definition of Calder\'on-Zygmund operator below. Nevertheless, we shall see that these operators will play the role of \lq\lq dilated Markov balls" as it is the case in classical CZ theory. In fact, in the classical case our last condition trivially holds for doubling measures, and also for measures of polynomial or even exponential growth provided we find a Markov metric with large enough $\gamma_{j,t}$'s. Our assertions will be illustrated below. The first condition takes the form in the classical case of a couple of easy consequences of Jensen's inequality, namely
\begin{equation} \label{Eq-Jensen}
\begin{array}{rcl}
\displaystyle \mean_{\mathrm{B}_1} |f|^2 d\mu - \Big| \mean_{\mathrm{B}_1} f d\mu \Big|^2 & \le & \displaystyle  \mean_{\mathrm{B}_1 \times \mathrm{B}_1} \big| f(y) - f(z) \big|^2 d\mu(y) d\mu(z), \\
\displaystyle \Big| \ \mean_{\mathrm{B}_1} f d\mu \ - \ \mean_{\mathrm{B}_2} f d\mu \ \Big|^2 & \le & \displaystyle  \mean_{\mathrm{B}_1 \times \mathrm{B}_2} \big| f(y) - f(z) \big|^2 d\mu(y) d\mu(z).
\end{array}
\end{equation}

\subsection{CZ extrapolation} \label{SubsectCZE}

Now we introduce CZOs in this context. As we already mentioned, we consider a priori densely defined (unbounded) maps $T: \A_\M \to \M$ whose amplified maps are right $\mathcal{B}$-modules according to our algebraic assumptions above. In addition, we impose three conditions generalizing $L_2$ boundedness, the size and the smoothness conditions for the kernel. %For that we recall the operators $a_t, a_{j,t}$ introduced in the metric/measure growth condition and the norm $$\|f\|_{L_\infty^c(\NN; \E_\M)} \, = \, \Big\| \big( \E_\M(f^*f) \big)^\frac12 \Big\|_\M.$$ 

\vskip5pt

\noindent \textbf{Calder\'on-Zygmund type conditions:}
\begin{itemize}
\item[i)] \emph{Boundedness condition} $$\widehat{T}: L_\infty^c(\NN_\pi; \E_\pi) \to L_\infty^c(\NN_\rho; \E_\rho).$$

\vskip5pt

\item[ii)] \emph{Size \lq\lq kernel" condition} 

\vskip3pt

\begin{itemize}
\item[$\bullet$] $ \displaystyle \widehat{M}_{t} \Big( \big| \widehat{T} (\pi_2(f) (A_{j,t} - a_{t})) \big|^2 \Big) \, \lesssim \, \gamma_{j,t}^2 \|f\|_\infty^2$, 

\vskip3pt

\item[$\bullet$] $\displaystyle \widehat{R}_{j,t} \Big( \big| \widehat{T} (\pi_2(f) (A_{j,t} - a_{j,t})) \big|^2 \Big) \, \lesssim \,  \gamma_{j,t}^2 \|f\|_\infty^2$,
\end{itemize}

\vskip3pt

\noindent for a family of operators $A_{j,t} \in \NN_\pi$ with $A_{j,t} \ge a_{j,t}, a_t$ to be determined.

\vskip5pt

\item[iii)] \emph{Smoothness \lq\lq kernel" condition}

\vskip3pt

\begin{itemize}
\item[$\bullet$] $\displaystyle \Phi_{j,t} \Big( \big| \delta \big( \widehat{T}(\pi_2(f)(\mathbf{1} \hskip1.5pt - \hskip1.5pt a_{j,t})) \big) \big|^2 \Big) \, \lesssim \, \gamma_{j,t}^2 \|f\|_\infty^2$,

\vskip3pt

\item[$\bullet$] $\displaystyle \Psi_{j,t} \Big( \big| \delta \big( \widehat{T}(\pi_2(f)(\mathbf{1} - A_{j,t})) \big) \big|^2 \Big) \, \lesssim \, \gamma_{j,t}^2 \|f\|_\infty^2$.
\end{itemize}

%\vskip3pt

%\noindent This completes the conditions to be hold by the operators $a_t, a_{j,t}$ and $A_{j,t}$.
\end{itemize}

Let $T: \A_\M \to \M$ be a densely defined map which admits an amplification $\widehat{T}$ satisfying \eqref{Amplification}. Any such $T$ will be called an \emph{algebraic column} CZO whenever the amplification map is right $\mathcal{B}$-modular and satisfies the CZ conditions we have given above. At first sight, our boundedness assumption might appear to be unrelated to the classical condition. The reader could have expected the $L_2$ boundedness of $T$, but our assumption is formally equivalent to it in the classical case and gives the right condition for more general algebras. On the other hand, our size and smoothness conditions are intrinsic in the sense that the kernel is not specified under this degree of generality. We shall recover classical kernel type estimates from our conditions in our examples below. As explained above, the operators $a_t, a_{j,t}$ and $A_{j,t}$ play the role of dilated Markov balls and our conditions were somehow modeled by Tolsa's arguments in \cite{To}. Perhaps a significant difference ---in contrast to Tolsa's approach--- is that our smoothness condition is analog to a H\"ormander type condition, more than the (stronger) Lipschitz regularity assumption. 

\begin{theorem} \label{Extrapolation}
Let $(\M,\tau)$ be a noncommutative measure space equipped with a Markov semigroup $\S = (S_t)_{t \ge 0}$ with associated Markov metric $\Q$ which fulfills our algebraic and analytic assumptions. Then, any algebraic column \emph{CZO} $T$ defines a bounded operator $$T: \A_{\M} \to \mathrm{BMO}_\Q^c.$$
\end{theorem}

\dem The first goal is to estimate the norm of $$\mathrm{A} \, = \, \gamma_{j,t}^{-1} \Big( R_{j,t}|Tf|^2 - |R_{j,t}Tf|^2 \Big) \gamma_{j,t}^{-1}.$$ The map $\Pi_{j,t}: a \otimes b \in \M \bar\otimes_{R_{j,t}} \M \mapsto \mathbf{1} \otimes R_{j,t}(a)b \in \mathbf{1} \otimes \M$ extends to a right $(\mathbf{1} \otimes \M)$-module projection, which is well-defined in the sense that $\langle \xi, \xi \rangle_{R_{j,t}} = 0$ implies $\Pi_{j,t}(\xi) = 0$. Now, since $$\mathrm{A} \, = \, \gamma_{j,t}^{-1} \Big\langle Tf \otimes \mathbf{1} - \mathbf{1} \otimes R_{j,t}Tf \, , \, Tf \otimes \mathbf{1} - \mathbf{1} \otimes R_{j,t}Tf \Big\rangle_{R_{j,t}} \gamma_{j,t}^{-1},$$ we may use $\Pi_{j,t}$ to deduce the following identity
$$\mathrm{A} \, = \, \Big\langle (id - \Pi_{j,t}) (Tf \otimes \gamma_{j,t}^{-1}) \, , \, (id - \Pi_{j,t}) (Tf \otimes \gamma_{j,t}^{-1}) \Big\rangle_{R_{j,t}}.$$ Consider the amplification maps $\widehat{R}_{j,t}$ and $\widehat{\Pi}_{j,t}$ determined by $$R_{j,t} = \widehat{R}_{j,t} \circ \rho_2 \quad \mbox{and} \quad \Pi_{j,t} = \widehat{\Pi}_{j,t} \circ (\rho_2 \otimes id).$$ By \eqref{Amplification}, it turns out that $\mathrm{A} = \langle \mathbf{a}, \mathbf{a} \rangle_{\widehat{R}_{j,t}}$ where 
\begin{eqnarray*}
\mathbf{a} & = & (id - \widehat{\Pi}_{j,t}) (\rho_2 Tf \otimes \gamma_{j,t}^{-1}) \\ & = & (id - \widehat{\Pi}_{j,t}) (\widehat{T} \pi_2 f \otimes \gamma_{j,t}^{-1}) \\ & = & (id - \widehat{\Pi}_{j,t}) \big( \widehat{T} (\pi_2 (f) a_{j,t}) \otimes \gamma_{j,t}^{-1} \big) \\ & + & (id - \widehat{\Pi}_{j,t}) (\widehat{T} \big( \pi_2 (f) (\mathbf{1} - a_{j,t})) \otimes \gamma_{j,t}^{-1} \big) = \mathbf{a}_1 + \mathbf{a}_2
\end{eqnarray*}
According to Lemma \ref{ModuleLemma} i), we may estimate $\mathrm{A}$ as follows $$\mathrm{A} \, \lesssim \, \big\langle \mathbf{a}_1, \mathbf{a}_1 \big\rangle_{\widehat{R}_{j,t}} + \big\langle \mathbf{a}_2, \mathbf{a}_2 \big\rangle_{\widehat{R}_{j,t}} = \, \mathrm{A}_1 + \mathrm{A}_2.$$ Since $\widehat{\Pi}_{j,t} (n \otimes b) = \mathbf{1} \otimes \widehat{R}_{j,t}(n) b$, the Kadison-Schwarz inequality yields $$\Big\langle \widehat{\Pi}_{j,t} (n \otimes b), \widehat{\Pi}_{j,t} (n \otimes b) \Big\rangle_{\widehat{R}_{j,t}} \, \lesssim \, \big\langle n \otimes b, n \otimes b \big\rangle_{\widehat{R}_{j,t}}.$$ In conjunction with Lemma \ref{ModuleLemma} i) again, we deduce the following estimate for $\mathrm{A}_1$
$$\mathrm{A}_1 \! \lesssim \! \Big\langle \widehat{T} (\pi_2 (f) a_{j,t}) \otimes \gamma_{j,t}^{-1}, \widehat{T} (\pi_2 (f) a_{j,t}) \otimes \gamma_{j,t}^{-1} \Big\rangle_{\widehat{R}_{j,t}} \!\!\!\! = \gamma_{j,t}^{-1} \widehat{R}_{j,t} \Big( \big| \widehat{T} (\pi_2(f) a_{j,t}) \big|^2 \Big) \gamma_{j,t}^{-1}.$$ In order to bound the term in the right hand side, we apply \eqref{Initial} and the properties of the operator-valued weight $\mathsf{E}_\rho$ together with our algebraic conditions. Indeed, we first use the $\Q$-monotonicity of $\mathsf{E}_\rho$; then the fact that it commutes with the left/right multiplication by elements  affiliated to $\M$ (like $\gamma_{j,t}^{-1}$ or $\mathsf{E}_\rho(q_{j,t})^{-1/2}$); finally we use the right $\mathcal{B}$-modularity of the amplification of $T$:
\begin{eqnarray*}
\lefteqn{\hskip-15pt \gamma_{j,t}^{-1} \widehat{R}_{j,t} \Big( \big| \widehat{T} (\pi_2(f) a_{j,t}) \big|^2 \Big) \gamma_{j,t}^{-1}} \\ & \le & \gamma_{j,t}^{-1} \mathsf{E}_\rho(q_{j,t})^{-\frac12} \mathsf{E}_\rho \Big( \big| \widehat{T} (\pi_2(f) a_{j,t}) \big|^2 \Big) \mathsf{E}_\rho(q_{j,t})^{-\frac12} \gamma_{j,t}^{-1} \\ & = & \mathsf{E}_\rho \Big( \gamma_{j,t}^{-1} \mathsf{E}_\rho(q_{j,t})^{-\frac12} \big| \widehat{T} \big( \pi_2(f) a_{j,t} \big) \big|^2 \mathsf{E}_\rho(q_{j,t})^{-\frac12} \gamma_{j,t}^{-1} \Big) \\ & = & \mathsf{E}_\rho \, \Big| \widehat{T} \Big( \underbrace{\pi_2(f) a_{j,t} \, \pi_1 \rho_1^{-1} \big( \mathsf{E}_\rho(q_{j,t})^{-\frac12} \gamma_{j,t}^{-1} \big)}_{\xi_1} \Big) \Big|^2 \, = \, \mathsf{E}_\rho |\widehat{T}(\xi_1)|^2.
\end{eqnarray*}
Now, our first CZ condition i) gives the boundedness we need since 
\begin{eqnarray*}
\|\mathrm{A}_1\|_\M & \le & \big\| \widehat{T}(\xi_1) \big\|^2_{L_\infty^c(\NN_\rho; \mathsf{E}_\rho)} \ \lesssim \ \big\| \xi_1 \big\|^2_{L_\infty^c(\NN_\pi; \mathsf{E}_\pi)} \\ & = & \Big\| \mathsf{E}_\pi \Big( \big| \pi_2(f) a_{j,t} \, \pi_1 \rho_1^{-1} \big( \mathsf{E}_\rho(q_{j,t})^{-\frac12} \gamma_{j,t}^{-1} \big) \big|^2 \Big) \Big\|_\M \\ & \le & \Big\| \pi_1 \rho_1^{-1} \big( \gamma_{j,t}^{-1} \E_\rho(q_{j,t})^{-\frac12} \big)^* \E_\pi (a_{j,t}^* a_{j,t}) \pi_1 \rho_1^{-1} \big( \E_\rho(q_{j,t})^{-\frac12} \gamma_{j,t}^{-1} \big) \Big\|_\M \|f\|_\infty^2.
\end{eqnarray*}
The last term on the right is dominated by $\|f\|_\infty^2$ according to our second analytic condition on metric/measure growth. The estimate for $\mathrm{A}_2$ is simpler. Indeed, if we set $\xi_2 = \widehat{T} (\pi_2(f) (\mathbf{1} - a_{j,t}))$ then
\begin{eqnarray*}
\mathrm{A}_2 & = & \Big\langle (id - \widehat{\Pi}_{j,t}) (\xi_2 \otimes \gamma_{j,t}^{-1}) \, , \, (id - \widehat{\Pi}_{j,t}) (\xi_2 \otimes \gamma_{j,t}^{-1}) \Big\rangle_{\widehat{R}_{j,t}} \\ & = & \gamma_{j,t}^{-1} \Big( \widehat{R}_{j,t} |\xi_2|^2 - \big| \widehat{R}_{j,t} (\xi_2) \big|^2 \Big) \gamma_{j,t}^{-1} \ \le \ \gamma_{j,t}^{-1} \Phi_{j,t} \big( |\delta \xi_2|^2 \big) \gamma_{j,t}^{-1} \ \lesssim \ \|f\|_\infty^2 \mathbf{1},
\end{eqnarray*}
where the first inequality holds for some derivation $\delta: \NN_\rho \to \NN_\sigma$ and some cpu map $\Phi_{j,t}: \NN_\sigma \to \rho_1(\M)$ by our first analytic condition on mean differences. Then our CZ condition iii) on kernel smoothness justifies our last estimate. Our estimates so far prove the desired estimate $$\sup_{t > 0} \sup_{j \ge 1} \Big\| \Big( \gamma_{j,t}^{-1} \big[ R_{j,t}|Tf|^2 - |R_{j,t}Tf|^2 \big] \gamma_{j,t}^{-1} \Big)^{\frac12} \Big\|_\M \, \lesssim \, \|f\|_\infty.$$ Therefore, it remains to estimate the norm of $$\mathrm{B} \, = \, \gamma_{j,t}^{-1} \Big( \big| R_{j,t} Tf - M_tTf \big|^2 \Big) \gamma_{j,t}^{-1}.$$ To do so, we decompose the middle term using \eqref{Amplification} as follows
\begin{eqnarray*}
\lefteqn{\hskip-5pt R_{j,t} Tf - M_tTf} \\ [3pt] & = & \widehat{R}_{j,t} (\rho_2 T f) - \widehat{M}_t (\rho_2 T f) \\ [3pt] & = & \widehat{R}_{j,t} \Big( \widehat{T} \big( \pi_2 (f) a_{j,t} \big) \Big) - \widehat{M}_t \Big( \widehat{T} \big( \pi_2 (f) a_t \big) \Big) \\ & + & \Big[ \widehat{R}_{j,t} \Big( \widehat{T} \big( \pi_2 (f) (\mathbf{1} - a_{j,t}) \big) \Big) - \widehat{M}_t \Big( \widehat{T} \big( \pi_2 (f) (\mathbf{1} - a_t) \big) \Big) \Big] \ = \ \mathbf{b}_1 - \mathbf{b}_2 + \mathbf{b}_3.
\end{eqnarray*}
Letting $\mathrm{B}_j = \gamma_{j,t}^{-1} |\mathbf{b}_j|^2 \gamma_{j,t}^{-1}$ we get $\mathrm{B} \lesssim \mathrm{B}_1 + \mathrm{B}_2 + \mathrm{B}_3$. By Kadison-Schwarz we get $$\mathrm{B}_1 \, \le \, \gamma_{j,t}^{-1} \widehat{R}_{j,t} \Big( \big| \widehat{T} (\pi_2(f) a_{j,t}) \big|^2 \Big) \gamma_{j,t}^{-1} \, \lesssim \, \|f\|_\infty^2,$$ where the last inequality was justified in our estimate of $\mathrm{A}_1$ above. Replacing $q_{j,t}$ by $q_t$, the same argument serves to control the term $\mathrm{B}_2$. To estimate $\mathrm{B}_3$ we decompose $\mathbf{b}_3$ as follows 
\begin{eqnarray*}
\mathbf{b}_3 \!\!\! & = & \!\!\! \Big[ \widehat{R}_{j,t} \Big( \widehat{T} \big( \pi_2 (f) (\mathbf{1} - A_{j,t}) \big) \Big) - \widehat{M}_t \Big( \widehat{T} \big( \pi_2 (f) (\mathbf{1} - A_{j,t}) \big) \Big) \Big]  \\ \!\!\! & + & \!\!\! \hskip2pt \widehat{R}_{j,t} \Big( \widehat{T} \big( \pi_2 (f) (A_{j,t} - a_{j,t}) \big) \Big) - \widehat{M}_t \Big( \widehat{T} \big( \pi_2 (f) (A_{j,t} - a_t) \big) \Big) = \mathbf{b}_{31} + \mathbf{b}_{32} - \mathbf{b}_{33}.
\end{eqnarray*}
Taking $\xi_3 = \widehat{T} \big( \pi_2 (f) (\mathbf{1} - A_{j,t}) \big)$ and applying our analytic condition i) on mean differences together with our CZ condition iii) on kernel smoothness, we obtain that $$\gamma_{j,t}^{-1} |\mathbf{b}_{31}|^2 \gamma_{j,t}^{-1} \, \lesssim \, \gamma_{j,t}^{-1} \Psi_{j,t} \big( |\delta \xi_3|^2 \big) \gamma_{j,t}^{-1} \, \lesssim \, \|f\|_\infty^2.$$ It remains to estimate the terms $\mathrm{B}_{32}$ and $\mathrm{B}_{33}$. Applying the Kadison-Schwarz inequality, it is easily checked that these terms are also dominated by $\|f\|_\infty^2$ by means of our CZ size kernel condition ii). Altogether, we have justified that 
$$\sup_{t > 0}  \inf_{\stackrel{\null}{M_t}\mathrm{ cpu}}  \sup_{j \ge 1}\, \Big\| \Big( \gamma_{j,t}^{-1} \big[ |R_{j,t} f - M_t f|^2 \big] \gamma_{j,t}^{-1} \Big)^\frac12 \Big\|_\M \, \lesssim \, \|f\|_\infty.$$ Combining our estimates for $\mathrm{A}$ and $\mathrm{B}$, we deduce that $T: \A_\M \to \mathrm{BMO}_{\Q}^c$. \fin

The $\A_\M \to \mathrm{BMO}_\Q^r$ boundedness of the map $T$ is equivalent to the $\A_\M \to \mathrm{BMO}_\Q^c$ boundedness of the map $T^\dag(f) = T(f^*)^*$. According to this, an \emph{algebraic} CZO is any column CZO $T$ for which $T^\dag$ remains a column CZO. By Theorem \ref{Extrapolation}, we know that any algebraic CZO $T$ associated to $(\M,\tau,\Q)$ as above is automatically $\A_\M \to \mathrm{BMO}_\Q$ bounded. Assuming $L_2$ boundedness and regularity of the Markov semigroup, we may interpolate via Theorem \ref{MetricBMOInterp}. Under the same assumptions for $T^*$, we may also dualize and obtain the following extrapolation result. 

\begin{corollary} \label{Extrapolation2}
Let $(\M,\tau)$ be a noncommutative measure space equipped with a Markov regular semigroup $\S = (S_t)_{t \ge 0}$ and a Markov metric $\Q = (R_{j,t}, \sigma_{j,t}, \gamma_{j,t})$ fulfilling our algebraic and analytic assumptions. Then, every  $L_2$-bounded algebraic \emph{CZO} $T$ satisfies that $J_pT: L_p(\M) \to L_p^\circ(\M)$ for $p > 2$. Applying duality, similar conditions for $T^*$ yield $L_p$-boundedness of $TJ_p$ for every $1 < p < 2$.
\end{corollary}

\begin{rmk}
Theorem \ref{Extrapolation} admits a completely bounded version in the category of operator spaces. Since the operator space structure \cite{ER, P3} of $\mathrm{BMO}$ is determined by
$$M_m(\mathrm{BMO}_\T(\M)) = \mathrm{BMO}_{\widehat{\T}}(M_m(\M))$$ for $\widehat{\T}=(id_{M_m} \ten S_t)_{t\geq 0}$, we just need to replace $\M$ by $M_m(\M)$ everywhere, amplify all the involved maps by tensorizing with $id_{M_m}$ and require that the hypotheses hold with constants independent of $m$. Then, we obtain the cb-boundedness of $T$. 
\end{rmk}

\begin{remark}
\emph{As noticed in the Introduction, a common scenario is given by the choice $\NN_\rho = \M \bar\otimes \M$ with $\rho_1(f) = f \otimes \mathbf{1}$ and $\rho_2(f) = \mathbf{1} \otimes f$, together with $\mathsf{E}_\rho = id \otimes \tau$ and $\pi_j = \rho_j$ for $j=1,2$. In this case, it is clear that the amplification map is given by $$\widehat{T} = id_\M \otimes T \quad \mbox{so that} \quad \widehat{T} \pi_2 = \rho_2 T.$$ In particular, it turns out that the $L_2$ boundedness of $T$ in Corollary \ref{Extrapolation2} follows automatically from our CZ boundedness condition i). This is the case in classical Calder\'on-Zygmund theory. It is also true when $\NN_\rho = \M \bar\otimes \A$ for an auxiliary algebra $\A$ and $\rho_2 = \mathrm{flip} \circ \sigma$, where $\sigma: \M \to \A \bar\otimes \M$ is a  $*$-homomorphism satisfying $\mathsf{E}_\rho \circ \rho_2(f) = \tau(f) \1_\M$. This leads to another significant family of examples. It is however surprising that in general, the $L_2$ boundedness and the CZ boundedness assumptions are a priori unrelated. Thus, CZ extrapolation requires in this context to verify two boundedness conditions. It would be quite interesting to explore the corresponding \lq\lq $T(1)$ problems\rq\rq${}$ that arise naturally.}
\end{remark}

\subsection{The classical theory revisited} \label{ClassicalCZT}

We now illustrate our algebraic approach in the classical context of Euclidean spaces with the Lebesgue measure. This will help us to understand some of our conditions and will show how some others are automatic in a commutative framework. Take $\M = L_\infty(\R^n)$ with the Lebesgue measure and $\S = (H_t)_{t \ge 0}$ the heat semigroup $H_t = \exp(t \Delta)$. In Paragraph \ref{SectEuclideanMetric} we introduced the Markov metric $\Q$ given by 

\vskip8pt

\begin{itemize}
\item $\sigma_{j,t}^2 \equiv \frac{2e}{\sqrt{\pi}}  j^{\frac{n}{2}}e^{-k}$ and $\gamma_{j,t}^2 \equiv  j^{\frac{n}{2}}\ge 1$,

\vskip5pt

\item $R_{j,t}f(x) = \displaystyle \frac{1}{|\mathrm{B}_{\sqrt{4jt}}(x)|} \int_{\mathrm{B}_{\sqrt{4jt}}(x)} f(y) \, dy$.
\end{itemize}
Moreover, as we explained at the beginning of this section 
$$R_{j,t}f \, = \, \E_\rho(q_{j,t})^{-\frac12} \E_\rho \big( q_{j,t} \rho_2(f) q_{j,t} \big) \E_\rho(q_{j,t})^{- \frac12}$$ satisfies our basic assumption \eqref{Initial}. Here the amplification von Neumann algebra is $\NN_\rho = L_\infty(\R^n \times \R^n)$, the $*$-homomorphisms $\rho_j f (x_1,x_2) = f(x_j)$, the projections $q_{j,t}(x,y) = \chi_{\mathrm{B}_{\sqrt{4jt}}(x)}(y)$, and the operator-valued weight $\mathsf{E}_\rho$ is the integration map with respect to the second variable.  The cpu map $M_t$ which appears in the definition of $\mathrm{BMO}_\Q$  is still taken by $M_t=R_{1,t}$ as in Paragraph \ref{SectEuclideanMetric}.  

Taking $\NN_\pi = \NN_\rho$ and $T$ a standard CZO in $\R^n$, the algebraic conditions trivially hold in this case. Let $E_{j,k,t}^x = \mathrm{B}_{\sqrt{4jt}}(x) \times \mathrm{B}_{\sqrt{4kt}}(x)$. Taking $\Phi_{j,t}$ and $\Psi_{j,t}$ to be the averaging maps over $E_{j,j,t}^x$ and $E_{j,1,t}^x$ respectively and the family of dilated balls $(A_{j,t}(x,y), a_{j,t}(x,y)) = (\chi_{\alpha \mathrm{B}_{\sqrt{4jt}}(x)}(y), \chi_{5\mathrm{B}_{\sqrt{4jt}}(x)}(y))$ with $a\geq 5$, we may recover the   conditions as we explained right after stating them. Let us now show  how our algebraic CZ conditions hold from the classical ones. The boundedness condition reduces to the classical one, see Remark \ref{Multirem} A). Our size conditions can be rewritten as follows:
\begin{itemize}[leftmargin=1cm]
\item $\displaystyle \esssup \limits_{x\in \R^n}\mean_{\mathrm{B}_{\sqrt{4t}}(x)}\Big| \int_{a\mathrm{B}_{\sqrt{4jt}}(x)\setminus 5\mathrm{B}_{\sqrt{4t}}(x)}k(y,z)f(z)dz\Big|^2dy\lesssim j^{\frac{n}{2}}\| f\|^2_\infty$,

\item $\displaystyle \esssup \limits_{x\in \R^n}\mean_{\mathrm{B}_{\sqrt{4jt}}(x)}\Big| \int_{a\mathrm{B}_{\sqrt{4jt}}(x)\setminus 5\mathrm{B}_{\sqrt{4jt}}(x)}k(y,z)f(z)dz\Big|^2dy\lesssim j^{\frac{n}{2}}\| f\|^2_\infty$.
\end{itemize}
The above conditions follow from the usual size condition 
$$
|k(y,z)|\lesssim \frac{1}{|z-y|^n}.
$$ 
Next, taking $E_{j,k,t}^x$ as above, our smoothness conditions are:
\begin{itemize}[leftmargin=1cm]
\item $\displaystyle \esssup\limits_{x\in \R^n}\mean_{E_{j,1,t}^x}\Big( \int_{(5\mathrm{B}_{\sqrt{4jt}}(x))^c}\big( k(y_1,z)-k(y_2,z)\big)f(z)dz \Big)^2dy_1dy_2\lesssim j^{\frac{n}{2}}\| f\|^2_\infty$,

\item $\displaystyle \esssup\limits_{x\in \R^n}\mean_{E_{j,j,t}^x}\Big( \int_{(a\mathrm{B}_{\sqrt{4jt}}(x))^c}\big( k(y_1,z)-k(y_2,z)\big)f(z)dz \Big)^2dy_1dy_2\lesssim j^{\frac{n}{2}}\| f\|^2_\infty$.
\end{itemize}
The above conditions easily follow from the usual H\"ormander condition
$$\esssup\limits _{y_1,y_2\in \R^n}\int_{|y_1-z|\geq 2|y_1-y_2|} |k(y_1,z)-k(y_2,z)| \, dz < \infty.$$
Note that our algebraic CZ conditions are slightly weaker than the classical CZ conditions, but still sufficient to provide $L_\infty \to \mathrm{BMO}$ boundedness of the CZ map.

\begin{rmk}\label{rem: no size condition}
Our size condition is only used to estimate B in the proof of Theorem \ref{Extrapolation}. We saw in Paragraph \ref{SectEuclideanMetric} that $\mathrm{B}\lesssim \mathrm{A}$ for the Euclidean metric. Thus, our size condition is not necessary here, as it also happens in the classical formulation.
\end{rmk}

\section{{\bf Applications I | Commutative spaces}} \label{Sect-App1}

In this section we give specific constructions of Markov metrics on two basic commutative spaces: Riemannian manifolds with nonnegative Ricci curvature and Gaussian measure spaces. Beyond the Euclidean-Lebesguean setting considered above, these are the most relevant settings over which Calder\'on-Zygmund theory has been studied. As a good illustration of our algebraic method, we shall recover the extrapolation results. Noncommutative spaces will be explored later on.

\subsection{Riemannian manifolds}

Let $(\Omega,\mu)$ be a measure space equipped with a Markov semigroup, so that we may construct the corresponding semigroup type BMO space. In order to study the $L_\infty \to \mathrm{BMO}$ boundedness of CZOs in $(\Omega,\mu)$ it is essential to identify a Markov metric to work with. Now we provide sufficient conditions for a semigroup on a Riemannian manifold to yield a Markov metric satisfying our algebraic/analytic conditions, so that Theorem \ref{Extrapolation} is applicable. Let us consider an $n$-dimensional complete Riemannian manifold $(M,g)$ equipped with the geodesic distance $d$ determined by the Riemannian metric $g$. Denote the volume of a geodesic ball centered at $x$ with radius $r$ by $\mathrm{vol}_g(\mathrm{B}_r(x))$. Let $\T_M$ be a Markov semigroup on $M$ given by
$$S_{M,t}f(x) = \int_M s_t(x,y) f(y) \, dy.$$

\begin{proposition}
Assume that
\begin{itemize}
\item[\emph{i)}]$M$ has Ricci curvature $\ge 0$.

\vskip5pt

\item[\emph{ii)}] The kernel admits an upper bound 
$$
s_t(x,y) \lesssim  \frac{\phi(t)^{n+\varepsilon}}{\mathrm{vol}_g(\mathrm{B}_{\phi(t)}(x))(d(x,y) + \phi(t))^{n + \varepsilon}},
$$ 
for some strictly positive function $\phi$ and some parameter $\varepsilon > 0$.

\end{itemize}
Then $\T_M$ admits a Markov metric satisfying the algebraic/analytic conditions.
\end{proposition}

\dem If $\Sigma_{j,t}(x) = \mathrm{B}_{2^j \phi(t)}(x)$, our assumption. gives 
\begin{equation}\label{eq: estimate of kernel on R manifold}
s_t(x,y) \, \lesssim \, \sum_{j=1}^\infty \frac{2^{-j(n+\varepsilon)}}{\mathrm{vol}_g(\Sigma_{0,t}(x))} \chi_{\Sigma_{j,t}(x)}(y).
\end{equation}
According to Davies \cite[Theorem 5.5.1]{Dav2}, non-negative Ricci curvature implies
\begin{eqnarray*}
\mathrm{vol}_g(\mathrm{B}_r(x)) & \le & c_n r^n, \\ \mathrm{vol}_g(\mathrm{B}_{\gamma r}(x)) & \le & \gamma^n \mathrm{vol}_g(\mathrm{B}_{r}(x))
\end{eqnarray*}
for all $x \in M$, $r > 0$ and $\gamma > 1$. In particular, $\mathrm{vol}_g(\Sigma_{j,t}(x)) \le 2^{j n} \mathrm{vol}_g(\Sigma_{0,t}(x))$. By \eqref{eq: estimate of kernel on R manifold}, this implies that $(\sigma_{j,t}, \gamma_{j,t}) = (2^{-j \varepsilon/2}, 1)$ forms a Markov metric for $\T_M$ in conjunction with the averaging maps $$R_{j,t}f(x) = \mean_{\Sigma_{j,t}(x)} f(y) dy \quad \mbox{for} \quad (j,t)\in \mathbb{Z}_+\times \mathbb{R}_+.$$ 

Our construction for $\M = L_\infty(M)$ follows the basic model in the Introduction and the one used above in the Euclidean setting: $\Mn_\rho = \Mn_\pi = \M \bar\otimes \M$ with $\rho_j$ the canonical inclusion maps and $q_{j,t}(x,y) = \chi_{\Sigma_{j,t}(x)}(y) = \chi_{\Sigma_{j,t}(y)}(x)$. Then, the algebraic conditions for the  Markov metric  are obviously satisfied. Let us now check the analytic conditions. Taking $\Mn_\sigma = \M \bar\otimes \M \bar\otimes \M$, the derivation $\delta: \Mn_\rho \to \Mn_\sigma$ given by $\delta(a \otimes b) = a \otimes (\1 \otimes b - b \otimes \1)$ and the maps $M_t = R_{1,t}$, it turns out that the mean difference conditions follow from Jensen's inequality on normalized balls of $(M,g)$ as it follows from our comments after the definition of the analytic conditions. It remains to consider the metric/measure growth conditions. By taking $a_{j,t}(x,y)=\chi_{\Sigma_{j+1,t,}(x)}(y)$ and $(q_t,a_t) = (q_{1,t}, a_{1,t})$, these conditions reduce to show that $$\mathrm{vol}_g(\mathrm{B}_{2^{j+1}\phi(t)}(x)) \approx \mathrm{vol}_g(\mathrm{B}_{2^j\phi(t)}(x)).$$ This follows in turn from the fact that $M$ has a non-negative Ricci curvature. \fin

Let $(M,g)$ be a complete Riemannian manifold with non-negative Ricci curvature and let $\Delta$ be the
Laplace-Beltrami operator. The heat semigroup $\T_\Delta$ generated by $\Delta$ admits a kernel on $(M,g)$ satisfying the upper bound estimate mentioned in the above proposition. We know  from Davies \cite[Theorem 5.5.11]{Dav2} that the heat kernel satisfies
\begin{equation}\label{heat kernel upper bound}
h_t(x,y) \le
\frac{a_\delta}{\mathrm{vol}_g(\mathrm{B}_{\sqrt{t}}(x))} \exp
\Big( -\frac{d(x,y)^2}{4(1+\delta)t} \Big)
\end{equation}
for any $\delta > 0$
and certain constant $a_\delta$. This implies that
\begin{eqnarray*}
h_t(x,y) & \lesssim &  \frac{a_\delta}{\mathrm{vol}_g(\mathrm{B}_{\sqrt{t}}(x))} \frac{(4(1+\delta)t)^\frac{n+\varepsilon}{2}}{(d(x,y)^2+4(1+\delta)t)^\frac{n+\varepsilon}{2}}\\
& \lesssim & \frac{(\sqrt{4(1+\delta)t})^{n+\varepsilon}}{\mathrm{vol}_g(\mathrm{B}_{\sqrt{4(1+\delta)t}}(x))(d(x,y)+\sqrt{4(1+\delta)t} )^{n+\varepsilon}},
\end{eqnarray*}
which gives the expected upper bound with $\phi(t)=\sqrt{4(1+\delta)t}$. 

\begin{remark}
\emph{Once we have confirmed that algebraic and analytic conditions hold for the Markov process generated by the Laplace-Beltrami operator $\Delta$, it should be noticed that our CZ conditions are again implied by the classical ones. Arguing as in Remarks \ref{Rem-NDTerm} and \ref{rem: no size condition}, we see that the Ricci curvature assumption allows us to ignore our size kernel conditions. Next, it is straightforward to check that the boundedness condition reduces in this case to standard $L_2$-boundedness. Finally, our discussion in section \ref{ClassicalCZT} shows that our smoothness kernel condition is guaranteed under the classical H\"ormander kernel condition. Note in addition that our conditions also hold in the row case. In particular, classical CZOs in $(M,g)$ become algebraic CZOs. Moreover, the gaussian upper estimate \eqref{heat kernel upper bound} indicates that in $(M,g)$ with the heat semigroup $\T_{\Delta}$ we have $L_p^\circ(M)=L_p(M,g)$ for $1< p< \infty$.}
\end{remark}
By the discussion above, we have all the ingredients to apply Theorem \ref{Extrapolation} and Corollary \ref{Extrapolation2}. Let us illustrate it for the Riesz transforms on $(M,g)$. Consider the Riemannian gradient $\nabla = (\partial_1, \partial_2, \ldots, \partial_n)$ on $(M,g)$. The Riesz transform on $(M,g)$ is formally defined by 
$$
R = (R_j) = \nabla(-\Delta)^{-\frac{1}{2}} \quad \mbox{with} \quad R_j = \partial_j (-\Delta)^{-\frac{1}{2}}.
$$
Then we may recover Bakry's theorem \cite{B} using our algebraic approach. Indeed integration by parts gives $\| |\nabla f|\|_2 = \|\Delta^{\frac12} f\|_2$ which implies $L_2$-boundedness of Riesz transforms. Moreover, the H\"ormander condition follows from \cite{Chen,Li}. 

\begin{corollary}
Let $(M,g)$ be a complete $n$-dimensional Riemannian manifold with non-negative Ricci curvature. Then for all $1< p< \infty$, there exists a constant $C_p>0$ such that 
$$
\| R_j f \|_{L_p(M,g)} \le C_p \| f\|_{L_p(M,g)} \quad \mbox{for all} \quad 1 \le j \le n.
$$
\end{corollary}

\subsection{The Gaussian measure} 

Now we study the Ornstein-Uhlenbeck semigroup on the Euclidean space equipped with its Gaussian measure. We shall first construct a Markov metric for it. Then, we shall prove that our algebraic/analytic and Calder\'on-Zygmund conditions hold for the standard CZOs in this setting. The infinitesimal generator of the Ornstein-Uhlenbeck semigroup $O=(O_t)_{t\geq0}$ is the operator 
$$
L=\frac{\Delta}{2}-x\cdot \nabla
$$
on $(\R^n,\mu)$ with $d\mu(y)=\exp(-|y|^2)dy$. We have 
\begin{eqnarray*}
O_t f(x) & = &\frac{1}{(\pi -\pi e^{-2t})^\frac{n}{2}}\int_{\R^n} \exp \Big( -\frac{|e^{-t}x-y|^2}{1-e^{-2t}} \Big) f(y) \, dy \\
& = & \frac{1}{(\pi -\pi e^{-2t})^\frac{n}{2}}\int_{\R^n} \exp \Big( |x|^2-\frac{|e^{t}x-y|^2}{e^{2t}-1} \Big) f(y) \, d\mu(y).
\end{eqnarray*}

First, we establish a useful lemma showing that the local behavior |i.e. for small values of $t$| of the semigroup type BMO norm for the Ornstein-Uhlenbeck semigroup determines it completely. 

\begin{lemma}
Given $\delta > 0$, there exists $C_\delta > 0$ such that
$$
\sup_{t \ge 0} \big\| O_t |f|^2-|O_t f|^2 \big\|_\infty \le C_\delta \sup_{t < \delta} \big\| O_t |f|^2-|O_t f|^2 \big\|_\infty.
$$
\end{lemma}

\dem
It is easy to check that 
\begin{equation}\label{eq: expression of OU}
O_t f(x) = H_{v(t)}f (e^{-t}x),
\end{equation}
for $v(t)=\frac14 (1-e^{-2t})$ and the heat semigroup $H_t=\exp (t\Delta)$. Given $t>0$ and $f\in L_\infty(\R^n)$, let $F(s) = H_s |H_{t-s} f|^2$ for $0\leq s\leq t$. According to the definition of $H_t$, we obtain the following identity 
\begin{eqnarray*}
\partial_s F &= & (\partial_s H_s)|H_{t-s}f|^2+ H_s [(\partial_s H_{t-s}f)^*(H_{t-s}f)] +H_s [( H_{t-s}f)^*(\partial_s H_{t-s}f)] \\
& =& \Delta H_s  |H_{t-s}f|^2 -H_s [(\Delta H_{t-s}f)^*(H_{t-s}f)] -H_s [( H_{t-s}f)^*(\Delta H_{t-s}f)] \\
& =& H_s [ \Delta  |H_{t-s}f|^2 -(\Delta H_{t-s}f)^*(H_{t-s}f) - ( H_{t-s}f)^*(\Delta H_{t-s}f) ]\\
& = & 2H_s |\nabla H_{t-s} f |^2.
\end{eqnarray*}
Kadison-Schwarz inequality gives for $0<u<s$
$$
H_u |\nabla H_{t-u} f |^2 = H_u | H_{s-u}\nabla_x H_{t-s} f |^2 \leq H_s  | \nabla_x H_{t-s} f |^2
$$
which implies that $\partial_s F$ is increasing and $F$ is convex. Rearranging the inequality $F(s) \le \frac12 (F(0) + F(2s))$, we get $H_{2t}|f|^2-|H_{2t}f |^2\leq 2 H_t (H_{t}|f|^2-|H_{t}f |^2)$ for any $t\geq 0$. Then, the $L_\infty$ contractivity of $H_t$ gives
\begin{equation}\label{eq: compare of BMO T}
\big\| H_{2^k t}|f|^2-|H_{2^k t}f |^2 \big\|_\infty \leq 2^k \big\| H_{t}|f|^2-|H_{t}f |^2 \big\|_{\infty}.
\end{equation}
On the other hand, choosing $k_\delta$ such that $2^{k_\delta} v(\delta) \geq \frac{1}{4}$ and applying \eqref{eq: expression of OU} and 
\eqref{eq: compare of BMO T}
\begin{eqnarray*}
\hskip48pt \sup_{t \ge 0} \big\| O_t |f|^2-|O_t f|^2 \big\|_\infty & = & \sup_{t<\frac{1}{4}} \big\| H_{t}|f|^2-|H_{t}f |^2 \big\|_{\infty} \\ & \le & \sup_{t<2^{k_\delta} v(\delta)} \hskip-2pt \big\| H_{t}|f|^2-|H_{t}f |^2 \big\|_{\infty} \\ & \le & 2^{k_\delta} \hskip2pt \sup_{t< \delta} \hskip2.5pt \big\| O_t |f|^2-|O_t f|^2 \big\|_\infty. \hskip49pt \square
\end{eqnarray*}
%\null \vskip-145pt \null \fin

By the lemma above, it suffices to construct a Markov metric for $(O_t)_{t\geq 0}$ with $0<2t< \frac{1}{18}$. Let $v=\sqrt{e^{2t}-1}$ and consider the following family of balls and coronas in the gaussian space for $(j,t) \in \Z_+ \times \R_+$ 
$$ \Sigma_{j,t}(x) = \mathrm{B}(e^t x,\sqrt{j}v) \quad \mbox{and} \quad \Omega_{j,t}(x) = \Sigma_{j,t}(x) \setminus \Sigma_{j-1,t}(x).$$
Let $j_0 = j_0(x,t)$ be the smallest possible integer $j$ satisfying that $0 \in \Sigma_{j,t}(x)$.

\noindent \textbf{The case $n=1$.} If $1\leq  j<j_0$, let 
\begin{eqnarray*}
 D_{j,t}^-(x) &= & \Big\{ y \in \Omega_{j,t}(x) : e^t |x|-\sqrt{j}v \leq |y |\leq e^t |x|-\sqrt{j-1} v \Big\}, \\  D_{j,t}^+(x) & = & \Big\{ y \in \Omega_{j,t}(x) : e^t |x|+\sqrt{j-1}v \leq |y |\leq e^t |x|+ \sqrt{j} v \Big\}.
\end{eqnarray*}
Then, $D_{j,t}^-(x) \cup D_{j,t}^+(x) = \Omega_{j,t}(x)$ and we get 
\begin{equation} \label{eq: estimate of K}
O_tf(x)\lesssim \frac{1}{v} \Big( \sum_{\begin{subarray}{c} \varepsilon = \pm \\ 1 \le j<j_0 \end{subarray}} \exp (|x|^2-j)\int_{D_{j,t}^\varepsilon(x)} f d\mu + \sum_{j\geq j_0} \exp (|x|^2-j) \int_{\Sigma_{j,t}(x)} f d\mu \Big)
\end{equation}
for any positive $f \in L_\infty(\R,\mu)$. The above estimate indicates the natural candidates for the cpu maps $R_{j,t}$ and $\sigma_{j,t} \in L_\infty(\R,\mu)$. When $1\leq j<j_0$ and $\varepsilon = \pm$, we define 
$$R_{j,t,\varepsilon}f(x) = \frac{1}{\mu(D_{j,t}^\varepsilon(x))} \int_{D_{j,t}^\varepsilon(x)} f d\mu \quad \mbox{and} \quad \sigma_{j,t,\varepsilon}^2(x) = \frac{1}{v} \,\exp (|x|^2-j) \mu (D_{j,t}^\varepsilon(x)).$$ Note here we need an extra index for $R_{j,t}$'s when $j<j_0$. This is consistent with the assumptions (i) and (ii) in our definition of Markov metric, since we only need  the index-set of $R_{j,t}$'s to be countable. On the other hand, if $j\geq j_0$, we set
$$R_{j,t}f(x) = \frac{1}{\mu(\Sigma_{j,t}(x))}\int_{\Sigma_{j,t}(x)} f d\mu \quad \mbox{and} \quad \sigma_{j,t}^2(x) = \frac{1}{v}\,\exp (|x|^2-j) \mu (\Sigma_{j,t}(x)).$$
In order to find $\gamma_{j,t,\varepsilon}$ and $\gamma_{j,t}$ satisfying the metric integrability condition, we need to estimate $\mu (D_{j,t}^\varepsilon(x))$ and $\mu (\Sigma_{j,t}(x))$ respectively. Since the density function $\mu$ is monotone on $D_{j,t}^\varepsilon(x)$ we get 
\begin{eqnarray}\label{eq: estimate of D_k}
\mu(D_{j,t}^-(x)) & = & \int_{D_{j,t}^-(x)} e^{-y^2} dy \ \le \ \exp \big( - \big| e^t |x|-\sqrt{j} v \big|^2 \big) \frac{v}{\sqrt{j}}, \\
\mu(D_{j,t}^+(x)) & = & \int_{D_{j,t}^+(x)} e^{-y^2} dy \ \le \ \exp \big(- \big| e^t |x|+ \sqrt{j-1} v \big|^2 \big) \frac{v}{\sqrt{j}}.
\end{eqnarray}
When $j \ge j_0$ we use the trivial estimate
\begin{equation*}
\mu(\Sigma_{j,t}(x))  =  \int_{\Sigma_{j,t}(x)} e^{-y^2} dy \le 2 \sqrt{j}v.
\end{equation*}
Combining the estimates obtained above, we deduce for $1\leq j<j_0$
$$
\sigma_{j,t,-}^2\leq \frac{1}{\sqrt{j}} \, \exp \big( -\big| v |x|- e^t \sqrt{j} \big|^2 \big) \quad \mbox{and} \quad \sigma_{j,t,+}^2 \leq \frac{1}{\sqrt{j}} \, \exp \big( - \big| v|x|+  e^t \sqrt{j-1} \big|^2 \big).
$$
When $j\geq j_0$, we have $e^t |x|\leq \sqrt{j}v$ and $|x|^2\leq jv^2e^{-2t}< j/4$. Therefore
$$
\sigma_{j,t}^2 \leq 2\sqrt{j}\,\exp (|x|^2-j)<2\sqrt{j}\,\exp \Big( -\frac{3}{4}j \Big).
$$
Now we are ready to choose the optimal $\gamma$'s for the metric integrability condition in the definition of Markov metric. We respectively define for $1\leq j<j_0$ and $j \ge j_0$ 
$$
\gamma_{j,t,\varepsilon}^2(x) = \exp \Big( \frac{|v |x|+ e^t \sqrt{j-1} |^2}{2} \Big) \quad \text{} \mbox{and} \quad \gamma_{j,t}^2(x) =  \frac{1}{\sqrt{j}} \exp \Big( \frac{j}{4} \Big).
$$
Then it turns out that 
\begin{eqnarray*}
\sup_{\begin{subarray}{c} x \in \R \\ 0 < t < \frac{1}{36} \end{subarray}} \sum_{1\leq j<j_0 } \sigma_{j,t,-}^2(x) \gamma_{j,t,-}^2(x) & \le & \sup_{\begin{subarray}{c} x \in \R \\ 0 < t < \frac{1}{36} \end{subarray}} \sum_{1\leq j<j_0} \frac{1}{\sqrt{j}} \exp \Big( -\frac{| 6 v |x|-  e^t \sqrt{j} |^2}{4} \Big) \\ & \le &  \sup_{\begin{subarray}{c} x \in \R \\ 0 < t < \frac{1}{36} \end{subarray}} 2 \int_\R \exp \Big( -\frac{ | 6 v |x|- e^t u|^2}{4} \Big) \, du \ < \ \infty.
\end{eqnarray*}

\begin{eqnarray*}
\sup_{\begin{subarray}{c} x \in \R \\ 0 < t < \frac{1}{36} \end{subarray}} \sum_{1\leq j<j_0 } \sigma_{j,t,+}^2(x) \gamma_{j,t,+}^2(x) & \le & \sup_{\begin{subarray}{c} x \in \R \\ 0 < t < \frac{1}{36} \end{subarray}} \sum_{1\leq j<j_0} \frac{1}{\sqrt{j}} \exp \Big( -\frac{|v |x|+ e^t \sqrt{j-1} |^2}{2} \Big) \\ & \le & 1+  \sup_{\begin{subarray}{c} x \in \R \\ 0 < t < \frac{1}{36} \end{subarray}}2 \int_\R \exp \Big( -\frac{ |v |x|+e^t u|^2}{2} \Big) \, du \ < \ \infty.
\end{eqnarray*}
On the other hand, it is clear that
$$
\sup_{\begin{subarray}{c} x \in \R \\ 0 < t < \frac{1}{36} \end{subarray}}\sum_{j\geq j_0}\sigma_{j,t}^2(x) \gamma_{j,t}^2(x)  \, \le \,  \sum_{j\geq j_0} 2\, \exp(-\frac{j}{2}) \, < \, \infty.
$$
We have constructed a Markov metric for the Ornstein-Uhlenbeck semigroup. 

Let us now verify the analytic conditions, since the algebraic conditions are trivial by commutativity.  As we mentioned in Subsection \ref{Sect-AAConditions}, the first condition is an easy consequence of Jensen's inequality for the Gaussian measure. By the definition of $R_{j,t,\varepsilon}$ and $R_{j,t}$, we get $q_{j,t,\varepsilon}(x,y)= \chi_{D_{j,t}^\varepsilon(x)}(y)$ and $q_{j,t}(x,y)=\chi_{\Sigma_{j,t}(x)}(y)$. Thus, it remains to find proper $a_{j,t,\varepsilon}$ and $a_{j,t}$ to make sure that 
$$
\mu(D_{j,t}^\varepsilon(x)) \lesssim \int_\R a_{j,t,\varepsilon}^2(x,y) \, d\mu(y) \lesssim \gamma_{j,t,\varepsilon}^2(x) \mu(D_{j,t}^\varepsilon(x))
$$
and similarly for the pairs $(\Sigma_{j,t}(x)), a_{j,t}^2(x,y))$. When $j<j_0$ we consider the functions $a_{j,t,\varepsilon}(x,y) = \chi_{2D_{j,t}^\varepsilon(x)}(y)$, so the lower estimates are trivial. Denote by $c_{j,t,\varepsilon}(x)$ the center of $D_{j,t}^\varepsilon(x)$. Let  $\beta = v (\sqrt{j} - \sqrt{j-1})$. Arguing as in \eqref{eq: estimate of D_k}, we get
\begin{equation*}
\exp \Big( - \big| |c_{j,t,\varepsilon}(x)| + \frac{\beta}{2} \big|^2 \Big) \le \frac{\mu(D_{j,t}^\varepsilon(x))}{\beta} \le \exp \Big( - \big| |c_{j,t,\varepsilon}(x)| - \frac{\beta}{2} \big|^2 \Big),
\end{equation*}
\begin{equation*}
\exp \Big( - \big| |c_{j,t,\varepsilon}(x)| + \beta \big|^2 \Big) \le \frac{\mu(2D_{j,t}^\varepsilon(x))}{2 \beta} \le \exp \Big( - \big| |c_{j,t,\varepsilon}(x)| - \beta \big|^2 \Big).
\end{equation*}
Since $c_{j,t,\pm}(x) = e^t |x| \pm \frac12 v (\sqrt{j}+ \sqrt{j-1})$, this implies that 
$$ \frac{\mu(2D_{j,t}^\varepsilon(x))}{\mu(D_{j,t}^\varepsilon(x))} \, \le \, 2 \exp \Big( 3|c_{j,t,\varepsilon}(x)| \beta - \frac{3}{4} \beta^2 \Big) \, \le 2 \gamma_{j,t,\varepsilon}^2(x).$$
The estimate for $j \ge j_0$ is easier. Take $a_{j,t}(x,y)=\chi_{2\Sigma_{j,t}(x)}(y)$. Note
$$\begin{array}{rcccl}
2\sqrt{j}v\,\exp \big( -|2\sqrt{j}v|^2 \big) & \le & \mu (\Sigma_{j,t}(x)) & \le & 2\sqrt{j}v, \\
4\sqrt{j}v\,\exp \big(-|3\sqrt{j}v|^2 \big) & \le & \mu (2 \Sigma_{j,t}(x)) & \le & 4\sqrt{j}v.
\end{array}$$
Then, since $0<2t< \frac{1}{18}$,  we get
$$\frac{ \mu (2\Sigma_{j,t}(x)) }{ \mu (\Sigma_{j,t}(x)) } \leq 2\,\exp(4jv^2)\leq \, 2\exp (\frac{j}{4})=2\gamma_{j,t}^2(x).$$ Our choice for $q_t$ and $a_t$ correspond to the balls $\Sigma_{1,t}(x)$ and $2 \Sigma_{1,t}(x)$.

\noindent \textbf{The case $n>1$.} The argument is similar, so we  just point out the necessary modifications. Since $v (\sqrt{j} - \sqrt{j-1}) < v j^{-\frac12}$, we may pick $c_n j^{n-1}$ balls $D_{j,t}^s(x)$ for $1 \le s \le c_n j^{n-1}$ with radius $\frac{v}{2\sqrt{j}}$, centered on the sphere $$\Big\{ y : |y-e^tx |=\frac{v(\sqrt{j}+\sqrt{j-1})}{2} \Big\}$$ such that $$\Omega_{j,t}(x) \subset \bigcup _{s \ge 1} D_{j,t}^s(x)$$ and each $D_{j,t}^{s_0}(x)$ overlaps with at most $c_n'$ other balls $D_{j,t}^s(x)$. Then, if $f \ge 0$
$$
O_tf(x) \lesssim \frac{1}{v^n} \Big( \hskip-5pt \sum_{\begin{subarray}{c} j<j_0 \\ 1 \le s \le c_n j^{n-1} \end{subarray}} \exp ( |x|^2-j) \int_{D_{j,t}^s(x)} f d\mu + \sum_{j\geq j_0} \exp (|x|^2-j) \int_{\Sigma_{j,t}(x)} f d\mu \Big).
$$
Then, we may consider the following Markov metric 
\begin{eqnarray*}
\lefteqn{\hskip-10pt \Big( \sigma_{j,t,s}^2(x), \gamma_{j,t,s}^2(x), R_{j,t,s}f(x) \Big)} \\ & = & \Big(\frac{\exp (|x|^2-j) }{v^n} \mu (D_{j,t}^s(x)), j^{-\frac{n-1}{2}} \exp \Big( \frac{|v |x|+ e^t \sqrt{j-1} |^2}{2} \Big), \mean_{D_{j,t}^s(x)} f d\mu \Big)
\end{eqnarray*}
for $j < j_0$ and $1 \le s \le c_n j^{n-1}$. When $j \ge j_0$, we set
$$\Big( \sigma_{j,t}^2(x), \gamma_{j,t}^2(x), R_{j,t,}f(x) \Big) \, = \, \Big( \frac{\exp (|x|^2-j)}{v^n} \mu (\Sigma_{j,t}(x)), j^{-\frac{n}{2}}\exp(\frac{j}{4}), \mean_{\Sigma_{j,t}(x)} f d\mu \Big).$$
The analytic conditions hold under the same choices we made for $n=1$.

\begin{corollary}\label{thm: CZ Gaussian}
Let $\mathcal{O}=(O_t)_{t\geq 0}$ be the Ornstein-Uhlenbeck semigroup and $T$ be a singular integral operator defined on $L_\infty(\R^n,d\mu)$ with kernel $k$. More precisely, we have the kernel representation 
$$
Tf(x) \, = \, \int_{\R^n} k(x,y)f(y) \, d\mu(y) \quad \mbox{for} \quad x \notin \mathrm{supp} f.
$$
Suppose $T$ is bounded on $L_2(\R^n,\mu)$  and it satisfies 
\begin{equation}\label{eq: size k O}
\sup_{\mathrm{B \ ball}} \, \sup_{\begin{subarray}{c} z \in \mathrm{B} \\ j \ge 1\end{subarray}} \ \int_{2^{j+1} \mathrm{B}\setminus 2^j \mathrm{B}} | k(z,y)|d\mu(y)<\infty,
\end{equation}
\begin{equation}\label{eq: smooth k O}
\sup_{\mathrm{B \ ball}} \, \sup_{z_1, z_2 \in \mathrm{B}} \ \int_{(5\mathrm{B})^c} | k(z_1,y)-  k(z_2,y) | d\mu(y)<\infty.
\end{equation} 
Then $T$ is a bounded map from $L_\infty(\R^n,\mu)$ to the semigroup $\mathrm{BMO}_{\mathcal{O}}$ space.
\end{corollary}

\dem It suffices to prove that our CZ conditions hold. The row and column boundedness conditions reduce to $L_2$-boundedness. Let $M_t$ be the averaging map in $L_\infty(\R^n,\mu)$ over the ball $\Sigma_{1,t}(x)$. 
Given $1\leq  j<j_0$, define $A_{j,t,s}(x,\cdot)$ as the characteristic function over the ball $\Sigma_{j+1,t}(x)$. Then $A_{j,t,s} \leq  \chi_{2^{r} D_{j,t}^s(x) \cap \Sigma_{1,t}(x)}$ with $r= [2 \log_2(j+1)]+3$, where $[ \ ]$ stands for the integer part. Therefore, applying \eqref{eq: size k O}, we have 
\begin{eqnarray*}
 \sup_{z \in \Sigma_{1,t}(x)} \int_{\Sigma_{j+1,t}(x) \setminus 2 \Sigma_{1,t}(x)} |k(z,y)| \, d\mu(y) \!\! & \lesssim & \!\! r  \ \ \lesssim \ \gamma_{j,t,s}^2(x), \\
 \sup_{z \in D_{j,t}^s(x)} \int_{\Sigma_{j+1,t}(x) \setminus 2D_{j,t}^s(x)}  |k(z,y)| \, d\mu(y) \!\! & \lesssim & \!\! r  \ \ \lesssim \ \gamma_{j,t,s}^2(x). 
\end{eqnarray*}
For $j\geq j_0$, let  $A_{j,t}(x,\cdot) = \chi_{\Sigma_{4j,t}(x)} \le \chi_{2^u \Sigma_{1,t}(x)}$ with $u= [\log_2 (2\sqrt{j})]+1$. Applying \eqref{eq: size k O} as above, we see that $T$ satisfies our size conditions. Moreover,  \eqref{eq: smooth k O} implies our smoothness conditions as in the Euclidean-Lebesguean setting, Section \ref{ClassicalCZT}. \fin

\begin{rmk}
Since the Gaussian measure is non-doubling, the term $R_{j,t} f-M_t f$ in the Markov metric BMO space $\mathrm{BMO}_{\mathcal{Q}}$ is essential to characterise the changes of the mean values of the function $f$. This explains the relevance of the size kernel condition in the Calder\'on-Zygmund theory for the gaussian measure.
\end{rmk}

\section{{\bf Applications II | Noncommutative spaces}} \label{Sect-App2}

In this section we apply our algebraic approach to study Calder\'on-Zygmund operators in flag von Neumann algebras which originally motivated us and include matrix algebras, quantum Euclidean spaces and quantum groups. We start by reconstructing and refining the semicommutative theory, which deals with  tensor and crossed products with metric measure spaces. 

\subsection{Operator-valued theory}\label{Subsec-Opvalued app}

Let $(\Omega,\mu)$ be a doubling metric space |as in Remark \ref{Rem-NDTerm}| and consider a Markov semigroup $S_t: L_\infty(\Omega) \to L_\infty(\Omega)$. Let $\M$ be a semifinite von Neumann algebra with a n.s.f. trace $\tau$. Then we call the semigroup $\T = (S_t \otimes id_\M)_{t \ge 0}$ a \emph{semicommutative Markov semigroup}. Consider the algebra of essentially bounded functions $f: \Omega \to \M$ equipped with the trace $$\varphi(f) = \int_{\Omega} \tau(f(y)) \, d\mu(y).$$ Its weak-$*$ closure $\RR = L_\infty(\Omega) \bar\otimes \M$ is a von Neumann algebra. Assume that there exists a Markov metric  $\Q = \{ ( R_{j,t}, \sigma_{j,t}, \gamma_{j,t} ) : (j,t) \in \Z_+ \times \R_+ \}$ associated to the original Markov semigroup on $L_\infty(\Omega)$. Let $q_{j,t}(x,y) = \chi_{\Omega_{j,t}^x}(y)$ stand for the projections determined by $\Q$ via \eqref{Initial}. We assume in addition that $\Q$ satisfies the metric/measure growth condition  
\begin{equation}\label{eq: metric growth}
\frac{\mu(\Sigma_{j,t}^x)}{\mu(\Omega_{j,t}^x)}\leq \gamma_{j,t}(x)
\end{equation}
by choosing $a_{j,t}(x,y)=\chi _{\Sigma_{j,t}^x}(y)$. The remaining algebraic and analytic conditions trivially hold in this case. Indeed, the algebraic conditions follow by commutativity and analytic conditions just require to pick the right averaging maps according to Jensen's inequality, as explained in \eqref{Eq-Jensen}. 
Note that $\Q$ satisfies an operator-valued generalization of the Hilbert module majorization in the line of Remark \ref{ossMarkovBMO}. Thus $\Q$ extends to a Markov metric in $\RR$ by tensorizing with $id_\M$ and $\1_\M$ respectively.

\noindent Our goal is to study CZO's formally given by $$Tf(x) = \int_\Omega k(x,y) \, (f(y)) \, d\mu(y) \quad \mbox{with} \quad \begin{cases} f: \Omega \to \M_1 \ \mbox{and} \ x \notin \mathrm{supp}_\Omega f, \\ k(x,y) \in \mathcal{L}(L_0(\M_1), L_0(\M_2)). \end{cases}$$ That is, $k(x,y)$ is linear from $\tau_1$-measurable to $\tau_2$-measurable operators. If we set $\RR_j = L_\infty(\Omega) \bar\otimes \M_j$, we should emphasize that $L_p(\RR_j) = L_p(\Omega; L_p(\M_j))$. In particular, this framework does not fall in the vector-valued theory because we take values in different Banach spaces for different values of $p$, see \cite{Pa1} for further explanations. This class of operators is inspired by two distinguished examples with $\M_1 = \M = \M_2$:
\begin{itemize}
\item Operator-valued case $$Tf(x) = \int_\Omega k_{ov}(x,y) \cdot
\! f(y) \ d\mu(y).$$

\item Noncommutative model $$Tf(x) = \int_\Omega (id_\M \otimes
\tau) \, \big[ k_{nc}(x,y) \cdot \! (\mathbf{1}_\M \otimes f(y)) \big] \,
d\mu(y).$$
\end{itemize}
In the first case, the kernel takes values in $\M$ or even in the
complex field and acts on $f(y)$ by left multiplication
$k(x,y)(f(y)) = k_{ov}(x,y) \cdot f(y)$. It is the canonical map
when $L_p(\RR)$ is regarded as the Bochner space
$L_p(\Omega;L_p(\M))$. On the contrary if we simply think of
$L_p(\RR)$ as a noncommutative $L_p$ space, a \emph{natural} CZO
should be an integral map with respect to the full trace $\varphi
= \int_\Omega \otimes \tau$ and the kernel should be a $\varphi
\otimes \varphi$-measurable operator $k: \Omega \times
\Omega \to \M \bar\otimes \M$. The noncommutative model provides the
resulting integral formula. Note that this model also falls in our
general framework by taking $k(x,y)(f(y)) = (id_\M \otimes \tau) [
k_{nc}(x,y) \cdot (\1_\M \otimes f(y)) ]$.

\begin{theorem} \label{CZ2}
Let $\T = (S_t)_{t \ge 0}$ be a Markov
semigroup on $(\Omega,\mu)$ which admits a Markov metric $\Q = \{ ( R_{j,t}, \sigma_{j,t}, \gamma_{j,t} ) : (j,t) \in \Z_+ \times \R_+ \}$ satisfying the above assumptions. Let $a_{j,t}(x,y) = \chi_{\Sigma_{j,t}^x(y)}$ be the projections determined by $\Q$ via \eqref{eq: metric growth}.
Consider the \emph{CZO} formally given by $$Tf(x) =
\int_\Omega k(x,y) (f(y)) \, d\mu(y).$$ Then, $T$ maps
$L_\infty(\RR_1)$ to $\mathrm{BMO}_\T^c(\RR_2)$ provided the
conditions below hold
\begin{itemize}
\item[\emph{i)}] $L_2^c$-boundedness condition, $$\Big\| \Big( \int_\Omega
|Tf|^2 \, d\mu \Big)^\frac12 \Big\|_{\M_2} \ \lesssim \,
\Big\| \Big( \int_\Omega |f|^2 \, d\mu \Big)^\frac12
\Big\|_{\M_1}.$$

\vskip3pt

\item[\emph{ii)}] Smoothness condition for the kernel, $$\hskip25pt
\int_{\Omega \setminus \Sigma_{j,t}^x} \big\| \big( k(y_1,z) - k(y_2,z) \big) \big( f(z)  \big) \big\|_{\M_2}
\, d\mu(z) \lesssim \|f\|_{\RR_1}$$ uniformly in $j \ge 1$, $t > 0$, $x \in \Omega$ and $y_1, y_2 \in \Omega_{j,t}^x$.
\end{itemize}
\end{theorem}

\dem The proof follows from Theorem \ref{Extrapolation}. Since the underlying space $(\Omega,\mu)$ is a doubling metric space, the size  kernel condition is unnecessary. Thus, it remains to check the $L_2^c$-boundedness condition and the kernel smoothness condition.  Consider $\Mn_\pi=L_\infty(\Omega\times \Omega)\bar\otimes \M_1$, $\Mn_\rho=L_\infty(\Omega\times \Omega)\bar\otimes \M_2$, $\omega
(\varphi)(x,y)=\varphi(y)$ for $\varphi \in L_\infty(\Omega)$ and $(\pi_2, \rho_2) = (\omega \otimes id_{\M_1},\omega \otimes id_{\M_2})$. Let $\widehat{T}=id_{\Omega}\ten T$, $\Phi_{j,t}$ be  the averaging map  over $\Omega_{j,t}^x \times \Omega_{j,t}^x$ and $\Delta =\delta \ten id_{\M_2}$ with $\delta \varphi (x,y) = \varphi(x) - \varphi(y)$. Then condition i) yields the $L_2^c$-boundedness condition. It is also easy to see that condition ii) implies our kernel smoothness condition. 
Thus, the result follows from Theorem \ref{Extrapolation}. \fin

\begin{rmk} \label{Multirem}
We continue with a few comments:

\textbf{A)} When $\M_1 = \M_2=\M$ and the kernel $k(x,y) (f(y)) = k(x,y) \cdot f(y)$ acts by left multiplication, the boundedness condition i) becomes equivalent to the usual $L_2$ boundedness. Indeed, using that $\M \subset \mathcal{B}(L_2(\M))$
we obtain
\begin{eqnarray*}
\lefteqn{\hskip-10pt\Big\| \Big( \int_{\R^n} |Tf(y)|^2 \, dy
\Big)^\frac12 \Big\|_\M} \\ & = & \sup_{\|h\|_2 \le 1} \Big(
\int_{\R^n} \big\langle h, |Tf(y)|^2 h \big\rangle \, dy
\Big)^{\frac12} \ = \ \sup_{\|h\| \le 1} \big\| (Tf) \, (\1_{\R^n}
\otimes h) \big\|_{L_2(\RR)} \\ [7pt] & = & \sup_{\|h\| \le 1}
\big\| T(f \, (\1_{\R^n} \otimes h)) \big\|_{L_2(\RR)} \ \le \
\|T\|_{\mathcal{B}(L_2(\RR))} \, \Big\| \Big( \int_{\R^n}
|f(y)|^2 \, dy \Big)^{\frac12} \Big\|_\M.
\end{eqnarray*}

\textbf{B)} We have used so far semigroup type $\mathrm{BMO}$'s. When $(\Omega,\mu)$ comes equipped with a doubling metric, we may replace it by other standard (equivalent) forms of BMO, as pointed in Remark \ref{Remark-EquivBMO}. By well-known arguments \cite{Pa1}, our kernel smoothness condition reduces to 
\renewcommand{\theequation}{$\mathrm{Sm}_\lambda$}
\addtocounter{equation}{-1}
\begin{equation} \label{SmoothSC}
\displaystyle \sup_{R > 0}
\esssup_{y_1, y_2 \in
\mathrm{B}_R}\Big\| \int_{(\mathrm{B}_{\lambda R})^c} \big( k(y_1,z) - k(y_2,z)
\big)( f(z)) \, dz \Big\|_\M \lesssim \|f\|_\RR.
\end{equation}
for $\lambda>1$. The classical H\"ormander condition 
\renewcommand{\theequation}{$\mathrm{Hr}_\lambda$}
\addtocounter{equation}{-1}
\begin{equation} \label{HorSC}
\esssup_{y_1,y_2\in \R^n} \int_{d(y_1,z) > \lambda d(y_1,y_2)} \big\|
k(y_1,z) - k(y_2,z) \big\|_\M\, dz <\infty.
\end{equation}
satisfies $(\mathrm{Hr}_\lambda) \Rightarrow
(\mathrm{Sm}_{2\lambda+1})$. In fact, an even weaker condition suffices $$\displaystyle \sup_{R > 0}\Big\|\mean_{\mathrm{B}_R\times\mathrm{B}_R} \Big|\int_{(\mathrm{B}_{\lambda R})^c}\big( k(y_1,z)-k(y_2,z)\big)f(z)dz \Big|^2 dy_1dy_2\Big\|_\M \lesssim \| f\|^2_\RR.$$

\textbf{C)} We recall that $L_\infty(\RR) \to \mathrm{BMO}_\T$ boundedness requires that $T^\dag f = T(f^*)^*$ satisfies the same assumptions as $T$. If $k(x,y) \in \M$ is given by left multiplication the only effect in $T^\dag$ is that $k(x,y)$ is replaced by $k(x,y)^*$ and now operates by right multiplication. This left/right condition was formulated in \cite{Pa1} in terms of $\M$-bimodular maps. Moreover, a counterexample was constructed to show that the bimodularity is indeed essential. It is also quite interesting to note that in the \lq noncommutative model\rq${}$ we have $$\int_{\R^n} (id_\M \otimes \tau) \, \big[ k(x,y) \cdot \! (\1_\M \otimes f(y)) \big] \, dy = \int_{\R^n} (id_\M \otimes \tau) \, \big[ (\1_\M \otimes f(y)) \cdot k(x,y) \big] \, dy$$ by traciality and this pathology does not occur. Finally, the $L_p$ boundedness is guaranteed for $2 < p < \infty$ since the classical heat semigroup has a regular Markov metric and $J_p = id_{L_p(\R^n)}$ in this case. As for $1 < p < 2$, it suffices to take adjoints which leads to H\"ormander smoothness in the second variable $$\esssup_{z_1,z_2\in \R^n} \int_{|y-z_1| > \lambda |z_1-z_2|} \big\| k(y,z_1) - k(y,z_2) \big\|_\M \, dy <\infty.$$ Of course, this is still consistent with the classical CZ theory $\M = \C$.

\textbf{D)}  Our analysis of the semicommutative case from our basic Theorem
\ref{CZ2} does not recover the weak type $(1,1)$ inequality from \cite{Pa1}. It requires quasi-orthogonality methods which are still missing for general von Neumann algebras.
\end{rmk}

We now study the $L_\infty \to \mathrm{BMO}$ boundedness of twisted CZO's on homogeneous spaces. Given a discrete group $\G$ with left regular representation $\lambda :\G\mapsto \mathcal{B}(\ell_2(\G))$ let $\mathcal{L}(\G)$ denote its group von Neumann algebra. Let $(\M,\tau)$ with $\M \subset \mathcal{B}(\mathcal{H})$ be a noncommutative probability space and $\alpha: \G \mapsto \mathrm{Aut}(\M)$ be a trace preserving action. Consider two $*$-representations 
$$
\rho: \M \ni f \mapsto \sum_{h\in \G} \alpha_{h^{-1}}(f)\otimes e_{h,h} \in \M \bar{\otimes}\mathcal{B}(\ell_2(\G)),$$ $$\Lambda: \G \ni g \mapsto \sum_{h\in \G} \1_\M \otimes e_{gh,h} \in \M \bar{\otimes}\mathcal{B}(\ell_2(\G)),$$
where $e_{g,h}$ is the matrix unit for $\mathcal{B}(\ell_2(\G))$. Now we define the crossed product algebra $\M \rtimes_\alpha \G$ as the weak operator closure in $\M \otimes  \mathcal{B}(\ell_2(\G))$ of the $*$-algebra generated by $\rho(\M)$ and $\Lambda(\G)$. A generic element of $\M  \rtimes_\alpha \G$ can be formally written as $\sum_{g\in \G} f_g  \rtimes_\alpha \lambda (g)$ with $f_g\in \M$. With this convention, we may embed the crossed product algebra $\M  \rtimes_\alpha \G$ into $\M \bar{\otimes}\mathcal{B}(\ell_2(\G))$ via the map $j=\rho \rtimes \Lambda$. Indeed, we have
\begin{eqnarray*}
j \Big( \sum_{g\in \G} f_g  \rtimes_\alpha \lambda (g) \Big) & = & \sum_{g\in \G}  \rho(f_g)\Lambda(g) \\
& =& \sum_{g\in \G} \Big(\sum_{h,h'\in \G} (\alpha_{h^{-1}}(f_g)\otimes e_{h,h})(\1_{\M}\otimes e_{gh',h'}) \Big)\\
&=& \sum_{g\in \G} \Big(\sum_{h\in \G} \alpha_{h^{-1}}(f_g)\otimes e_{h,g^{-1}h} \Big)\\
&=& \sum_{g\in \G} \Big(\sum_{h\in \G} \alpha_{(gh)^{-1}}(f_g)\otimes e_{gh,h} \Big).
\end{eqnarray*}

Since the action $\alpha$ will be fixed, we relax the terminology and write $\sum_{g\in \G}f_g \lambda(g)$ instead of $\sum_{g\in \G}f_g \rtimes _\alpha\lambda(g)$. We say that a Markov semigroup $\T = (S_t)_{t \ge 0}$ in $\M$ is \emph{$\G$-equivariant} if 
$$
\alpha _g S_t =S_t \alpha_g \;\; \text{for} \;\; (t,g)\in \mathbb{R}_+\times \G.
$$
If $\T$ is a $\G$-equivariant Markov semigroup on $\M$, let $\T_\rtimes = (S_t \rtimes id_\G)_{t\geq 0}$ and $\T_\otimes =(S_t\otimes id_{\mathcal{B}(\ell_2(\G))})_{t\geq 0}$ denote the crossed/tensor product amplification of our semigroup on $\M \rtimes \G$ and $\M\bar\otimes \mathcal{B}(\ell_2(\G))$ respectively. Note that $\T_\rtimes$ is Markovian due to the $\G$-equivariance of $\T$. In the following result, our CZO's are of the form $$Tf(x) = \int_\Omega k(x,y)( f(y)) \, d\mu(y)$$ for all $f \in (\RR_1, \varphi_1)$, where $(\RR_j, \varphi_j) = L_\infty(\Omega,\mu) \bar\otimes (\M_j, \tau_j)$ and $k(x,y): \M_1 \to \M_2$. In other words, we keep the same terminology as for Theorem \ref{CZ2}. We shall also use the notation $$\widehat{\M}_j = \M_j \bar\otimes \mathcal{B}(\ell_2(\G)) \quad \mbox{and} \quad  \widehat{\RR}_j = \RR_j \bar\otimes \mathcal{B}(\ell_2(\G)).$$

\begin{corollary} \label{PrNonequiv}
Let $\G \curvearrowright L_\infty(\Omega,\mu)$ be an action $\alpha$ which is implemented by a measure preserving
transformation $\beta$, so that $\alpha_gf(x) = f (\beta_{g^{-1}}x)$. Let $\T = (S_t)_{t \ge 0}$ be a $\G$-equivariant Markov
semigroup on $(\Omega,\mu)$ which admits a Markov metric $\Q = \{ ( R_{j,t}, \sigma_{j,t}, \gamma_{j,t} ) : (j,t) \in \Z_+ \times \R_+ \}$ satisfying the assumptions above. Let us consider a family of \emph{CZO's} formally given by $$T_gf(x) = \int_\Omega k_g(x,y)( f(y)) \, d\mu(y) \quad \mbox{for} \quad g \in \G.$$ Then, $\summ_g f_g \lambda(g) \mapsto \summ_g T_g(f_g) \lambda(g)$ is bounded $\RR_1 \rtimes \G \to \mathrm{BMO}_{\T_\rtimes}^c(\RR_2 \rtimes \G)$ if 
\begin{itemize}
\item[\emph{i)}] $L_2^c$-boundedness condition, $$\Big\| \Big( \int_\Omega
\big| (T_{gh^{-1}}) \bullet \xi \big|^2 \, d\mu \Big)^\frac12
\Big\|_{\widehat{\M}_2} \ \lesssim  \Big\| \Big(
\int_\Omega |\xi|^2 \, d\mu \Big)^\frac12
\Big\|_{\widehat{\M}_1},$$
where $\bullet $ stands for the generalized Schur product of matrices. In other words, the \emph{CZO} $T_{gh^{-1}}$ only acts on the $(g,h)$-th entry of $\xi$ for each $g,h \in \G$.

\vskip3pt

\item[\emph{ii)}] Smoothness condition for the kernel, $$\hskip35pt \int_{\Omega} \big\| \big( K(y_1,z) - K(y_2,z) \big) \bullet \big( \xi(z) (\1 - a_{j,t}(x,z)) \big) \big\|_{\widehat{\M}_2} \, d\mu(z) \lesssim \|\xi\|_{\widehat{\RR}_1},$$ 
uniformly on $j \ge 1$, $t > 0$, $x \in \Omega$ and $y_1, y_2 \in \Omega_{j,t}^x$. Here, the \emph{CZ} kernel $K(y,z) = \summ_{g,h} k_{gh^{-1}}(\beta_gy, \beta_gz) \otimes e_{g,h}$ acts once more as a Schur multiplier.
\end{itemize}
\end{corollary}

\dem 
Letting $\xi = \sum_{g,h} a_{g,h} \otimes e_{g,h} \in \widehat{\RR}_1$, we define the map  
$$\Phi: \widehat{\RR}_1 \to \mathrm{BMO}_{\T_\otimes}^c(\widehat{\RR}_2),$$
$$
\Phi(\xi)(x) = \sum_{g,h}
\alpha_{g^{-1}} \int_\Omega k_{gh^{-1}}(x,y)
(a_{g,h}(\beta_{g^{-1}}(y))) \hskip1pt d\mu(y) \otimes e_{g,h}. 
$$
By the definition of $j$, it is easy to check that $$j \Big( \sum_g T_g(f_g) \lambda(g) \Big) \ = \ \Phi
\Big( j \big( \sum_g f_g \lambda(g) \big) \Big).$$ 
 Since $\T$ is $\G$-equivariant, according to \cite[Lemma 2.1]{JMP4}, we have 
 $$
 \| g\|_{\rm{BMO}_{\T\rtimes}^c(\RR_2\rtimes G)}=\sup_{t\geq 0} \Big\|\Big(S_{\otimes,t} |j(g)|^2- |S_{\otimes,t} j(g)|^2\Big)^\frac{1}{2} \Big\|_{\widehat{\RR}_2}. 
 $$
Therefore,  it suffices to show that $\Phi$ is $\widehat{\RR}_1 \to \mathrm{BMO}_{\T_\otimes}^c(\widehat{\RR}_2)$ bounded. We find
$$\Phi(\xi)(x) = \int_\Omega K(x,y) (\xi(y)) \hskip1pt d\mu(y).$$
Thus, we may regard $\Phi$ as a semicommutative CZO and apply Theorem \ref{CZ2}  where $\M_j$ is replaced by $\widehat{\M}_j$.  Since $\Phi(\xi) = \sum_{g,h}  (\alpha_{g^{-1}}) \bullet  (T_{gh^{-1}}) \bullet (\alpha_g) \bullet \xi$ and $\beta$ is measure preserving, we immediately find that the $L_2^c$-boundedness assumption implies that the map $$\Phi: L_2^c(\Omega) \bar\otimes \widehat{\M}_1 \to L_2^c(\Omega) \bar\otimes \widehat{\M}_2$$ is bounded. Moreover, the smoothness condition matches that of Theorem \ref{CZ2}. \fin
 
\begin{remark} \label{GralNonEquiv}
\emph{ Our work so far yields sufficient conditions for the $L_\infty \to \mathrm{BMO}$ boundedness of $T \rtimes id_\G$ in more general settings.  In particular, if $T_g=T$ and $\alpha_g T=T\alpha_g$ for all $g\in \G$, then we find for any  $T$ fulfilling the assumption of Theorem \ref{CZ2},  $T\rtimes id_\G: \RR_1\rtimes \G \mapsto \mathrm{BMO}_{\T_\rtimes}^c(\RR_2 \rtimes \G)$ is  bounded.}
\end{remark}

\subsection{Matrix algebras}

In this paragraph, we introduce a Markov metric for the matrix algebra $ \mathcal{B}(\ell_2)$. The triangular truncation plays the noncommutative form of the Hilbert transform on $\mathcal{B}(\ell_2)$. We shall reprove the $L_p$-boundedness of the triangular truncation for $1 < p < \infty$ and a new $\mathrm{BMO} \to \mathrm{BMO}$ estimate by means of this Markov metric and our algebraic approach. Consider the $*$-homomorphism $u: \mathcal{B}(\ell_2) \to L_\infty(\mathbb{R}) \bar\otimes \mathcal{B}(\ell_2)$ determined by $$u (e_{mk}) = e^{2\pi i (m-k) \, \cdot} e_{mk}.$$ 
Given $A = \sum_{m,k} a_{mk} e_{mk}$, define the semigroup 
$$S_t(A) = \summ_{m,k} e^{- t |m-k|^2} a_{mk} e_{mk}.$$ 
It is not difficult to see that it defines a Markov semigroup of convolution type. In fact, $u$ is a corepresentation of $L_\infty(\mathbb{R})$ (equipped with its natural comultiplication map $\Delta f(x,y) = f(x+y)$) in $\mathcal{B}(\ell_2)$ and it turns out that $\T = (S_t)_{t \ge 0}$ is the transferred semigroup associated to the heat semigroup on $\mathbb{R}$ 
$$
u\circ S_t = (H_t \otimes id_{\mathcal{B}(\ell_2)}) \circ u.
$$ 
Define the cpu map $R_{j,t}$ on $\mathcal{B}(\ell_2)$ by $u\circ R_{j,t}=(\widetilde{R}_{j,t}\otimes id_{\mathcal{B}(\ell_2)})\circ u$, where $\widetilde{R}_{j,t}f(x)$ denotes the average of $f \in L_\infty(\R)$ over the interval $\mathrm{B}_{\sqrt{4jt}}(x)$. Now, given a matrix $A = \sum_{m,k} a_{mk} e_{mk}$ we find 
\begin{eqnarray*}
u\circ R_{j,t}(A)(x) &= & \mean_{\mathrm{B}_{\sqrt{4jt}}(x)}u(A)(y)dy\\
& =&  \mean_{\mathrm{B}_{\sqrt{4jt}}(x)} \sum_{m,k}e^{2\pi i(m-k)y}a_{mk}e_{mk}dy\\
& = & \sum_{m,k} \frac{\sin (4\sqrt{jt}\pi (m-k))}{4\sqrt{jt}\pi (m-k)}e^{2\pi i(m-k)x}a_{mk}e_{mk}.
\end{eqnarray*}
Thus, we find the following identity $$R_{j,t}(A) = \sum_{m,k} \frac{\sin (4\sqrt{jt}\pi (m-k))}{4\sqrt{jt}\pi (m-k)} a_{mk}e_{mk}.$$ Taking $\sigma_{j,t}^2 = 2e \sqrt{j/\pi} e^{-j} \1_{\mathcal{B}(\ell_2)}$ and $\gamma_{j,t}^2 = \sqrt{j} \1_{\mathcal{B}(\ell_2)}$, we obtain a Markov metric in $\mathcal{B}(\ell_2)$. Indeed, the metric integrability condition holds trivially, as for the Hilbert module majorization it reduces to prove that $B_1 \le B_2$ with 
\begin{eqnarray*}
B_1 & = & u \big( \langle \xi, \xi \rangle_{S_t} \big) \, = \, \big\langle u \otimes u (\xi), u \otimes u (\xi) \big\rangle_{H_t \otimes id_{\mathcal{B}(\ell_2)}}, \\
B_2 & = & \sum_j \sigma_{j,t}^2 u \big( \langle \xi, \xi \rangle_{R_{j,t}} \big) \, = \, \sum_j \sigma_{j,t}^2 \big\langle u \otimes u (\xi), u \otimes u (\xi) \big\rangle_{\widetilde{R}_{j,t} \otimes id_{\mathcal{B}(\ell_2)}}. 
\end{eqnarray*}
In other words, it suffices to note that the canonical Markov metric in $\R$ |which recovers the Euclidean metric, as proved in Paragraph \ref{SectEuclideanMetric}| admits a matrix-valued extension, as it was justified in Remark \ref{ossMarkovBMO}. Let us now consider the triangular truncation 
$$
\triangle (A)=\sum_{m>k} a_{mk}e_{mk}.
$$

\begin{corollary} \label{CorTT}
We have $$\|\triangle (A)\|_{\mathrm{BMO}_\T} \lesssim \| A \|_{\mathrm{BMO}_\T}.$$
In particular, given $1 < p < \infty$ we obtain $\displaystyle \| \triangle (A)\|_{S_p} \lesssim \frac{p^2}{p-1}\| A \|_{S_p}$.
\end{corollary}

\dem Recall that 
$$
u\circ \triangle = (L \otimes id_{ \mathcal{B}(\ell_2)})\circ u.
$$
for $L = \frac{1}{2}(id+i H)$ and $\widehat{Hf}(\xi) = -i \mathrm{sgn}(\xi) \widehat{f}(\xi)$, the Hilbert transform in the real line. We may also regard $u: \mathcal{B}(\ell_2) \to L_\infty(\mathbb{T}) \bar\otimes \mathcal{B}(\ell_2)$ as a corepresentation of $\mathbb{T}$ instead of $\R$ and the above identity holds replacing $H$ by the Hilbert transform in the torus. In this case, $u$ becomes a trace preserving $*$-homomorphism and the well-known $S_p$ inequalities for $\triangle$ reduce to the boundedness of the Hilbert transform in $L_p(\mathbb{T}; S_p(\ell_2))$, which is also well-known and follows in passing from the semicommutative theory in the previous paragraph. Alternatively, the second assertion follows from the first one by interpolation and duality. According to Remark \ref{Remark-EquivBMO}, to prove the first assertion it suffices to show that the map $$T=  i(id_{B(\ell_2)}-2\triangle)$$ is $\mathrm{BMO} \rightarrow \mathrm{BMO}$ bounded for the semigroup BMO space which is associated to the transferred Poisson semigroup $P_t$ on $B(\ell_2)$ given by $${P_t}: (a_{ij})  \mapsto ({e^{-t|i-j|}}a_{ij}).$$ Given $A=(a_{jk})_{j,k}$ in $B(\ell_2)$ then
 \begin{eqnarray*}
|A|^2 & = & A^*A \ = \ \Big( \sum_k \overline{a_{ki}}a_{kj} \Big)_{i,j}\\
T(A) & = & i \Big( \mathrm{sgn}(k-j)a_{jk} \Big)_{j,k}\ \ (TA)^* = i \Big( \mathrm{sgn}(k-j)\overline{a_{kj}} \Big)_{j,k}
\end{eqnarray*} 
Then $\big( P_t|A|^2-|P_tA|^2 \big)_{ij} = \sum_k (e^{-t|i-j|} - e^{-t|k-j|}e^{-t|i-k|} )\overline{a_{ki}} a_{kj}$ and
\begin{eqnarray*}
\lefteqn{\hskip-20pt \big( P_t|T(A)|^2 - |P_tT(A)|^2 \big)_{ij}} \\ \!\!\! & = & \!\!\! \sum_k (e^{-t|i-j|} - e^{-t|k-j|} \, e^{-t|i-k|}) \, \mathrm{sgn}(k-i) \, \mathrm{sgn}(k-j) \, \overline{a_{ki}} a_{kj}.
\end{eqnarray*}
Since $\mathrm{sgn}(k-i) \mathrm{sgn}(k-j) \neq 1$ iff $e^{-t|i-k|}e^{-t|k-j|}=e^{-t|i-j|}$, we get $$P_t|A|^2 - |P_tA|^2 \, = \, P_t|T(A)|^2 - |P_tT(A)|^2.$$ The last identity implies that $T$ is an isometry on the Poisson BMO space. \fin

\subsection{Quantum Euclidean spaces}

Given an integer $n\geq 1$, fix an anti-symmetric $\mathbb{R}$-valued $n\times n$ matrix $\Theta$. We define $\rm{A}_{\Theta}$ as the universal C*-algebra generated by a family $u_1(s), u_2(s),\cdots, u_n(s)$ of one-parameter unitary groups in $s\in \mathbb{R}^n$ which are strongly continuous and satisfy the following $\Theta$-commutation relations 
$$
u_j(s)u_k(t)=e^{2\pi i \Theta_{jk}st}u_k(t)u_j(s).
$$
If $\Theta =0$, by Stone's theorem we can take $u_j(s)=\exp (2\pi is\langle e_j,\cdot \rangle)$ and $\rm{A}_{\Theta}$ is the space of bounded continuous functions on $\R^n$. In general, given $\xi \in \R^n$, define the unitaries $\lambda_\Theta (\xi)=u_1(\xi_1)u_2(\xi_2)\cdots u_n(\xi_n)$. Let $\rm{E}_{\Theta}$ be the closure in $\rm{A}_{\Theta}$ of $\lambda_\Theta (L_1(\R^n))$ with 
$$
f=\int_{\R^n}  \check{f}_\Theta(\xi) \lambda_\Theta (\xi)\,d\xi.
$$
If $\Theta=0$, $\rm{E}_{\Theta}=\mathcal{C}_0(\R^n)$. Define
$$
\tau_\Theta(f) = \tau_\Theta \left( \int_{\R^n} \check f_\Theta(\xi) \lambda_\Theta (\xi)\, d\xi \right) = \check f_\Theta (0)
$$
for $\check f_\Theta: \R^n \rightarrow \mathbb{C}$ integrable and smooth. $\tau_\Theta$ extends to a normal faithful semifinite trace on $\rm{E}_{\Theta}$. Let $\mathcal{R}_\Theta=\rm{A}_{\Theta}''=\rm{E}_{\Theta}''$  be the von Neumann algebra generated by $\rm{E}_{\Theta}$ in the GNS representation of $\tau_\Theta$. Note that if $\Theta=0$,  $\mathcal{R}_\Theta=L_\infty(\R^n)$. In general we call $\mathcal{R}_\Theta$ a quantum Euclidean space. There are two maps which play important roles while doing analysis over quantum Euclidean spaces. The first one is the corepresentation map $\sigma_\Theta: \mathcal{R}_\Theta \rightarrow L_\infty(\R^n)\bar{\otimes}\mathcal{R}_\Theta
$, given by $\lambda_\Theta(\xi) \mapsto \exp_\xi \otimes \lambda _\Theta(\xi)$ where $\exp_\xi$ stands for the Fourier character $\exp (2\pi i \langle \xi, \cdot \rangle)$. Note that $\sigma_\Theta$ is a normal injective $*$-homomorphism. The second map is $\pi_\Theta: \exp_\xi \mapsto \lambda_\Theta(\xi)\otimes \lambda_\Theta (\xi)^*$, which extends to a normal $*$-homomorphism from $L_\infty(\R^n)$ to $ \mathcal{R}_\Theta \bar{\otimes} \mathcal{R}_\Theta^{\rm{op}}$,  where $\mathcal{R}_\Theta^{\rm{op}}$ is the apposite algebra of $\mathcal{R}_\Theta$, which is obtained by preserving the linear and adjoint structures but reversing the product. We refer the readers to \cite{GPJP} for more detailed information of quantum Euclidean spaces and these two maps.

\noindent \textbf{BMO and Markov metric.} Our first goal is to construct a natural Markov metric for quantum Euclidean spaces. Let us recall the heat semigroup on $\R^n$ acting on $\varphi :\R^n \rightarrow \mathbb{C}$ admits the following form 
$$
H_t \varphi(x)=\int_{\R^n}\widehat{\varphi}(\xi)e^{-t |\xi |^2}\exp _\xi (x)\, d\xi.
$$
This induces a semigroup on $\mathcal{R}_\Theta$ determined by $$\sigma_\Theta \circ S_{\Theta,t} = (H_t \otimes id_{\mathcal{R}_\Theta})\circ \sigma_\Theta.$$ $S_{\Theta,t}$ gives a Markov semigroup on $\RR_\Theta$ which formally acts as 
\begin{equation} \label{Eq-QESSemigroup}
S_{\Theta,t}(f) = \int_{\R^n}  \check{f}_\Theta(\xi) e^{-t |\xi |^2} \lambda_\Theta (\xi)\,  d\xi.
\end{equation}
The corresponding semigroup column $\rm{BMO}$ norm is given by
\begin{eqnarray*}
\| f \|_{\rm{BMO}_c(\RR_\Theta)} \!\!\! & = & \!\!\! \sup_{t>0} \Big\| \Big( S_{\Theta,t}(|f |^2)-|S_{\Theta,t}(f) |^2\Big)^
\frac{1}{2} \Big\|_{\RR_{\Theta}} \\
\!\!\! & \approx & \!\!\! \hskip-5pt \sup_{\begin{subarray}{c} \rm{B} \, \mathrm{ball \, in} \, \R^n \end{subarray}} \Big\|  \Big( \mean _{\rm{B}} |\sigma_\Theta(f)-\sigma_\Theta (f)_{\rm{B}} |^2 d\mu \Big)^
\frac{1}{2}\Big\|_{\RR_{\Theta}} \hskip-5pt = \| \sigma _\Theta (f) \|_{\rm{BMO}_c(\R^n; \RR_\Theta)}.
\end{eqnarray*}
According to Remark \ref{ossMarkovBMO}, the semicommutative extension $H_t\otimes id_{\mathcal{R}_\Theta}$ of the heat semigroup, together with the extension of the corresponding Markov metric from Paragraph \ref{ClassicalCZT} still satisfies the Hilbert module majorization 
\renewcommand{\theequation}{\arabic{section}.\arabic{equation}}
\begin{equation}\label{eq: opvalued upper estimate}
\langle \xi, \xi \rangle _{H_t \otimes id_{\mathcal{R}_\Theta}} \leq \sum_{j \geq 1} \sigma _{j,t}^* \langle \xi, \xi \rangle _{R_{j,t} \otimes id_{\mathcal{R}_\Theta}}   \sigma _{j,t} 
\end{equation}
as well as the integrability condition, where $\sigma_{j,t}^2 \equiv 2e \sqrt{j^n/\pi} e^{-j}$, $\gamma_{j,t}^2 \equiv j^{\frac{n}{2}}$ and $R_{j,t}f(x)$ is the average of $f$ over $\mathrm{B}_{\sqrt{4jt}}(x)$. Then we can easily produce a Markov metric on $\RR_\Theta$. Let $\mathrm{B}_{j,t}$ be the Euclidean ball in $\R^n$ centered at the origin with radius $\sqrt{4jt}$ and consider the projections $q_{j,t} = \chi_{\mathrm{B}_{j,t}} \ten \1_{\RR_\Theta}$. Define the cpu maps  
$$
R_{\Theta,j,t}(f) = \frac{1}{|\mathrm{B}_{j,t}|} \int_{\mathrm{B}_{j,t}}  \sigma_\Theta(f)(x) \, dx = \frac{1}{|\mathrm{B}_{j,t}|} \int_{\R^n} \widehat{\chi}_{\mathrm{B}_{j,t}}(\xi) \check{f}_\Theta(\xi) \lambda_\Theta (\xi)\, d\xi.
$$
It is easy to check that 
\begin{equation} \label{Eq-QESQM}
\sigma_\Theta \circ R_{\Theta,j,t}= (R_{j,t}\otimes id_{\mathcal{R}_\Theta})\circ \sigma_\Theta.
\end{equation}
The Hilbert module majorization
$$
\langle \xi, \xi \rangle_{S_{\Theta,t} } \le \sum_{j\geq 1} \sigma _{j,t}^* \langle \xi, \xi \rangle _{R_{\Theta,j,t} } \sigma _{j,t} 
$$
for $\xi \in \mathcal{R}_\Theta \bar\otimes_{S_{\Theta,t}} \mathcal{R}_\Theta$ is equivalent to the same inequality after  composing with the $*$-homomorphism $\sigma_\Theta$, which follows in turn by the intertwining identities \eqref{Eq-QESSemigroup} and \eqref{Eq-QESQM}, together with the majorization \eqref{eq: opvalued upper estimate}. Therefore, we obtain a Markov metric on $\mathcal{R}_\Theta$ associated to $S_\Theta$ $$\Q_\Theta = \big\{ (R_{\Theta, j,t}, \sigma_{j,t}, \gamma_{j,t})\mid (j,t)\in \mathbb{Z}_+ \times \mathbb{R}_+ \big\}.$$

\noindent \textbf{The algebraic structure.} We start with the kernel representation of our CZOs over the (fully noncommutative) von Neumann algebra $\mathcal{R}_\Theta$. Given a kernel $k$ affiliated to $\mathcal{R}_\Theta \bar{\otimes}\mathcal{R}_\Theta^{\rm op}$, the linear map associated to it is formally given by $$
T_k f = (id_{\RR_\Theta ^{\rm}}\otimes \tau_\Theta) \big( k(\1_{\RR_\Theta}  \otimes f) \big) = (id_{\RR_\Theta ^{\rm}}\otimes \tau_\Theta) \big( (\1_{\RR_\Theta} \otimes f)k \big).
$$
The reader is referred to \cite{GPJP} for more details. Our goal is to provide sufficient conditions for the $L_\infty\to\rm BMO$ boundedness of $T_k$. Consider the $*$-homomorphism $\sigma_\Theta:$ $\RR_\Theta \rightarrow L_\infty(\R^n)\bar\otimes \RR_\Theta$. In the case of quantum Euclidean spaces, we need the full algebraic skeleton introduced in Section \ref{Sect-CZ}. In Table 1 there is a little dictionary to identify the main objects.
\begin{table}[htbp]
\begin{center}
\begin{tabular}{|l|l|}
\hline
Generic algebraic objects & Quantum Euclidean spaces \\
\hline \hline
$\M$ & $\RR_\Theta$ \\ \hline
$\mathcal{N}_\rho$ & \hskip-5pt $\begin{array}{l} L_\infty(\R^n) \bar\otimes \RR_\Theta \\ \rho_1 = \1 \ten \cdot, \rho_2 = \sigma_\Theta \\ \mathsf{E}_\rho = \mathrm{Lebesgue \ integral} \end{array}$ \\ \hline
$\mathcal{N}_\pi$ & \hskip-5pt $\begin{array}{l} \RR_\Theta \bar\otimes \RR_\Theta^{\mathrm{op}} \\ \pi_1 = \1 \ten \cdot, \pi_2 = \cdot \ten \1 \\ \mathsf{E}_\pi = \tau_\Theta \ten id_{\RR_\Theta^{\mathrm{op}}} \\ [-10pt] \null \end{array}$ \\ \hline
$\mathcal{N}_\sigma$ & \hskip-5pt $\begin{array}{l} L_\infty(\R^n) \bar\otimes L_\infty(\R^n) \bar\otimes \RR_\Theta \\ \mathcal{N}_\sigma = (\delta \ten id_{\RR_\Theta}) (\mathcal{N}_\rho) \\ \delta \varphi(x,y) = \varphi(x) - \varphi(y) \end{array}$ \\ \hline
\end{tabular}
\vskip5pt
\caption{Algebraic skeleton for $\RR_\Theta$.}
\label{tabla:sencilla}
\end{center}
\end{table}
Next, note that
$$
 \sigma_\Theta \circ T_k (f)= (id_{\R^n} \ten id_{\RR_\Theta} \ten \tau _\Theta) \big( k_\sigma (\1_{\R^n} \otimes \1_{\RR_\Theta}\ten f) \big),
$$
where $k_\sigma =(\sigma_\Theta \otimes id_{\RR_\Theta^{\rm op}})(k)$. Denote $\sigma_\Theta \circ T_k$ by $T_{k_\sigma}$. Define  
$$
\widehat{T}_k : \RR_\Theta \bar\ten \RR_\Theta \ni  f \ten a \mapsto T_{k_\sigma} (f)(\1_{\R^n} \ten a)  \in L_\infty (\R^n) \bar \ten \RR_\Theta.
$$
Then it is clear that the compatibility condition \eqref{Amplification} holds since $\widehat{T}_k \circ \pi_2 = \sigma_\Theta \circ T_k$.

\begin{lemma}\label{lem: L2 norm control}
If $T_k$ is bounded on $L_2(\RR_\Theta)$, then 
$$
\big\|  \widehat{T}_k : L_2^c(\RR_\Theta) \bar \ten \RR_\Theta \rightarrow L_2^c(\R^n)\bar \ten \RR_\Theta \big\| \leq \big\|   T_k : L_2(\RR_\Theta)\rightarrow L_2(\RR_\Theta) \big\|.
$$
\end{lemma}

\dem We need to introduce two maps:
$$
j_\Theta: L_2(\R^n) \ni \int_{\R^n} \varphi(\xi) \exp_{\xi} \, d\xi \mapsto \int_{\R^n} \varphi(\xi) \lambda_\Theta (\xi) \, d\xi \in L_2(\RR_\Theta),
$$
$$
W: L_2^c(\R^n) \bar \ten \RR_\Theta \ni \int_{\R^n}\exp_{\xi} \ten a(\xi)\,d\xi \mapsto  \int_{\R^n}\exp_{\xi} \ten \lambda_\Theta(\xi) a(\xi)\,d\xi \in  L_2^c(\R^n) \bar \ten \RR_\Theta. 
$$
It is straightforward to show that $W$ extends to an isometry. Moreover, $j_\Theta$ is also an $L_2$-isometry, we refer the reader to \cite[Section 1.3.2]{GPJP} for the proof. Observe that 
\begin{eqnarray*}
\sigma_\Theta (f)(\1_{\R^n} \ten a ) & = & \int_{\R^n}  \check{f}_\Theta(\xi) \exp_{\xi} \ten  \lambda_\Theta (\xi)a \,d\xi\\
& = & W ( \int_{\R^n}  \check{f}_\Theta(\xi) \exp_{\xi} \ten \,a  \,d\xi) \ = \ W \circ (j_\Theta^*  \ten id_{\RR_\Theta} )(f \ten a ).
\end{eqnarray*}
Letting $f = T_k g$, we get $$\widehat{T}_{k} (g \ten a) = W(j_\Theta^* T_k \ten id_{\RR_\Theta})(g \ten a ).$$ The properties of the maps $j_\Theta$ and $W$ readily imply the assertion. \fin

Now let us introduce a weak-$*$ dense subalgebra of $\RR_\Theta$, which is the analogue of the classical Schwartz class. Let $\T(\R^n)$ denote the classical Schwartz class in the Euclidean space $\R^n$ and define 
$$
\T_\Theta =\big\{ f\in \RR_\Theta: \check f_\Theta \in \T(\R^n) \big\}.
$$
Let $\T_\Theta '$ denote the space of continuous linear functionals on $\T_\Theta$, which is the quantum space of tempered distributions. Consider a continuous linear operator $T\in \mathcal{L}(\mathcal{S}_\Theta,  \mathcal{S}_\Theta ')$. By using the unitary map $$j_\Theta: \T(\R^n) \to \T_\Theta$$ as defined in the proof of Lemma \ref{lem: L2 norm control}, we get $j_\Theta^* T j_\Theta \in \mathcal{L}(\mathcal{S}(\R^n),  \mathcal{S}(\R^n) ')$. By a result of Schwartz, there exists a unique kernel $K \in \T '(\R^{2n})= (\T(\R^n)\ten_\pi \T(\R^n))'$ such that $T$ admits the kernel $k=(j_\Theta \ten j_\Theta)(K)\in (\mathcal{S}_\Theta \ten_\pi \mathcal{S}_\Theta)'$. Actually, the kernel representations $T_k$ satisfying the Calder\'on-Zygmund type conditions in the following theorem belong to $\mathcal{L}(\mathcal{S}_\Theta,  \mathcal{S}_\Theta ')$. It provides sufficient conditions for the $L_\infty(\RR_\Theta)\rightarrow\rm {BMO}_c(\RR_\Theta)$ boundedness of CZO operators associated to kernels in $ (\mathcal{S}_\Theta \ten_\pi \mathcal{S}_\Theta)'$.
We shall use the quantum analogue of the bands around the diagonal $$a_\mathrm{B} = \pi_\Theta(\chi_{5\mathrm{B}}) = \int_{\R} \widehat{\chi}_{5\mathrm{B}}(\xi) \lambda_\Theta (\xi) \otimes \lambda_\Theta (\xi)^* \, d\xi.$$ 

\begin{theorem}\label{thm: CZ QE}
Let $T_k\in \mathcal{L}(\mathcal{S}_\Theta,  \mathcal{S}_\Theta')$ and assume
\begin{itemize}
\item[\emph{i)}] Cancellation $$\big\| T_k: L_2(\RR_\Theta) \to L_2(\mathcal{R}_\Theta) \big\| < \infty.$$

\vskip3pt

\item[\emph{ii)}] For any $f \in \RR_\Theta$ and any Euclidean ball $\mathrm{B}$ centered at the origin 
$$
\mean_{\mathrm{B} \times \mathrm{B}} \big| \Sigma_{\Theta, k, f, \mathrm{B}}(y_1) - \Sigma_{\Theta, k, f, \mathrm{B}}(y_2) \big|^2 dy_1dy_2 
\lesssim \| f\|_{\RR_\Theta}^2,
$$
where $\Sigma_{\Theta, k, f, \mathrm{B}} = (id_{\R^n} \otimes id_{\RR_\Theta} \ten \tau_\Theta) \big[ k_\sigma(\1_{\R^n} \otimes \1_{\RR_\Theta} \ten f)(\1_{\R^n} \otimes a_{\mathrm{B}}^\perp) \big]$.
\end{itemize}
Then, the Calder\'on-Zygmund operator $T_k$ is bounded from $L_\infty(\RR_\Theta)$ to ${\rm {BMO}}_c(\RR_\Theta)$.
\end{theorem}

\dem By Theorem \ref{MetricBMOInterp}, it suffices to prove $$T_k: L_\infty(\RR_\Theta) \to \mathrm{BMO}_{\Q_\Theta}^c.$$ Arguing as in Paragraph \ref{SectEuclideanMetric}, the Markov metric BMO norm takes the simpler form 
$$
\|f\|_{\mathrm{BMO}_{\Q_\Theta}^c} = \sup_{t > 0} \sup_{j \geq 1} \Big\| \Big( \gamma_{j,t}^{-1} \big[ R_{\Theta,j,t}( | f |^2 ) - | R_{\Theta,j,t}   (f) |^2 \big] \gamma_{j,t}^{-1} \Big)^\frac{1}{2} \Big\|_{\RR_\Theta}.
$$
In other words, the extra term in the definition of BMO is dominated by the above expression as in \eqref{Eq-BMODoubling}. As noticed in Remark \ref{rem: no size condition}, the size kernel condition is then superfluous. This also reduces the analytic conditions and the smooth kernel conditions to be checked. In summary, according to the proof of Theorem \ref{Extrapolation}, the assertion will follow if we can justify: 

\begin{itemize}
\item[C0)] Initial condition $$T_k: \A_\Theta \to \RR_\Theta \quad \mbox{for $\A_\Theta \subset \RR_\Theta$ weak-$*$ dense.}$$

\item[Al1)] $\Q_\Theta$-monotonicity of $\mathsf{E}_\rho$ $$\hskip20pt \mathsf{E}_\rho(q_{j,t} |\xi|^2 q_{j,t}) \le \mathsf{E}_\rho(|\xi|^2).$$

\vskip3pt

\item[Al2)] Right modularity of $\widehat{T}_k$ $$\hskip20pt \widehat{T}_k(\eta \pi_1(b)) \hskip1pt = \hskip1pt \widehat{T}_k(\eta) \rho_1(b).$$

\vskip3pt

\item[An1)] Mean differences $$\hskip20pt  \widehat{R}_{\Theta,j,t}(\xi^*\xi) - \widehat{R}_{\Theta,j,t}(\xi)^*\widehat{R}_{\Theta,j,t}(\xi) \, \le \, \Phi_{j,t} \big( \delta(\xi)^*\delta(\xi)) \big) \ \mbox{for some cpu $\Phi_{j,t}$}.$$

\vskip3pt

\item[An2)] Metric/measure growth $$\hskip20pt  \mathbf{1} \, \le \, \pi_1 \rho_1^{-1} \E_\rho(q_{j,t})^{-\frac12} \E_\pi (a_{j,t}^* a_{j,t}) \pi_1 \rho_1^{-1} \E_\rho(q_{j,t})^{-\frac12} \, \lesssim \, \pi_1 \rho_1^{-1}(\gamma_{j,t}^2).$$

\vskip3pt

\item[CZ1)] $L_2^c$-boundedness condition $$\hskip20pt \widehat{T}_k: L_\infty^c(\mathcal{N}_\pi; \mathsf{E}_\pi) \to L_\infty^c(\mathcal{N}_\rho; \mathsf{E}_\rho).$$

\item[CZ2)] Kernel smoothness condition $$\hskip20pt \Phi_{j,t} \Big( \big| \delta \big( \widehat{T}_k(\pi_2(f)(\mathbf{1} \hskip1.5pt - \hskip1.5pt a_{j,t})) \big) \big|^2 \Big) \, \lesssim \, \gamma_{j,t}^2 \|f\|_\infty^2.$$
\end{itemize}

The initial condition trivially holds for good kernels $k \in \S_\Theta \ten_\mathrm{alg} \S_\Theta$. In \cite{GPJP} it was required to extend the main result from this class of kernels to general ones in $\S_{\Theta \oplus \Theta}'$, by reproving certain auxiliary results in the context of distributions. In our case, this is much simpler. Indeed, when dealing with general kernels, we just note that $T_k(f) \in L_2(\RR_\Theta)$ for all $f \in \S_\Theta$ by assumption. Given the form of $R_{\Theta,j,t}$, it trivially follows that $R_{\Theta,j,t} (|T_k f|^2)$ and $R_{\Theta,j,t} T_k f$ are well-defined operators in $L_1(\RR_\Theta)$ and $L_2(\RR_\Theta)$ respectively. In particular, the proof of Theorem \ref{Extrapolation} follows exactly as it was written there under this more flexible assumption. Therefore, the initial condition can be relaxed to the condition $$T_k: \S_\Theta \to L_2(\RR_\Theta).$$ In fact, according to \cite[Proposition 2.17]{GPJP}, every algebraic column CZO is normal. Thus, it suffices |as we did in Theorem \ref{Extrapolation}| to justify that $T_k: \S_\Theta \to \mathrm{BMO}_{\Q_\Theta}^c$ is bounded, as we shall do by justifying the remaining conditions.   

Al1 holds trivially since $q_{j,t} = \chi_{\mathrm{B}_{j,t}} \ten \1_{\RR_\Theta}$ lives in the center of $\mathcal{N}_\rho$. On the other hand, according to the definition of $\rho_1, \pi_1$ from Table 1, the algebraic condition Al2 can be rewritten as follows $$\widehat{T}_k \big( \eta (\1_{\RR_\Theta} \ten b) \big) = \widehat{T}_k (\eta) (\1_{\R^n} \ten b).$$ This is clear from the definition of $\widehat{T}_k$. Next, condition An1 reads as $$\mean_{\mathrm{B}_{j,t}} |\xi|^2 d\mu - \Big| \mean_{\mathrm{B}_{j,t}} \xi d\mu \Big|^2 \le \mean_{\mathrm{B}_{j,t} \times \mathrm{B}_{j,t}} \big| \xi(x) - \xi(y)\big|^2 d\mu(x) d\mu(y)$$ for $\RR_\Theta$-valued functions, when $\Phi_{j,t}$ is chosen to be the average over $\mathrm{B}_{j,t} \times \mathrm{B}_{j,t}$. As in \eqref{Eq-Jensen}, this is a consequence of the operator-valued Jensen's inequality. Next recalling that $a_{j,t} = \pi_\Theta(\chi_{5 \mathrm{B}_{j,t}})$, condition An2 takes the form
$$
|\mathrm{B}_{j,t}| \1_{\RR_\Theta} \, \le \, (\tau_\Theta \ten id_{\RR_\Theta^{\mathrm{op}}}) (\pi_\Theta(\chi_{5 \mathrm{B}_{j,t}})) \, \lesssim \, j^{\frac{n}{2}} |\mathrm{B}_{j,t}| \1_{\RR_\Theta}.
$$
To verify it we note that 
$$
(\tau_\Theta \ten id_{\RR_\Theta^{\mathrm{op}}}) (\pi _\Theta (\varphi)) = (\tau_\Theta \ten id_{\RR_\Theta^{\mathrm{op}}}) \Big(\int_{\R^n}\widehat{\varphi}(\xi)\lambda_\Theta(\xi)\ten\lambda_\Theta(\xi)^* d\xi \Big)=\widehat{\varphi}(0)\1_{\RR_\Theta }.
$$
Then we get $(\tau_\Theta \ten id_{\RR_\Theta^{\mathrm{op}}}) (\pi_\Theta(\chi_{5 \mathrm{B}_{j,t}})) = 5|\mathrm{B}_{j,t}|\1_{\RR_\Theta }$. Condition CZ1 reduces to our $L_2$-boundedness assumption by Lemma \ref{lem: L2 norm control}. Finally, the smoothness condition ii) in the statement readily implies condition CZ2 for all values of $j,t$. \fin 

The smoothness condition in Theorem \ref{thm: CZ QE} is of H\"ormander type, while the one in the main result of \cite{GPJP} is a gradient condition. As expected, we shall show that our condition in this paper is more flexible than that of \cite[Theorem 2.6]{GPJP}. We use $\bullet$ for the product in $\M \bar\ten \M^{\rm op}$, so that 
$$
(a \ten b)\bullet (a'\ten b')= (a a')\ten (b'b).
$$
 The quantum analogue of the metric  is defined by 
$$
{\rm d}_\Theta =\pi _{\Theta}(| \cdot |)
$$
for the Euclidean norm $| \cdot |$. Moreover, we also introduce the $\Theta$-deformation of the  free gradient. Let $\mathcal{L}(\mathbb{F}_n)$ denote the group von Neumann algebra associated to the free group  with $n$ generators $\mathbb{F}_n$. It is well-known from (say) \cite{VDN} that $\mathcal{L}(\mathbb{F}_n)$ is generated by $n$ semicircular random variables $s_1, s_2, \ldots, s_n$. Note that there exist  derivations $\partial_\Theta^j$ in $\T_\Theta$ which are determined by 
$$
\partial_\Theta^j (\lambda_\Theta(\xi))=2\pi i \xi_j \lambda_\Theta(\xi)
$$
for $1\leq j\leq n$. Define the $\Theta$-deformed free gradient as
$$
\nabla_\Theta =\sum_{j=1}^n s_j \ten \partial_\Theta ^j : \T_\Theta \rightarrow \mathcal{L}(\mathbb{F}_n) \bar \ten \RR_\Theta.
$$
If $\nabla$ denotes the free gradient for $\Theta =0$, it is easy to check that 
\begin{eqnarray} \label{EqFreenonFreeGrad}
( id_{\mathcal{L}(\mathbb{F}_n)}\ten \sigma_{\Theta})\circ \nabla_\Theta & = & \sum_{j=1}^n s_j \ten (\sigma_\Theta \circ  \partial_\Theta ^j) \\ \nonumber & = & \sum_{j=1}^n s_j \ten (\partial_j \circ \sigma_\Theta ) \ = \ ( \nabla\ten id_{\RR_\Theta})\circ \sigma_\Theta.
\end{eqnarray}
For the convenience of the reader, we cite Theorem 2.6 from \cite{GPJP} below. 
\begin{theorem}\label{thm: cite}
Let  $T_k\in \mathcal{L}(\mathcal{S}_\Theta,  \mathcal{S}_\Theta ')$ and assume:
\begin{itemize}
\item[\emph{i)}] Cancellation 
$$
\| T_k : L_2(\RR_\Theta) \rightarrow L_2(\RR_\Theta) \| \leq A_1.
$$
\item[\emph{ii)}] Gradient condition. There exists 
$$
\alpha < \frac{n}{2} <\beta< \frac{n}{2}+1
$$
satisfying the gradient conditions below for $\rho =\alpha, \beta$
$$
\Big| {\rm d}_\Theta ^\rho \bullet (\nabla _\Theta \ten id_{\RR_\Theta ^{\rm op}}) (k) \bullet {\rm d}_\Theta ^{n+1-\rho} \Big| \leq A_2. 
$$
\end{itemize}
Then, we find the following $L_\infty \rightarrow {\rm{BMO}}_c$ estimate 
$$
\big\| T_k: L_\infty(\RR_\Theta)\rightarrow {\rm {BMO}}_c(\RR_\Theta)\big\| \leq C_n(\alpha,\beta)(A_1+A_2).
$$
\end{theorem}

To simplify notation, we shall write in what follows $\Sigma$ for $\Sigma_{\Theta,k,f,\mathrm{B}}$. According to the semicommutative Poincar\'e type inequality introduced in \cite[Proposition 1.6]{GPJP} we obtain 
$$
\Big\| \mean_{\mathrm{B} \times \mathrm{B}} \big|\delta (\Sigma) \big|^2 d\mu \times \mu \Big\| _{\RR_\Theta} \leq 16 \mathrm{R}^2 \Big\| (\1\ten \chi_{\mathrm{B}} \ten \1) (\nabla\ten id_{\RR_\Theta}){(\Sigma)} \Big\|^2_{\mathcal{L}(\mathbb{F}_n) \bar\ten L_\infty(\R^n)\bar\ten \RR_\Theta}
$$
for $\mathrm{R} =$ radius of $\mathrm{B}$. By \eqref{EqFreenonFreeGrad}, we may rewrite 
\begin{eqnarray*}
\lefteqn{(\1\ten \chi_\mathrm{B} \ten \1) (\nabla\ten id_{\RR_\Theta}) (\Sigma)} \\
\!\!\!\! & = & \!\!\!\! (id^{\ten 3}\ten \tau_\Theta) \Big((\1\ten \chi_{\mathrm{B}}\ten \1^{\ten 2})(\nabla \ten id^{\ten 2})(k_\sigma) (\1^{\ten 3} \ten f)(\1^{\otimes 2} \ten a_{j,t}^\perp)\Big)\\
\!\!\!\! & = & \!\!\!\! (id^{\ten 3}\ten\tau_\Theta) \Big((\1\ten \chi_{\mathrm{B}}\ten \1^{\ten 2}) (id\ten  \sigma_{\Theta}\ten id)(\nabla_\Theta \ten id)(k)(\1^{\ten 3} \ten f)(\1^{\ten 2}\ten a_{j,t}^\perp) \Big)\\ [3pt]
\!\!\!\! & = & \!\!\!\! (id^{\ten 3}\ten \tau_\Theta) \big( \mathbf{K} \bullet (\1^{\ten 3} \ten f) \big)
\end{eqnarray*}
with $$\mathbf{K} = (\1\ten \chi_{\mathrm{B}}\ten \1^{\ten 2}) (id\ten  \sigma_{\Theta}\ten id)(\nabla_\Theta \ten id)(k) \bullet (\1^{\ten 2}\ten a_{j,t}^\perp)$$ in $\mathcal{L}(\mathbb{F}_n) \bar\ten ( \T(\R^n)\ten_\pi \T_\Theta\ten_\pi \T_\Theta)'$. Thus, $(\1\ten \chi_{\mathrm{B}}\ten \1) (\nabla\ten id_{\RR_\Theta}) (\Sigma)= T_{\mathbf{K}}(f)$. We turn to the proofs of Theorem 2.6, Proposition 2.15 and Remark 2.16 (as the generalizations of Theorem 2.6) in \cite{GPJP}, they show that the condition ii) in Theorem \ref{thm: cite} implies 
$$
\big\|T_{\mathbf{K}}(f) \big\|_{\mathcal{L}(\mathbb{F}_n) \bar\ten L_\infty(\R^n)\bar\ten \RR_\Theta}\leq C_n(\alpha,\beta)\frac{A_2}{\mathrm{R}}\| f\|_{\RR_\Theta},
$$
which is inequality (2.2) in \cite{GPJP}. Combining the calculations above, we deduce that condition ii) in Theorem \ref{thm: cite} is stronger than condition ii) in Theorem \ref{thm: CZ QE}. In conclusion,  the Calder\'on-Zygmund extrapolation on $\RR_\Theta$ that we obtain by applying Theorem \ref{Extrapolation} improves the corresponding result in \cite{GPJP}.

\subsection{Quantum Fourier multipliers}

We now refine our abstract result for locally compact quantum groups. We shall need some basic notions from the theory of quantum groups, details can be found in Kustermans/Vaes' papers \cite{KV1, KV2}. Let us consider a von Neumann algebra $\Mn$ equipped with a comultiplication map, a normal injective unital $*$-morphism $\Delta: \Mn \to \Mn \bar\otimes \Mn$ satisfying the coassociativity law $$(id_\Mn \otimes \Delta) \Delta = (\Delta \otimes id_\Mn)
\Delta.$$ Assume also the existence of two n.s.f weights $\psi$ and $\varphi$ on $\Mn$ such that 
$$
(id_\Mn \ten \psi) \Delta(a) = \psi(a) \1_\Mn \;\; \text{and }\;\; ( \varphi \ten id_\Mn) \Delta(a) = \varphi(a) \1_\Mn  \;\; \text{for }\,a\in \Mn_+.
$$
We call $\psi$ and $\varphi$ the \emph{left-invariant  Haar weight} and the \emph{right-invariant  Haar weight} on $\Mn$ respectively.  Then the  quadruple $\mathbb{G}=(\Mn, \Delta, \psi, \varphi)$ is called a (von Neumann algebraic) \emph{locally compact quantum group} and we write $L_\infty(\mathbb{G})$ for the quantum group von Neumann algebra $\Mn$. Using the Haar weights, one can construct an antipode $S$ on $\Mn$ which is a densely defined anti-automorphism on $\Mn$ satisfying the identity
$$
(id_\Mn \ten \psi) \big( (\1_\Mn \ten a^*) \Delta(b) \big) = S \big( (id_\Mn \ten \psi) \big( \Delta(a^*)(\1_\Mn \ten b) \big) \big).
$$

The comultiplication map $\Delta$ determines a multiplication  on the predual $L_1(\mathbb{G})$ given by convolution $\varphi_1 \star \varphi_2 (a) = (\varphi_1 \otimes \varphi_2) \Delta(a)$. The pair $(L_1(\mathbb{G}), \star)$ forms a Banach algebra. In what follows, if not specified otherwise, the quantum groups $\mathbb{G}$ we shall work with admit a tracial left-invariant  Haar weight $\psi$. The simplest model of noncommutative quantum groups are group von Neumann algebras $\V$ associated to discrete groups. If $\lambda$ is the left regular representation of $\G$, the comultiplication  is determined by $\Delta(\lambda(g)) = \lambda(g) \otimes \lambda(g)$. Its isometric nature follows from Fell's absorption principle and the convolution is abelian. The standard trace on $\V$ is a left and right-invariant Haar weight. Moreover, in this case, the antipode is bounded and $S(\lambda(g))=\lambda(g^{-1})$. 

A convolution semigroup of states is a family $(\phi_t)_{t \ge 0}$ of normal states on $L_\infty(\mathbb{G})$ such that $\phi_{t_1} \star \phi_{t_2} = \phi_{t_1+t_2}$.  The corresponding semigroup of completely positive maps is given by
$$S_{\Delta,t} (a) = (\phi_t \otimes id_{\mathbb{G}}) \circ \Delta(a).$$ When $\T_\Delta = (S_{\Delta,t})_{t \ge 0}$ is a Markov semigroup, we call it a \emph{convolution semigroup}.

\begin{lemma} \label{Tstar}
Let $\mathbb{G}$ be a locally compact quantum group equipped with a convolution semigroup of states $(\phi_t)_{t \ge 0}$. Then, $\T_\Delta = ((\phi_t \otimes id_{\mathbb{G}}) \circ \Delta)_{t \ge 0}$ is a Markov semigroup on $L_\infty(\mathbb{G})$ whenever
\begin{itemize}
\item[\emph{i)}]$\phi_t \circ S = \phi_t$ for all $t \ge 0$,

\item[\emph{ii)}] $S_{\Delta,t} (a) \to a$ as $t \to 0^+$ in the weak-$*$
topology of $L_\infty(\mathbb{G})$.
\end{itemize}
\end{lemma}

\dem Let us begin with the self-adjointness
\begin{eqnarray*}
\psi \big( a^* S_{\Delta,t}(b) \big) & = & \psi \big( a^* (\phi_t
\ten id_{\mathbb{G}}) \Delta(b) \big) \\ [5pt] & = & \phi_t \ten \psi \big(
(\1_{\mathbb{G}} \ten a^*) \Delta(b) \big) \\ [2pt] & = & \phi_t \Big(
(id_{\mathbb{G}} \ten \psi) \big( (\1_{\mathbb{G}} \ten a^*) \Delta(b) \big)
\\ & = & \phi_t \Big( S \underbrace{(id_{\mathbb{G}} \ten \psi) \big(
\Delta(a^*) (\1_{\mathbb{G}} \ten b) \big)}_\rho \Big).
\end{eqnarray*}
This means that $\psi \big( a^* S_{\Delta,t}(b) \big) = \phi_t
(S(\rho)) = \phi_t(\rho)$ and we get
\begin{eqnarray*}
\psi \big( a^* S_{\Delta,t}(b) \big) & = & \phi_t \ten \psi \big(
\Delta(a^*) (\1_{\mathbb{G}} \ten b) \big) \\ & = & \psi \big( (\phi_t \ten
id_{\mathbb{G}}) \circ \Delta(a^*) b \big) \ = \ \psi \big( S_{\Delta,t}(a)^* b
\big).
\end{eqnarray*}
The remainder properties are straightforward. Indeed, identity $S_{\Delta,t}(\1_{\mathbb{G}}) = \1_{\mathbb{G}}$ is obvious. The weak-$*$ convergence of the $S_{\Delta,t}(a)$'s as $t \to 0^+$ is assumed and the complete positivity is clear. The normality follows from the weak-$*$ continuity of $\phi_t$ and $\Delta$. Finally, the semigroup law easily follows from coassociativity. \fin

In what follows, we shall assume that the hypotheses of Lemma \ref{Tstar} hold. Let us fix a quantum group $\mathbb{G} = (\Mn, \Delta,\psi,\varphi)$ and consider a convolution semigroup $\T_\Delta$ associated to it. A Markov metric $$\Q = \big\{ (R_{j,t}, \sigma_{j,t}, \gamma_{j,t}) : \  j,t \in \mathbb{Z}_+ \times \mathbb{R}_+ \big\}$$ in $L_\infty({\mathbb{G}}) = \Mn$ associated to $\T_\Delta$ will be called an \emph{intrinsic Markov metric} when there exists an increasing family of projections $p_{j,t}$ in $L_\infty(\mathbb{G})$ such that the cpu maps take the form 
\begin{equation} \label{Eq-IMM}
R_{j,t}f = \frac{1}{\psi (p_{j,t})}(\psi \ten id_{\mathbb{G}}) \big( (p_{j,t} \ten \1_\mathbb{G}) \Delta ( f) \big).
\end{equation}
In other words, we use the algebraic skeleton $$\big( \Mn_\rho = \Mn_\pi, \rho_1, \rho_2, \mathsf{E}_\rho, q_{j,t} \big) = \big( L_\infty(\mathbb{G}) \bar\otimes L_\infty(\mathbb{G}), \1 \ten \cdot, \Delta, \psi \ten id_\mathbb{G}, p_{j,t} \ten \1_\mathbb{G} \big).$$ 

\begin{remark} \label{Rem-NoSizeTerm}
\emph{Assume that $$\gamma_{j,t} \in \R_+ \quad \mbox{and} \quad \gamma_{j,t}^2 \ge \frac{\psi(p_{j,t})}{\psi(p_{1,t})}.$$ Then, the term $|R_{j,t}f-M_t f|$ in the metric BMO norm satisfies for $M_t=R_{1,t}$ that
\begin{eqnarray*}
\big| R_{j,t}f  - M_t f \big|^2 & = & \Big| \frac{1}{\psi (p_{1,t})} (\psi \ten id_{\mathbb{G}}) \big( (\Delta (f)- \1 \ten R_{j,t} f) (p_{1,t} \ten \1) \big) \Big|^2 \\ & \leq & \frac{1}{\psi (p_{1,t})} (\psi \ten id_{\mathbb{G}}) \Big[ \big| (\Delta (f)-\1 \ten R_{j,t} f) (p_{1,t} \ten \1) \big|^2 \Big] \\ & = & \frac{\psi (p_{j,t})}{\psi (p_{1,t})} \big( R_{j,t}|f|^2 -|R_{j,t} f|^2\big)  \leq \gamma_{j,t}\big( R_{j,t}|f|^2 -|R_{j,t} f|^2\big) \gamma_{j,t}.
\end{eqnarray*}
According to Theorem \ref{MetricBMOInterp}, this yields
\begin{equation}\label{eq: BMO_Q}
\| f \|_{{\rm BMO}^c_{\T_\Delta}} \, \lesssim \, \|f\|_{\mathrm{BMO}_\Q^c \, }\lesssim \, \sup _{t>0}\sup_{j\geq 1} \Big\| \Big(R_{j,t}|f|^2 -|R_{j,t} f|^2\Big)^\frac{1}{2}\Big\|_{L_\infty(\mathbb{G})}.
\end{equation}}
\end{remark}

Additionally, we may consider transferred Markov metrics in other von Neumann algebras. Consider a convolution semigroup of states $(\phi_t)_{t\geq 0}$ on a locally compact quantum group $L_\infty(\mathbb{G})$. A corepresentation $\pi : \M \rightarrow L_\infty(\mathbb{G}) \bar\ten \M$ is a normal injective $*$-representation satisfying the identity
$$
(id_{\mathbb{G}}\ten \pi)\circ \pi =(\Delta \ten id_\M)\circ \pi.
$$
Every such $\pi$ yields a \emph{transferred convolution semigroup} $\T_\pi =(S_{\pi, t})_{t\geq 0}$ with  
$$
S_{\pi,t}: \M \to \M,$$ $$S_{\pi,t}f =(\phi_t \ten id_\M)\circ \pi (f).
$$ 

\begin{lemma}\label{lem: trans convolution semigroup}
Assume that
\begin{itemize}
\item $\tau(S_{\pi,t}(f_1)^* f_2) = \tau(f_1^* S_{\pi,t}(f_2))$,

\item $S_{\pi,t}f \to f$ as $t \to 0^+$ in the weak-$*$ topology
of $\M$.
\end{itemize}
Then $\T_\pi$ defines a Markov semigroup on $\M$ such that $\pi \circ S_{\pi,t}=(S_{\Delta,t} \ten id_{\M})\circ \pi$.
\end{lemma}

\dem It is easy to check that $S_{\pi,t}$ is cpu and the normality follows from the weak-$*$ continuity of $\phi_t$ and $\pi$. Hence, it remains to show the identity $\pi \circ S_{\pi,t}=(S_t\ten id_{\mathbb{G}})\circ \pi$ and the semigroup law. We first observe that $\pi (\phi_t \otimes id_\M)=(\phi_t \ten id_{\mathbb{G}} \ten id_\M)(id_{\mathbb{G}}\ten \pi)$ as maps on $L_\infty(\mathbb{G})\bar\ten \M$. Indeed, by weak-$*$ continuity, it suffices to test the identity on elementary tensors $n\ten m$, for which the identity is trivial. Therefore, we have
\begin{eqnarray*}
(S_{\Delta,t} \ten id_{\M}) \pi &=& (\phi_t \ten id_{\mathbb{G}} \ten id_\M)(\Delta \ten id_\M) \pi \\
& = & (\phi_t \ten id_{\mathbb{G}} \ten id_\M)(id_{\mathbb{G}} \ten \pi )\pi \\
& = & \pi (\phi_t \otimes id_\M) \pi =\pi S_{\pi,t}.
\end{eqnarray*}
For the semigroup law we note that
\begin{eqnarray*}
\hskip-3pt S_{\pi,t_1}S_{\pi,t_2} \!\!\! &= & \!\!\! (\phi_{t_1} \otimes id_\M) (\phi_{t_2} \ten id_{\mathbb{G}} \ten id_\M)(id_{\mathbb{G}}\ten \pi) \pi \\
\!\!\! & = & \!\!\!  (\phi_{t_1} \otimes id_\M) (\phi_{t_2} \ten id_{\mathbb{G}} \ten id_\M)(\Delta \ten id_\M) \pi \\
\!\!\! & = & \!\!\! (\phi_{t_2}\ten \phi_{t_1} \otimes id_\M)(\Delta \ten id_\M) \pi = (\phi_{t_2} \star \phi_{t_1} \ten id_\M) \pi = S_{\pi, t_1+t_2}. \hskip15pt \square
\end{eqnarray*} 

In the sequel, we shall assume that the assumptions in Lemma \ref{lem: trans convolution semigroup} hold. Intrinsic Markov metrics on $L_\infty(\mathbb{G})$ yield \emph{transferred Markov metrics} on $\M$ associated to the transferred convolution semigroup $\T_\pi$.  Indeed, given any intrinsic Markov metric $\Q = \{ (R_{j,t}, \sigma_{j,t}, \gamma_{j,t}) \}$ in $\mathbb{G}$ with cpu maps $R_{j,t}$ given by \eqref{Eq-IMM}, the transferred cpu maps $R_{\pi,j,t}$ are given by 
$$
R_{\pi,j,t} f =\frac{1}{\psi (p_{j,t})}(\psi \ten id_\M) \big( (p_{j,t} \ten \1)\pi (f) \big).
$$
It is easy to check that $\pi \circ R_{\pi,j,t}= (R_{j,t}\ten id_\M) \circ \pi$. Assume in addition that $\sigma_{j,t} \in \R_+$. Then, arguing as we did before Corollary \ref{CorTT} for the corepresentation $u$ of $\R$ in $\mathcal{B}(\ell_2)$, we get a Markov metric in $\M$ $$\Q_\pi = \big\{ (R_{\pi,j,t}, \sigma_{j,t}, \gamma_{j,t} ) : j \in \Z_+, t \in \R_+ \big\}.$$ Let  $\alpha: \N \rightarrow \N$ be a strictly increasing function with $\alpha(j) > j$. This  Markov metric is called \emph{$\alpha$-doubling} if there exists some constant $c_\alpha$ such that $\psi (q_{\alpha(j),t}) \leq c_\alpha \psi (q_{j,t}).$

\begin{remark}
\emph{In what follows, we impose our Markov metrics to be $\alpha$-doubling for some function $\alpha: \N \to \N$, to satisfy $\sigma_{j,t} \in \R_+$ as well as the condition in Remark \ref{Rem-NoSizeTerm}. Altogether, this allows to eliminate the size CZ condition and reduce the number of analytic and smoothness CZ conditions to be checked for both the intrinsic Markov metric and the transferred one.}
\end{remark}

Observe that the transferred formulation above includes the intrinsic formulation by taking $(\M,\pi) = (\mathbb{G},\Delta)$. Let us now state the corresponding Calder\'on-Zygmund theory. Given $\A_\M$ a weakly dense $*$-subalgebra of $\M$, let $T$ be a (not necessarily bounded) operator $T: \A_\M \to \M$. We say $T$ is a  \emph{transferred map} if there exists  an amplification map $$\widehat{T}: \mathcal{D} \subset L_\infty(\mathbb{G}) \bar\otimes \M \to
L_\infty(\mathbb{G}) \bar\otimes \M$$ satisfying the identity 
\begin{equation} \label{Eq-TransferQFM}
\pi \circ T =\widehat{T} \circ \pi _{\mid_ {\A_\M}}.
\end{equation}
Again, $\mathcal{D}$ is a weakly dense $*$-subalgebra for which $\pi (\A_\M) \subset \mathcal{D}$. In the case $(\M,\pi) = (\mathbb{G},\Delta)$, we can always take the amplification $T \ten id_{\M}$ and condition above just imposes that $T$ is a quantum Fourier multiplier. In the following theorem, we provide sufficient conditions on the amplification map to make a given transferred CZO $T$ bounded from $\A_\M$ to $\mathrm{BMO}_{\T_\pi}^c$. 

\begin{theorem} \label{CZ1}
Let $$\pi: \M \to L_\infty(\mathbb{G}) \bar\ten \M$$ be a corepresentation of a locally compact quantum group $\mathbb{G}$ in a semifinite von Neumann algebra $(\M,\tau)$. Assume that $L_\infty(\mathbb{G})$ comes equipped with an $\alpha$-doubling intrinsic Markov metric $\Q$ determined by an increasing family of projections $p_{j,t}$ as above. Then, a transferred map $T$ $($with amplification for which \eqref{Eq-TransferQFM} holds$)$ will be bounded from $\A_\M$ to $\mathrm{BMO}_{\S_\pi}^c$ provided$\hskip1pt :$
\begin{itemize}
\item[\emph{i)}] $\widehat{T}: L_2^c(\mathbb{G}) \bar\otimes \M \to L_2^c(\mathbb{G})
\bar\otimes \M$ is bounded,

\vskip3pt

\item[\emph{ii)}] $\displaystyle \frac{(\psi\ten \psi \ten id_\M)}{\psi(p_{j,t})^2} \Big( (p_{j,t} \otimes p_{j,t}\otimes \1_{\M}) \big| \delta_{\mathbb{G}} \big(
\widehat{T}(\pi (f)p_{\alpha(j),t}^{\perp})\big)\big|^2 \Big) 
\lesssim  \|f\|_\M^2$.
\end{itemize}
\end{theorem}

\dem We use the algebraic skeleton $$\big( \M, \Mn_\rho = \Mn_\pi, \rho_1, \rho_2, \mathsf{E}_\rho, q_{j,t} \big) = \big( \M, L_\infty(\mathbb{G}) \bar\otimes \M, \1 \ten \cdot, \pi, \psi \ten id_\mathbb{G}, p_{j,t} \ten \1_\mathbb{G} \big).$$ Identity \eqref{Eq-TransferQFM} is the compatibility condition \eqref{Amplification}. Let us justify the algebraic conditions. The second one is trivial since both $\mathsf{E}_\rho(q_{j,t})$ and $\rho_1(\gamma_{j,t})$ belong to $\R_+$ in this case. For the first one, consider the product $\bullet$ in $L_\infty(\mathbb{G}) \bar\otimes \M_\mathrm{op}$. Then, we just observe that
\begin{eqnarray*}
\mathsf{E}_\rho(q_{j,t} |\xi|^2 q_{j,t}) & = & (\psi \ten id_\M) \big( (p_{j,t} \ten \1) \xi^* \xi \big) \\
& = & (\psi \ten id_\M) \big( \xi \bullet (p_{j,t} \ten \1) \bullet \xi^* \big) \\ & \le & (\psi \ten id_\M) ( \xi \bullet \xi^* ) \ = \ (\psi \ten id_\M) (\xi^* \xi) \ = \ \mathsf{E}_\rho (|\xi|^2).
\end{eqnarray*} 
Define the amplifications
$$
\widehat{R}_{\pi ,j,t}:L_\infty(\mathbb{G}) \bar\ten \M \ni \xi \mapsto \frac{1}{\psi (p_{j,t})} (\psi \otimes id_{\M}) \big( (p_{j,t}\otimes  \1_{\M})\xi \big) \in \M.
$$
Consider also the cpu maps $$\Phi_{j,t}: L_\infty(\mathbb{G})  \bar\ten L_\infty(\mathbb{G}) \bar\ten \M \rightarrow \M,$$
$$
\Phi_{j,t} (\eta)= \frac{1}{\psi (p_{j,t})^2} (\psi \otimes \psi \otimes id_{\M}) \big( (p_{j,t} \otimes p_{j,t}\otimes \1_{\M})\eta \big).
$$
Recalling that $\delta_\mathbb{G}(x) = x \ten \1 - \1 \ten x$, the identity
$$
\Phi_{j,t} (|\delta_{\mathbb{G}} (\xi)|^2)=2 \widehat{R}_{\pi ,j,t} (|\xi |^2)- 2 \big|\widehat{R}_{\pi ,j,t} (\xi )\big|^2.
$$
is straightforward. This readily implies the first analytic condition. On the other hand, since the auxiliary Markov metric is $\alpha$-doubling, the second analytic condition reduces to note that  
$$
(\psi \otimes id_{\M})(q_{\alpha(j),t}\ten \1_\M) \leq c_\alpha (\psi \otimes id_{\M})(q_{j,t}\ten \1_\M). 
$$
Thus, according to inequalitty \eqref{eq: BMO_Q}, the assertion follows from Theorem \ref{Extrapolation}.\fin

\begin{rmk}
As noticed, the main particular case of Theorem \ref{CZ1} arises for $(\M,\pi) = (\mathbb{G}, \Delta)$ with amplification $T \ten id_\mathbb{G}$. Condition \eqref{Eq-TransferQFM} becomes the identity $$\Delta \circ T =(T \ten id_\mathbb{G}) \circ \Delta.$$ In other words, these are translation invariant CZ operators. We also call them quantum Fourier multipliers in this paper and it can be checked, as expected, that these maps are of convolution type in the sense that there exists a kernel $k$ affiliated to $L_\infty(\mathbb{G})$ so that $$Tf = k \star f = ( id_{\mathbb{G}} \ten \psi ) \big( \Delta(k) (\1_\mathbb{G} \ten Sf) \big).$$ In this particular case, it is not difficult to prove that our conditions in Theorem \ref{CZ1} reduce to those in Theorem B2 from the Introduction. Of course, Theorem \ref{CZ1} also applies as well for nonconvolution CZ operators on quantum groups, or even for transferred forms of them to other von Neumann algebras $\M$.       
\end{rmk}

\begin{rmk}
 One may consider twisted  convolution CZO's on quantum groups applying Theorem \ref{CZ1}. As an illustration, assume that $\G \curvearrowright L_\infty(\mathbb{G})$ by a trace preserving action $\alpha$ and that $\mathbb{G}$ is a quantum group satisfying $(\alpha_g \otimes \alpha_g) \Delta = \Delta \alpha_g$ for all $g \in \G$. This property is quite natural in the commutative case, where quantum groups come from locally compact  groups and $\alpha$ is typically implemented by a measure preserving map $\beta$. Note that the underlying Haar measure is translation invariant and the condition above just imposes that $\beta$ is an homomorphism. Let us see what  we get for a map $$\summ_g f_g \lambda(g) \mapsto \summ_g T_g(f_g) \lambda(g),$$ where the $T_g$'s are normal convolution maps on $L_\infty(\mathbb{G})$. Assume $L_\infty(\mathbb{G})$ comes equipped with a convolution $\G$-equivariant semigroup $\T_\Delta$ which admits a $\eta$-doubling intrinsic Markov metric. Then, we get a bounded map $L_\infty(\mathbb{G}) \rtimes. \G \to \mathrm{BMO}_{\T_\rtimes}^c$ when the following conditions hold:
\begin{itemize}
\item[i)] We have a bounded map $$\hskip15pt L_2^c(\mathbb{G}) \bar\otimes \mathcal{B}(\ell_2(\G)) \ni \xi \mapsto \big( T_{gh^{-1}}  \big) \bullet \xi \in L_2^c(\mathbb{G}) \bar\otimes \mathcal{B}(\ell_2(\G)),$$
where $\bullet $ stands once more for the generalized Schur product of matrices.

\vskip3pt

\item[ii)] Letting $\RR = L_\infty(\mathbb{G}) \bar\ten \mathcal{B}(\ell_2(\G))$ and $\Psi(\xi) = \sum_{g,h}  (\alpha_{g^{-1}}) \bullet  (T_{gh^{-1}}) \bullet (\alpha_g) \bullet \xi,$
$$\hskip15pt \frac{(\psi\ten \psi \ten id_{\mathcal{B}(\ell_2(\G))})}{\psi(p_{j,t})^2} \Big( (p_{j,t} \otimes p_{j,t}\otimes \1) \big| \delta_{\mathbb{G}} \big( \Psi( \xi q_{\eta(j),t}^{\perp}\big)\big|^2 \Big) \lesssim  \|\xi \|_\RR^2.$$
\end{itemize}
\end{rmk}

\begin{remark}
\emph{All our results in this paragraph impose the additional assumption that our quantum groups admit a tracial Haar weight. We believe however that our results can be extended to the general non-tracial case. We leave this generalization  open to the interested reader.}
\end{remark}

\section{{\bf  Noncommutative transference}}

Originally motivated by Cotlar's paper \cite{Co} and the method of rotations, Calder\'on developed a circle of ideas \cite{C} which was called the \emph{transference method} after the systematic study of Coifman/Weiss in their monograph \cite{CW}. The fundamental work of K. de Leeuw \cite{dL} also had a big impact in this line of research. Let us consider an amenable locally compact group $\G$ with left Haar measure $\mu$, a $\sigma$-finite measure space $(\Omega,\nu)$ and a uniformly bounded representation $\beta: \G \to \mathcal{B}(L_p(\Omega))$. Roughly, Calder\'on's transference is a technique which allows to transfer the $L_p$ boundedness of a convolution operator $f \mapsto k \star f$
on $L_p(\G)$ to the corresponding transferred operator on $L_p(\Omega)$ $$V \! f (w) = \int_\G k(g) \, \beta_{g^{-1}}f(w) \,
d\mu(g),$$ for some compactly supported kernel $k$ in $L_1(\G)$. A case by case limiting procedure also allows to consider more general (singular) kernels. In the rest of this section we shall develop a noncommutative form of Calder\'on-Coifman-Weiss technique.

Our first task is to clarify what we mean by \lq representation\rq${}$ and  \lq
amenable\rq${}$  in the context of quantum groups. Using the
 commutative locally compact quantum group $L_\infty(\G)$ as above, a
representation $\beta: \G \to \mathrm{Aut}(\M)$ induces a
$*$-representation $\pi_\beta: \M \to L_\infty(\G; \M)$ by
$$\pi_\beta f(g) = \beta_{g^{-1}} f.$$ Note that we have
\begin{eqnarray*}
(id_\G \otimes \pi_\beta) (\pi_\beta f)(g,h) & = &
\pi_\beta(\beta_{g^{-1}} f)(h) \ = \ \beta_{h^{-1}} \beta_{g^{-1}}
f \\ & = & \beta_{(gh)^{-1}} f \ = \ (\Delta_\G \otimes id_\M)
(\pi_\beta f)(g,h).
\end{eqnarray*}
Given a semifinite von Neumann algebra $(\M, \tau)$ and a  locally compact quantum group $\mathbb{G}$, this leads us  to consider corepresentations $\pi : \M \mapsto L_\infty(\mathbb{G})\bar \ten \M$  satisfying $(id_{\mathbb{G}}\ten \pi)\circ \pi =(\Delta \ten id_\M)\circ \pi$. Note that comultiplication is a corepresentation by coassociativity.  To show what we mean by \lq uniformly bounded\rq${}$, let us go back to our motivating example $\beta: \G \to \mathrm{Aut}(\M)$, where we take $\M = L_\infty(\Omega)$ for some $\sigma$-finite measure space $(\Omega, \nu)$. In the classical case
$$\|\beta_g f\|_p \sim \|f\|_p \quad \mbox{for all} \quad g \in \G$$ 
up to an absolute constant independent of $f,g$. We say that a corepresentation $\pi: \M \to L_\infty(\mathbb{G}) \bar\ten \M$ is \emph{uniformly bounded in $L_p(\M)$} if for any $f \in \M \cap L_p(\M)$ we have $$\frac{1}{c_\pi} \|f\|_{L_p(\M)}^p \le (id_{\mathbb{G}} \ten \tau) \big( |\pi(f)|^p \big) \le c_\pi \|f\|_{L_p(\M)}^p$$ for some absolute constant $c_\pi$ independent of $f$. Note that our notion again reduces to the classical one on $L_\infty(\G)$. Note also that, since $|\pi(f)|^p = \pi(|f|^p)$, our definition reduces to the $p$-independent condition $$\frac{1}{c_\pi} \|f\|_{L_1(\M)} \le (id_{\mathbb{G}} \ten \tau) \big( \pi(f) \big) \le c_\pi \|f\|_{L_1(\M)} \quad \mbox{for all} \quad f \in \M_+ \cap L_1(\M).$$
Now we introduce what we mean by an \lq amenable\rq${}$ quantum group.  We say that $\mathbb{G}$ satisfies \emph{F\o lner's condition} if for every projection $q \in L_1(\mathbb{G})$ and every $\varepsilon > 0$, there exists
two projections $q_1, q_2 \in L_1(\mathbb{G})$ such that $$\Delta(q_1) (q \ten q_2) = q \ten q_2 \quad \mbox{and} \quad
\psi(q_1) \le (1 + \varepsilon) \, \psi(q_2).$$ In the standard example for a locally compact group $\G$, where $(L_\infty(\mathbb{G}),\psi)$ is $L_\infty(\G)$ equipped with the left Haar measure $\mu$ and $\Delta$ is given by $\Delta_\G(\xi)(g,h) = \xi(gh)$ the classical comultiplication, it turns out that $\mathbb{G}$ is amenable iff $\G$ is an amenable group. Indeed, our notion can be rephrased in this case by saying that for any compact set $\mathrm{K}$ in $\G$ and any $\varepsilon > 0$, there exists a neighborhood of the identity $\mathrm{W}$ of finite measure such
that $$\mu(\mathrm{KW}) \le (1 + \varepsilon) \, \mu(\mathrm{W}),$$ which corresponds to $(q,q_1,q_2) =
(\chi_\mathrm{K}, \chi_\mathrm{KW}, \chi_\mathrm{W})$ in our formulation. This is exactly the classical characterization of
amenability, known as F\o lner's condition, used by Coifman and Weiss in \cite{CW}. Given an amenable locally compact group $\G$ with left Haar measure $\mu$, it is clear that $L_\infty(\G,\mu)$ with its natural quantum group structure is amenable. On the other hand, as expected, any compact quantum group is amenable just by taking $q_1 = q_2 = \1_{\mathbb{G}}$.

Assume that $\mathbb{G}$ admits a corepresentation $\pi: \M \to L_\infty(\mathbb{G}) \bar\ten \M$. Given $\A_\M$ a weakly dense $*$-subalgebra of $\M$, we say that a linear operator $V: \A_\M \to \M$ is a \emph{transferred convolution map} if there exists $\Phi: \mathcal{D} \subset L_\infty(\mathbb{G}) \bar\ten \M \to L_\infty(\mathbb{G}) \bar\ten \M$, an auxiliary convolution map such that $\pi \circ V = \Phi \circ \pi_{\mid_{\A_\M}}$. The classical transferred operator $$V = \int_\G k(g) \, \beta_{g^{-1}}f(w) \, d\mu(g)$$ comes from $$\Phi(\xi)(g,w) = \int_\G k(h) \, \xi(hg,w) \, d\mu(h) = (\varphi \ten id_\G \ten id_\Omega) (\Delta_\G \ten id_\Omega).$$ If $\pi_\beta f(g) = \beta_{g^{-1}} f$ denotes the corresponding corepresentation, we may then apply the identities in the proof of Lemma \ref{lem: trans convolution semigroup} again to deduce the following identities $$\Phi \circ \pi_\beta = (\varphi \ten id_\G \ten id_\Omega) (\Delta_\G \ten id_\Omega) \pi_\beta = (\varphi \ten id_\G \ten id_\Omega) (id_\G \ten \pi_\beta) \pi_\beta = \pi_\beta (\varphi \ten id_\Omega)\circ \pi_\beta.$$ By injectivity of $\pi_\beta$, we must have $$V \! f(w) = (\varphi \ten id_\Omega) \, \pi_\beta f(w) = \int_\G k(g) \, \beta_{g^{-1}} f (w) \, d\mu(g)$$ as expected. This shows how we recover the classical construction.

Let us now settle the framework for our transference result. Assume that $\mathbb{G}$ is amenable and consider $\pi: \M \to L_\infty(\mathbb{G}) \bar\ten \M$ a uniformly bounded corepresentation in some noncommutative measure space $(\M,\tau)$.  We say that $T: L_p(\mathbb{G}) \to L_p(\mathbb{G})$ is a convolution map with finitely supported $L_1$ kernel when the map $T$ has the form $T = (\phi \ten id_{\mathbb{G}}) \circ \Delta$ for some functional $\phi = \psi(d \, \cdot)$, with $d$ an element in $L_1(\mathbb{G})$ whose left support $q$ satisfies $\psi(q) < \infty$. In the commutative case, this is the kind of operators which are transferred. Roughly, the goal is to show how a limit operator $T = \lim_\gamma T_\gamma$ of such maps which is bounded on $L_2(\mathbb{G})$ and  $L_\infty(\mathbb{G}) \to \mathrm{BMO}_{\T}$   can be transferred under suitable conditions to a  bounded map on $L_p(\M)$.

\begin{remark}
\emph{Young's inequality extends to this setting as
$$\|d \star f\|_p \, \mbox{\lq=\rq} \, \|(\phi \ten id_{\mathbb{G}})
\Delta(f)\|_p \le 4 \, \|d\|_1 \|f\|_p,$$ where $\phi = \psi(d
\, \cdot)$ and $1\leq p\leq\infty$. Indeed, when $d$ and $f$ are positive the inequality
holds with constant $1$. This can by justified by interpolation.
When $p=1$ we use Fubini and the left-invariance of $\psi$, while
for $p=\infty$ it follows from the fact that $(\phi \ten
id_{\mathbb{G}}) \Delta$ is a positive map with $\1_{\mathbb{G}} \mapsto \psi(d)$.
In the general case, we split $d,f$ into their positive parts and
obtain the constant $4$. In fact, the same argument still holds
after matrix amplification and we deduce that $(\phi \ten
id_{\mathbb{G}}) \Delta$ is completely bounded on $L_p(\mathbb{G})$ with cb-norm
$4 \|d\|_1$. This is however not enough for transference, since
the norms $\|d_\gamma\|_1$ might not be uniformly bounded.}
\end{remark}

\begin{theorem} \label{TransfTh}
Let $\mathbb{G}$ be an amenable quantum group and consider a uniformly bounded corepresentation 
$\pi: \M \to L_\infty(\mathbb{G}) \bar\ten \M$ in some
noncommutative measure space $(\M,\tau)$. Let $T: L_2(\mathbb{G}) \to
L_2(\mathbb{G})$ be a bounded map and assume that $$(T \ten id_\M) =
\mathrm{SOT}-\lim_\gamma (T_\gamma \ten id_\M)$$ for some net
$T_\gamma = (\phi_\gamma \ten id_{\mathbb{G}}) \circ \Delta$ of
convolution maps with finitely supported $L_1$ kernels and such that $\lim _\gamma \|T_\gamma \|_{\mathcal{B}(L_2(\mathbb{G}))}\leq \|T \|_{\mathcal{B}(L_2(\mathbb{G}))}$. Then, the
net of transferred operators $V_\gamma = (\phi_\gamma \ten
id_\M) \circ \pi$ satisfies the inequalities
$$\|V_\gamma\|_{\mathcal{B}(L_2(\M))} \le c_\pi
\|T_\gamma\|_{\mathcal{B}(L_2(\mathbb{G}))}.$$ We thus find a
$\mathrm{WOT}$-cluster point $V$ satisfying $\|V\|_{\mathcal{B}(L_2(\M))} \le c_\pi
\|T\|_{\mathcal{B}(L_2(\mathbb{G}))}.$
\end{theorem}

\dem Note that we have $$\pi V_\gamma = (\phi_\gamma \ten
id) (id_{\mathbb{G}} \ten \pi) \pi = (\phi_\gamma \ten
id) (\Delta \ten id_\M) \pi = (T_\gamma \ten
id_\M) \pi.$$
Hence, the uniform boundedness of $\pi$ yields
$$\frac{1}{c_\pi }\|V_\gamma  f\|_2^2 \ \le  (\rho \ten \tau) \big( | \pi
V_\gamma  (f)|^2 \big) \ =  (\rho \ten \tau) \big(
|(T_\gamma  \ten id_\M) \pi(f)|^2 \big)$$ 
for any state $\rho$ on
$L_\infty(\mathbb{G})$. On the other hand, if $\phi_\gamma = \psi(d_\gamma \,
\cdot)$ and $q_\gamma$ denotes the left support of $d_\gamma$, we
know from the amenability assumption that for any $\varepsilon >
0$ we may find projections $q_{1\gamma}$ and $q_{2\gamma}$ such
that $$\Delta(q_{1\gamma}) (q_\gamma \ten q_{2\gamma}) = q_\gamma
\ten q_{2\gamma} \quad \mbox{and} \quad \psi(q_{1\gamma}) \le (1 +
\varepsilon) \, \psi(q_{2\gamma}).$$ Taking $\rho =
\psi(q_{2\gamma} \, \cdot) / \psi(q_{2\gamma})$, we obtain the
inequality 
$$\frac{1}{c_\pi} \|V_\gamma  f\|_2^2 \ \le \ \frac{(\psi
\ten \tau)}{\psi(q_{2\gamma})} \Big( \Big| (\phi \ten
id) \big( (\Delta \ten id_\M) \pi(f) (\1_{\mathbb{G}}
\ten q_{2\gamma} \ten \1_\M) \big) \Big|^2 \Big) $$ 
since $\rho$ is
supported by $q_{2\gamma}$. Moreover, $d_\gamma \ten q_{2\gamma}$
is supported on the left by $q_\gamma \ten q_{2\gamma}$ and
amenability provides $d_\gamma \ten q_{2\gamma} =
\Delta(q_{1\gamma}) (d_\gamma \ten q_{2\gamma})$. Once we have
created $\Delta(q_{1\gamma})$, we can eliminate $q_{2\gamma}$.
Altogether gives $$\frac{1}{c_\pi} \|V_\gamma  f\|_2^2 \ \le \
\frac{1}{\psi(q_{2\gamma})} \, \Big\| (T_\gamma  \ten id) \big(
\pi(f) (q_{1\gamma} \ten \1_\M) \big) \Big\|_{L_2(L_\infty(\mathbb{G}) \bar\ten
\M)}^2.$$ 
Now we use the $L_2$ boundedness of $T_\gamma$ and
uniform boundedness of $\pi$ to conclude
\begin{eqnarray*}
\frac{1}{c_\pi} \|V_\gamma f\|_2^2 & \le &
\frac{1}{\psi(q_{2\gamma})} \,
\|T_\gamma \|_{\mathcal{B}(L_2(\mathbb{G}))}^2 \, \psi \Big( (q_{1\gamma}
\ten \1_\M) (id_{\mathbb{G}} \ten \tau)  \big( |\pi(f)|^2 \big) \Big) \\ &
\le & \frac{c_\pi}{\psi(q_{2\gamma})} \,
\|T_\gamma \|_{\mathcal{B}(L_2(\mathbb{G}))}^2 \, \psi (q_{1\gamma}) \,
\|f\|_2^2 \ \le \ c_\pi \, (1 + \varepsilon) \,
\|T_\gamma \|_{\mathcal{B}(L_2(\mathbb{G}))}^2 \, \|f\|_2^2.
\end{eqnarray*}
Letting $\varepsilon \to 0$, we prove the inequality $$\|V_\gamma \|_{\mathcal{B}(L_2(\M))} \ \le \ c_\pi \,
\|T_\gamma \|_{\mathcal{B}(L_2(\mathbb{G}))}.$$ Since $T$ is bounded on $L_2(\mathbb{G})$ and $\lim _\gamma  \|T_\gamma \|_{\mathcal{B}(L_2(\mathbb{G}))}\leq \|T \|_{\mathcal{B}(L_2(\mathbb{G}))}$, the operators $V_\gamma$ are eventually in a ball of radius $c_\pi (1 + \delta) \|T\|_{\mathcal{B}(L_2(\mathbb{G}))}$ for any $\delta > 0$. The closure of such ball is weak operator compact and thus we find our cluster point. \fin

We now study $L_\infty \to \mathrm{BMO}$ transference and then
interpolate/dualize to obtain $L_p$-transference. This approach seems to be new even in the classical theory and where our semigroup formulation becomes an essential ingredient.

\begin{corollary} \label{TransfCor}
Let $\mathbb{G}$ be a compact $($hence amenable$)$ quantum group equipped with a uniformly bounded corepresentation $\pi: \M \to L_\infty(\mathbb{G}) \bar\ten \M$. Let $(\phi_t)_{t \ge 0}$ be a convolution semigroup of states on $L_\infty(\mathbb{G})$, giving rise to Markov semigroups $\T_\Delta$ on $(\mathbb{G},\psi)$ and $\T_\pi$ on $(\M,
\tau)$. Let $T = \mathrm{SOT}-\lim_\gamma T_\gamma$ be as above and take $\A_\mathcal{M} = \mathcal{M} \cap L_2(\mathcal{M})$. Then, if $T: L_\infty(\mathbb{G}) \to \mathrm{BMO}_{\T_\Delta}$ is completely bounded, we find
that
$$V = \mathrm{WOT}-\lim_\gamma V_\gamma: \A_\M \to \mathrm{BMO}_{\T_\pi}$$ 
is completely bounded. Moreover, if $T_\pi$ is regular, the complete boundedness of $J_pV \hskip-2pt : L_p(\M) \to L_p(\M)$ follows for every $2 < p < \infty$ by interpolation. In addition the complete boundedness of $VJ_p: L_p(\M) \to L_p(\M)$ for $1 < p < 2$ holds under the same assumptions for $T^*$.
\end{corollary}

\dem By uniform boundedness of $\pi$ we have $$\Big\| (id_{\mathbb{G}} \ten
\tau) \big( |\pi(f)|^2 \big)^\frac12 \Big\|_{L_\infty(\mathbb{G})} \le c_\pi
\|f\|_2,$$ which implies that $\pi: L_2(\M) \to L_\infty(\mathbb{G}) \bar\ten
L_2^c(\M)$ is bounded by $c_\pi$. According to the finiteness of
$L_\infty(\mathbb{G})$, we deduce that in fact $\pi: L_2(\M) \to L_2(L_\infty(\mathbb{G})\bar\ten
\M)$ is still bounded with the same norm. This proves that $$\pi V
= \mathrm{WOT}-\lim_\gamma \pi V_\gamma = \mathrm{WOT}-\lim_\gamma
(T_\gamma \ten id_\M) \pi = (T \ten id_\M) \pi.$$ In particular,
$\pi V = (T \ten id_\M) \pi$ over $\A_\M$ and identity $\pi
S_{\pi,t} = (S_t\ten id_\M) \pi$ yields
\begin{eqnarray*}
\|Vf\|_{\mathrm{BMO}_{\T_\pi}^c} \!\!\! & = & \!\!\! \sup_{t \ge 0} \Big\|
S_{\pi,t} |Vf|^2 - |S_{\pi,t}Vf|^2 \Big\|_\M^\frac12 \\
\!\!\! & = & \!\!\! \sup_{t \ge 0} \Big\| \pi S_{\pi,t} |Vf|^2 - |\pi
S_{\pi,t}Vf|^2 \Big\|_{L_\infty(\mathbb{G}) \bar\ten \M}^\frac12 \\ [7pt] \!\!\! & = & \!\!\!
\big\| (T \ten id_\M) \pi (f) \big\|_{\mathrm{BMO}_{\T}^c}
\le \|T\|_{cb} \|\pi(f)\|_{L_\infty(\mathbb{G}) \bar\ten \M} = \|T\|_{cb}
\|f\|_\M
\end{eqnarray*}
for $f \in \A_\M$. Since the same inequality holds after matrix
amplification, we deduce that $V: \A_\M \to
\mathrm{BMO}_{\T_\pi}^c$ is completely bounded with cb-norm $\le
\|T\|_{cb}$. The row case is similar because $$\pi V^\dag = (\pi
V)^\dag = (T \pi)^\dag = T^\dag \pi.$$ The assertions on
$L_p$ boundedness follow as usual from Theorem
\ref{Interpolation}. \fin

\begin{rmk}
Under the above assumptions, we see that for $V = \mathrm{WOT}-\lim_\gamma V_\gamma$ we can find the concrete form of its amplification map $\Phi$ defined on $L_\infty(\mathbb{G}) \bar\ten \M$. In this case, by applying  Theorems \ref{CZ1} to $\Phi = T\ten id_\M$, we get Calder\'on-Zygmund extrapolation for the transferred convolution map $V$ on $\M$. 
\end{rmk}

\noindent \textbf{Acknowledgement.} Mei is partially supported by NSF grant DMS1700171. Javier Parcet is supported  by the Europa Excelencia Grant MTM2016-81700-ERC and the CSIC Grant PIE-201650E030. Javier Parcet and Runlian Xia are supported by ICMAT Severo Ochoa Grant SEV-2015-0554 (Spain).

\bibliographystyle{amsplain}

\begin{thebibliography}{99}

\bibitem {B} D. Bakry, \'Etude des transformations de Riesz dans les vari\'et\'es riemanniennes \`a courbure de Ricci minor\'ee. S\'eminaire de Probabilit\'es XXI. Lecture Notes in Math. \textbf{1247} (1987), 137-172. 

\bibitem {BK} S. Blunck and P. Kunstmann, Calder\'on-Zygmund theory for non-integral operators and the $H^\infty$--functional calculus. Rev. Mat. Iberoamericana \textbf{19} (2003), 919-942.

\bibitem {Bo3}  J. Bourgain, Vector valued singular integrals and the $H^1$-BMO duality. Probability Theory and Harmonic Analysis. (Eds. Chao and Woyczynski) Decker (1986), 1-19.

\bibitem {Ca} L. Cadilhac, Weak boundedness of Calder\'on-Zygmund operators on noncommutative $L_1$-spaces. J. Funct. Anal. \textbf{274} (2018), 769-796.

\bibitem {C} A.P. Calder\'on, Ergodic theory and translation-invariant operators. Proc. Nat. Acad. Sci. USA \textbf{59} (1968), 349-353.
 
\bibitem {CZ} A.P. Calder\'on and A. Zygmund, On the existence of certain singular integrals. Acta Math. \textbf{88} (1952), 85-139.

\bibitem {CPPR} M. Caspers, J. Parcet, M. Perrin and \'E. Ricard, Noncommutative de Leeuw theorems. Forum of Mathematics $\Sigma$ \textbf{3} (2015), e21.

\bibitem {CPSZ} M. Caspers, D. Potapov, F. Sukochev and D. Zanin, Weak type commutator and Lipschitz estimates: resolution of the Nazarov-Peller conjecture.
Amer. J. Math. To appear. arXiv: 1506.00778.

\bibitem {CS} M. Caspers and M. de la Salle, Schur and Fourier multipliers of an amenable group acting on non-commutative $L_p$-spaces. Trans. Amer. Math. Soc. \textbf{367} (2015), 6997-7013. 

\bibitem {Chen} J.Chen, Heat kernel on positively curved manifolds and their applications. Ph.D. Thesis, Hangzhou University, 1987. 

\bibitem {CXY} Z. Chen, Q. Xu and Z. Yin, Harmonic analysis on quantum tori. Comm. Math. Phys. \textbf{322} (2013), no. 3, 755-805.

\bibitem {Ch} M.D. Choi, A Schwarz inequality for positive linear maps on $C^*$-algebras. Illinois J. Math. \textbf{18} (1974), 565-574.

%\bibitem {CW2} R. Coifman and G. Weiss, Analyse harmonique non-commutative sur certains espaces homog{\`e}nes. Lecture Notes in Math. \textbf{242}, 1971.

\bibitem {CW} R. Coifman and G. Weiss, Transference methods in analysis. Regional Conference Series in Mathematics \textbf{31}. American Mathematical Society, 1976.

\bibitem {CW3} R. Coifman and G. Weiss, Extensions of Hardy spaces and their use in analysis. Bull. Amer. Math. Soc. \textbf{83} (1977), 569-645.

\bibitem {Co} M. Cotlar, A unified theory of Hilbert transforms and ergodic theorems. Rev. Mat. Cuyana \textbf{1} (1955), 105-167.

\bibitem {Dav} E.B. Davies, Explicit constants for Gaussian upper bounds on heat kernels. Amer. J. Math. \textbf{109} (1987), 319-333.

\bibitem {Dav2} E.B. Davies, Heat Kernels and Spectral Theory. Cambridge. Univ. Press, 1989.

\bibitem {dL}  K. de Leeuw, On $L_p$ multipliers. Ann. of Math. \textbf{81} (1965), 364-379.

\bibitem {DY1} X.T. Duong and L. Yan, New function spaces of BMO type, the John-Nirenberg inequality, interpolation, and applications. Comm. Pure Appl. Math. \textbf{58} (2005), 1375-1420.

\bibitem {DY2} X.T. Duong and L. Yan, Duality of Hardy and BMO spaces associated with operators with heat kernel bounds. J. Amer. Math. Soc. \textbf{18} (2005), 943-973.

\bibitem {ES} M. Enock and J.M. Schwartz, Kac algebras and duality of locally compact groups. With a preface by A. Connes and a postface by A. Ocneanu. Springer, 1992.

\bibitem {ER} E.G. Effros and Z.J. Ruan, Operator Spaces. London Math. Soc. Monogr. \textbf{23}, Oxford University Press, 2000.

\bibitem {MeietalAIM} T. Ferguson, T. Mei and B. Simanek, $H^\infty$ calculus for semigroup generators on BMO. Adv. Math. \textbf{347} (2019), 408-441.

\bibitem{GPJP} A. M. Gonz\'alez-P\'erez, M. Junge and J. Parcet, Singular integral in quantum Euclidean spaces. Mem. Amer. Math. Soc. To appear. arXiv:1705.01081. 

\bibitem{HLP18} H. Ha, G. Lee and R. Ponge. Pseudodifferential calculus on noncommutative tori I. arXiv:1803.03575.

\bibitem {H1} U. Haagerup, Operator valued weights in von Neumann algebras I. J. Funct. Anal. \textbf{32} (1979), 175-206.

\bibitem {H2} U. Haagerup, Operator valued weights in von Neumann algebras II. J. Funct. Anal. \textbf{33} (1979), 339-361.

\bibitem {HLMP} G. Hong, L.D. L\'opez-S\'anchez, J.M. Martell and J. Parcet, Calder\'on-Zygmund operators associated to matrix-valued kernels. Int. Math. Res. Not. 2014 \textbf{14}, 1221-1252.

\bibitem {Ho} L. H\"ormander, Estimates for translation invariant operators in $L^p$ spaces. Acta Math. \textbf{104} (1960), 93-140.
    
%\bibitem {Hy} T. Hyt\"onen, The vector-valued nonhomogeneous Tb theorem. Int. Math. Res. Not. IMRN 2014 \textbf{2}, 451-511. 

%\bibitem {HW1} T. Hyt\"onen and L. Weis, A $T1$ theorem for integral transformations with operator-valued kernel. J. reine angew. Math. \textbf{599} (2006), 155-200.

\bibitem {JLX} M. Junge, C. Le Merdy and Q. Xu, $H^\infty$-functional calculus and square functions on noncommutative $L^p$-spaces. Ast\'erisque \textbf{305}, 2006.

\bibitem {JM}  M. Junge and T. Mei, Noncommutative Riesz transforms--A probabilistic approach. Amer. J. Math. \textbf{132} (2010), 611-680.

\bibitem {JM2} M. Junge and T. Mei, BMO spaces associated with semigroups of operators. Math. Ann. \textbf{352} (2012), 691-743.

\bibitem {JMP4} M. Junge, T. Mei and J. Parcet, Smooth Fourier multipliers on group von Neumann algebras. Geom. Funct. Anal. \textbf{24} (2014), 1913-1980. 

\bibitem {JMP5} M. Junge, T. Mei and J. Parcet, Noncommutative Riesz transforms -- Dimension free bounds and Fourier multipliers. J. Eur. Math. Soc. \textbf{20} (2018), 529-595. 

\bibitem {JX2} M. Junge and Q. Xu, Noncommutative maximal ergodic theorems. J. Amer. Math. Soc. \textbf{20} (2007), 385-439.

\bibitem {KR} R.V. Kadison and J.R. Ringrose, Fundamentals of the Theory of Operator Algebras I and II. Grad. Stud. Math. \textbf{15} and \textbf{16}. American Mathematical Society, 1997.

\bibitem {KW} G. Kuperberg and N. Weaver, A von Neumann algebra approach to quantum metrics. Mem. Amer. Math. Soc. \textbf{215} (2012).

\bibitem {KV1} J. Kustermans and S. Vaes, Locally compact quantum groups, Ann. Sci.  \'Ecole Norm. Sup. \textbf{33} (2000), 837-934.

\bibitem {KV2} J. Kustermans and S. Vaes, Locally compact quantum groups in the von Neumann algebraic setting. 
Math. Scand.  \textbf{92} (2003), no. 1, 68-92. 

\bibitem{LSZ17} G. Levitina, F. Sukochev and D. Zanin. Cwikel estimates revisited. arXiv:1703.04254

\bibitem {Li} J. Li, Gradient estimate for the heat kernel of a complete Riemannian manifold and its applications. J. Funct. Anal. \textbf{97} (1991), 293-310.  

\bibitem{MSX19} E. McDonald, F. Sukochev and X. Xiong. Quantum differentiability on quantum tori. Comm. Math. Phys. 2019.

\bibitem {Me2} T. Mei, Operator valued Hardy spaces. Mem. Amer. Math. Soc. \textbf{188} (2007).

\bibitem {Me3} T. Mei, Tent spaces associated with semigroups of operators. J. Funct. Anal. \textbf{255} (2008), 3356-3406.

\bibitem {MP} T. Mei and J. Parcet, Pseudo-localization of singular integrals and noncommutative Littlewood-Paley inequalities. Int. Math. Res. Not. 2009 \textbf{9}, 1433-1487.

\bibitem {NTV1}  F. Nazarov, S. Treil and A. Volberg,  Weak type estimates and Cotlar inequalities for Calder\'on-Zygmund operators on nonhomogeneous spaces. Internat. Math. Res. Notices 1998, no. 9, 463-487.

\bibitem {NTV2}  F. Nazarov, S. Treil and A. Volberg, The Tb-theorem on non-homogeneous spaces. Acta Math. \textbf{190} (2003), no. 2, 151-239.

\bibitem {NR} S. Neuwirth and E. Ricard, Transfer of Fourier multipliers into Schur multipliers and sumsets in a discrete group. Canad. J. Math. \textbf{63} (2011), 1161-1187.

\bibitem {Pa1} J. Parcet, Pseudo-localization of singular integrals and noncommutative Calder\'on-Zygmund theory. J. Funct. Anal. \textbf{256} (2009), 509-593.

\bibitem {PRdlS} J. Parcet, \'E. Ricard and M. de la salle, Fourier multipliers in $\mathrm{SL}_n(\mathbf{R})$. Preprint. arXiv: 1811.07874.

\bibitem {Ps} W. Paschke, Inner product modules over $B^*$ algebras. Trans. Amer. Math. Soc. \textbf{182} (1973), 443-468.

\bibitem {P3} G. Pisier, Introduction to Operator Space Theory. Cambridge University Press, 2003.

\bibitem {PX2} G. Pisier and Q. Xu, Non-commutative $L_p$-spaces. Handbook of the Geometry of Banach Spaces II. North-Holland (2003), 1459-1517.

\bibitem {R2} \'E. Ricard, $L_p$-multipliers on quantum tori. J. Funct. Anal. \textbf{270} (2016), 4604-4613.

\bibitem {Ri1} M.A. Rieffel, Metrics on states from actions of compact groups. Doc. Math. \textbf{3} (1998), 215-229.

\bibitem {Ri2} M.A. Rieffel, Group $\mathrm{C}^*$-algebra as compact quantum metric spaces. Doc. Math. \textbf{7} (2002), 605-651.

\bibitem {RRT} J.L. Rubio de Francia, F.J. Ruiz and J.L. Torrea, Calder{\'o}n-Zygmund theory for operator-valued kernels. Adv. in Math. \textbf{62} (1986), 7-48.

%\bibitem {St1} E. M. Stein, Singular integrals and differentiability properties of functions. Princeton Mathematical Series, No. \textbf{30} Princeton Univ. Press,  N.J. 1970.

%\bibitem {St2} E.M. Stein, Harmonic Analysis: Real Variable Methods, Orthogonality, and Oscillatory Integrals. Princeton Math. Ser. \textbf{43}. Princeton Univ. Press, NJ, 1993.

\bibitem{SZ18} F. Sukochev and D. Zanin. Connes integration formula for the noncommutative plane. Commun. Math. Phys. \textbf{359} (2018), 449-466.

\bibitem {PS} D. Potapov and F. Sukochev, Operator-Lipschitz functions in Schatten-von Neumann classes. Acta Math. \textbf{207} (2011), 375-389.

\bibitem {Ta0} M. Takesaki, A characterization of group algebras as a converse of Tannaka-Stinespring-Tatsuuma duality theorem. Amer. J. Math. \textbf{91} (1969), 529-564. 

\bibitem {Ta} M. Takesaki, Theory of operator algebras I. Springer-Verlag, New York, 1979.

\bibitem {To} X. Tolsa, BMO, $H^1$, and Calder\'on-Zygmund operators for non doubling measures. Math. Ann. \textbf{319} (2001), 89-149.

\bibitem {T2} X. Tolsa, Littlewood-Paley theory and the $T(1)$ theorem with non-doubling measures. Adv. Math. \textbf{164} (2001), no. 1, 57-116.

\bibitem {VDN} D. V. Voiculescu,  K. J. Dykema and A. Nica, Free random variables.  CRM Monograph Series, \textbf{1}. American Mathematical Society, Providence, RI, 1992. 

\bibitem{XXY18} X. Xiong, Q. Xu and Z. Yin. Sobolev, Besov and Triebel-Lizorkin spaces on quantum tori. Mem. Amer. Math. Soc. \textbf{252} (2018), 1203.
\end{thebibliography}

\vskip-5pt

\enlargethispage{2cm}

\small \hfill \noindent \textbf{Marius Junge} \\
\null \hfill Department of Mathematics
\\ \null \hfill University of Illinois at Urbana-Champaign \\
\null \hfill 1409 W. Green St. Urbana, IL 61891. USA \\
\null \hfill\texttt{junge@math.uiuc.edu}

\vskip5pt

\hfill \noindent \textbf{Tao Mei} \\
\null \hfill Department of Mathematics
\\ \null \hfill Baylor University \\
\null \hfill 1301 S University Parks Dr, Waco, TX 76798, USA \\
\null \hfill\texttt{Tao\_Mei@baylor.edu}

\vskip5pt

\hfill \noindent \textbf{Javier Parcet} \\
\null \hfill Instituto de Ciencias Matem{\'a}ticas \\ \null \hfill Consejo Superior de
Investigaciones Cient{\'\i}ficas \\ \null \hfill C/ Nicol\'as Cabrera 13-15.
28049, Madrid. Spain \\ \null \hfill\texttt{javier.parcet@icmat.es}

\vskip5pt

\hfill \noindent \textbf{Runlian Xia} \\
\null \hfill Instituto de Ciencias Matem{\'a}ticas \\ \null \hfill Consejo Superior de
Investigaciones Cient{\'\i}ficas \\ \null \hfill C/ Nicol\'as Cabrera 13-15.
28049, Madrid. Spain \\ \null \hfill\texttt{runlian.xia@icmat.es}

\end{document}